\documentclass[11pt,a4paper]{amsart}
\usepackage{amsmath,caption,booktabs,lipsum}
\numberwithin{equation}{section}
\allowdisplaybreaks
\usepackage{amsthm}
\usepackage{graphicx}
\usepackage{pdfpages}
\usepackage{amsfonts}
\usepackage{amssymb}
\usepackage{mathrsfs}
\usepackage{tikz}
\usepackage{caption}
\usepackage{environ}
\usepackage{elocalloc}
\usepackage{url}
\usepackage{caption}
\usepackage{graphicx}
\usepackage[backend = bibtex, style=numeric-comp,doi=false,isbn=false,url=false]{biblatex}
\usepackage{CJKutf8}
\usepackage[normalem]{ulem}
\usepackage{xcolor}
\usepackage{etoolbox}
\usepackage{caption}
\usepackage{mathtools}
\usepackage{dcpic}
\usepackage{tikz-cd}
\usepackage[bottom]{footmisc}
\usepackage{hyperref}
\usepackage{enumitem}
\usepackage{stmaryrd}
\usetikzlibrary{positioning}
\usetikzlibrary{arrows,chains,positioning,scopes,quotes}
\usetikzlibrary{decorations.markings}
\DeclareFontFamily{U}{mathx}{}
\DeclareFontShape{U}{mathx}{m}{n}{<-> mathx10}{}
\DeclareSymbolFont{mathx}{U}{mathx}{m}{n}
\DeclareMathAccent{\widehat}{0}{mathx}{"70}
\DeclareMathAccent{\widecheck}{0}{mathx}{"71}
\newcommand{\hgd}{h_{\operatorname{glued}}}
\newcommand{\add}{\operatorname{Ad}}
\newcommand{\sfc}{\si^{d-5}(C)}
\newcommand{\kk}{{\mathbf{K}}}
\newcommand{\supp}{\operatorname{Supp}}
\newcommand{\reg}{{\operatorname{Reg}\,}}
\newcommand{\ai}{\alpha}
\newcommand{\dimh}{\dim_{\mathcal{H}}}

\newcommand{\sing}{{\operatorname{Sing}}\,}
\newcommand{\bff}{\mathbf{F}}
\newcommand{\sus}{{\operatorname{SU}(3)}}
\newcommand{\cgxyc}{\cg_c^{\frac{1}{c}}\xy}
\newcommand{\sut}{{\mathfrak{su}(3)}}

\newcommand{\A}{\mathbf{A}}
\newcommand{\ee}{{{\boldsymbol{\e}}_{\operatorname{Allard}}}}
\newcommand{\es}{\emptyset}
\newcommand{\rr}{{{\boldsymbol{\rho}}_{\operatorname{product}}}}
\newcommand{\be}{\beta}
\newcommand{\Ga}{\Gamma}
\newcommand{\ga}{\gamma}	
\newcommand{\cg}{\mathbf{g}}
\newcommand{\ct}{\mathbf{t}}
\newcommand{\de}{\delta}
\newcommand{\De}{\Delta}
\newcommand{\e}{\epsilon}

\newcommand{\lam}{\lambda}

\newcommand{\om}{\omega}
\newcommand{\fb}{\mathbf{B}}
\newcommand{\Om}{\Omega}
\newcommand{\si}{\sigma}
\newcommand{\Si}{\Sigma}

\newcommand{\hh}{\mathbf{H}}
\newcommand{\rh}{\rho}
\newcommand{\ta}{\theta}
\newcommand{\Ta}{\Theta}

\newcommand{\cms}{\operatorname{comass}}

\newcommand{\vv}{{\mathbf{v}}}
\newcommand{\mf}[1]{\mathfrak{#1}}
\newcommand{\ms}{\mathbf{M}}

\newcommand{\cd}{\cdots}
\newcommand{\T}{\mathbf{T}}
\newcommand{\ug}{U_{\operatorname{gluing}}}
\newcommand{\s}{\subset}
\newcommand{\ees}{{\e_{\operatorname{sum}}}}
\newcommand{\cp}{^\complement}

\newcommand{\cpt}{\mathbb{CP}^2}
\newcommand{\kp}{\ker \ed\pi}
\newcommand{\kn}{\ker^\perp \ed\pi}
\newcommand{\cgxyco}{\cg_c^{\frac{1}{c}-1}\xy}

\newcommand{\la}{\langle}
\newcommand{\ra}{\rangle}
\newcommand{\ov}[1]{\overline{#1}}
\newcommand{\no}[1]{\left\lVert#1\right\rVert}

\DeclarePairedDelimiter{\ri}{\la}{\ra}
\newcommand{\dvol}{\operatorname{dvol}}
\newcommand{\ed}{{\operatorname{d}}}
\newcommand{\pp}{\perp_{|\cdot|}}
\newcommand{\hkp}{h_{\ker \ed\pi}}
\newcommand{\hkn}{h_{\ker^\perp \ed\pi}}
\newcommand{\du}{^\ast}
\newcommand{\pf}{_\ast}

\newcommand{\ka}{\kappa}
\newcommand{\m}{^{-1}}

\newcommand{\cc}{\mathbf{c}}
\newcommand{\cs}{\mathbf{s}}
\newcommand{\w}{\wedge}

\newcommand{\pd}{\partial}

\newcommand{\dd}{{\boldsymbol\de}}
\newcommand{\nn}{{\boldsymbol\nu}}
\newcommand{\na}{\nabla}

\newcommand{\N}{\mathbb{N}}
\newcommand{\R}{\mathbb{R}}
\newcommand{\Z}{\mathbb{Z}}

\newcommand{\C}{\mathbb{C}}

\newcommand{\injrad}{\operatorname{InjRad}}
\newcommand{\focrad}{\operatorname{FocalRad}}
	
\newcommand{\kxp}{{\ker\ed\Pi_{X,L}}}
\newcommand{\kxn}{{\ker^\perp\ed\Pi_{X,L}}}

\newcommand{\cgxy}{\cg_c\bigg(\frac{|z|}{|x|}\bigg)}
\newcommand{\ccxy}{\cc_c\bigg(\frac{|z|}{|x|}\bigg)}
\newcommand{\csxy}{\cs_c\bigg(\frac{|z|}{|x|}\bigg)}
\newcommand{\xy}{\bigg(\frac{|z|}{|x|}\bigg)}

\newcommand{\fd}{\mathbf{d}}

\newcommand{\tr}{\textnormal{tr}}
\newcommand{\id}{\textnormal{id}}

\newcommand{\ssc}{\si^s(C)}

\newcommand{\su}{\mf{su}}

\newcommand{\ad}{\textnormal{ad}}
\newcommand{\po}{\Pi_{\R^3\times\{0\}^3}}\newcommand{\pt}{\Pi_{\{0\}^3\times\R^3}}
\newcommand{\sdc}{\si^{d-c}(C)}

\makeatletter
\def\thm@space@setup{%
	\thm@preskip=0.2cm plus 0cm minus 0cm
	\thm@postskip=\thm@preskip 
}
\makeatother
\theoremstyle{plain}
\newtheorem{thm}{Theorem}
\newtheorem{exam}{Example}[subsection]
\newtheorem{lem}[exam]{Lemma}

\newtheorem{fact}[exam]{Fact}
\newtheorem{claim}[exam]{Claim}
\newtheorem{consq}[exam]{Consequence}
\newtheorem{assump}[exam]{Assumption}

\newtheorem{defn}[exam]{Definition}

\theoremstyle{definition}
\newtheorem{conj}{Conjecture}
\newtheorem{ques}{Question}

\newtheorem{rem}{Remark}[subsection]

\title[Area-minimizing submanifolds are not generically smooth]{Area-minimizing submanifolds are not generically smooth}
\author{Zhenhua Liu\\Princeton University}
\dedicatory{Dedicated to Xunjing Wei}
\begin{document}
	\setlength{\abovedisplayskip}{5pt}
	\setlength{\belowdisplayskip}{5pt}
	\setlength{\abovedisplayshortskip}{5pt}
	\setlength{\belowdisplayshortskip}{5pt}
	\maketitle\vspace{-3em}
	\begin{abstract}
		We prove that area-minimizing submanifolds are not generically smooth, settling a conjecture of White that asks the generic smoothness of area-minimizing submanifolds. We furthermore establish a lower bound on the Hausdorff dimension of the singular sets of area-minimizing submanifolds with respect to  open sets of Riemannian metrics. The lower bound is $\max\{d-5,d-c\},$ where $d$ denotes the dimension of the submanifold and $c$ denotes the codimension. 
	\end{abstract}
	\section{Introduction}
	The general problem of finding area-minimizing representatives in homology classes is solved by Federer and Fleming (\cite{FF}):
	\begin{quote}
		\emph{Every integral homology class on a compact Riemannian manifold admits an area-minimizing representative.}
	\end{quote}
	In other words we can always find a representative of a homology class that has least area among all representatives of the same homology class. The representative in the above result is found in the category of integral currents, which roughly speaking is the closure of the set of algebraic topological polyhedron chains in the Whitney flat topology \cite{FF}.
	
	Calling the area-minimizing representatives submanifolds was justified due to the following theorem:
	\begin{quote}
		\emph{An area-minimizing representative of a homology class is a smooth submanifold with integer multiplicities on connected components outside of a singular set of codimension at least $2$ with respect to the submanifold.}
	\end{quote}
	Federer first proved the above in \cite{HFts} when the submanifold is a hypersurface. Almgren \cite{FA} proved the above statement in general, and De Lellis, Spadaro, Minter and Skorobogatova (\cite{DS1,DS2,DS3,DMS1,DMS2,DMS3}) strengthened the result by showing that the singular set is countably rectifiable. (See also the work of Krummel-Wickramasekera \cite{BK1,BK2,BK3}.) Thus, we will use interchangebly the terms area-minimizing representatives of homology classes and area-minimizing submanifolds.
	
	The above results on the singular sets of area-minimizing submanifolds are sharp, as singular sets do appear in very natural settings \cite[Section 4]{MR0168727}:
	\begin{quote}
		\emph{Complex algebraic subvarieties, including singular subvarieties, are area-minimizing in the Fubini-Study metric on complex projective spaces.}
	\end{quote}
	The presence of singular sets severely restrict the geometric and analytic tools available to obtain meaningful geometric or analytic statements. This raises the natural question of whether one can obtain everywhere smooth area-minimizing submanifolds in generic metrics.
	
	Thom's classical result \cite{RT} shows that there are integral homology classes that topologically cannot be represented by smooth submanifolds. For these homology classes, finding a smooth area-minimizing representative is impossible due to topological obstructions alone. For homology classes with smooth topological representatives, White raised the following question about area-minimizing representatives in \cite[Problem 5.16]{GMT}:
	\begin{ques}\label{cw}
		\emph{"One can also ask whether singularities in
			a homologically area minimizing cycle in a Riemannian manifold disappear after generic perturbations of the metric."}
	\end{ques}
	Area-minimizing $1$-dimensional submanifolds are embedded geodesics \cite{HFts}, so White's conjecture is trivially true in this case. Besides the $1$-dimensional case, there are only two kinds of progress towards the above question. White solved the above conjecture in the case of $2$-dimensional submanifolds of arbitrary codimension \cite{BWgt} based on \cite{SC,DSS1,DSS2,DSS3}:
	\begin{quote}
		\emph{In generic metrics, $2$-dimensional area-minimizing submanifolds of arbitrary codimension are smooth.}
	\end{quote}
	On the other hand, the masterful efforts of Hardt-Simon, Smale, Li-Wang, Chodosh-Mantoulidis-Schulze solved the conjecture in the case of low dimensional hypersurfaces \cite{HS,NS,LW,CMS1,CMS2}:
	\begin{quote}
		\emph{In generic metrics, area-minimizing hypersurfaces of dimension at most $9$ are smooth.}
	\end{quote}
	
	The above two masterful tour de force results indicate that both the dimension and the codimension of the homology class with respect to the ambient manifold play an important role in determining the generic smoothness of area-minimizing representatives. To this end, we adopt the following notation convention.
	\begin{defn}
		Let integer $d$ denote the dimension of an integral homology class, and integer $c$ denote the codimension of the class with respect to the ambient manifold.
	\end{defn}
	For instance, $d=2$ corresponds to $2$-dimensional homology classes, and $c=1$ corresponds to hypersurface homology classes. Following this convention, the dimension and codimension of an area-minimizing representative/submanifold will also be denoted by $d$ and $c$, respectively.
	
	We now state our main results. In sharp contrast to the two masterful results above, the opposite is true in general. 
	\begin{thm}\label{mts}
		Area-minimizing submanifolds are not generically smooth, provided one of the following holds,
\begin{itemize}
	\item either  $d\ge 5$ and $c\ge 3$,
	\item or $3\le c\le d\le 4$.
\end{itemize}
	\end{thm}
	In other words except for the six families of $c=1,$ $c=2,$ $d=1,$ $d=2,$ $d=3<c$, or $d=4<c$, area-minimizing submanifolds are not generically smooth. 
	
	Our results can be stated more precisely as follows.
	\begin{defn}\label{defnm}
		Let $M^{d+c}$ be a compact smooth oriented $(d+c)$-dimensional manifold and $[\Si]\in H_d(M,\Z)$ be a $d$-dimensional integral homology class.
		
		For a smooth Riemannian metric $g$ on $M^{d+c}$, we define $\inf \dimh\sing([\Si],g)$ to be the lowest possible dimension of singular sets of all area-minimizing submanifolds $T$ in $[\Si]$ with respect to the metric $g:$ 
		\begin{align*}
			\inf \dimh\sing([\Si],g)=\inf_{\{T|T\textnormal{ is area-minimizing in }[\Si]\textnormal{ in }g\}}\dim_{\mathcal{H}}\sing T,
		\end{align*} where $\sing T$ denotes the singular set of $T$, $\dimh$ is the Hausdorff dimension and we set
		\begin{align*}
			\dim_{\mathcal{H}}\sing T=-1,
		\end{align*}if $T$ is smooth.
	\end{defn}
	The number $\inf\dim_{\mathcal{H}}\sing([\Si],g)$ is non-negative if and only if every area-minimizing representative of $[\Si]$ is singular in ambient metric $g$.
	
	For homology classes $[\Si]$ that do not have smooth topological representatives, no smooth area-minimizing submanifolds exist, regardless of ambient metrics. For a zero homology class, any area-minimizing representative is necessarily zero as $d$-dimensional integral currents, and consequently never a smooth $d$-dimensional submanifold. Thus, we only have to focus our attention on the non-zero classes that do  have smooth representatives.
	\begin{thm}\label{mtc}
		If a $d$-dimensional non-zero integral homology class $[\Si]$ on $M^{d+c}$ admits a smoothly embedded representative with integer multiplicities on its connected components, then there is a non-empty open subset $\Om_{[\Si]}$ of the space of smooth Riemannian metrics on $M^{d+c}$ such that for all $g\in\Om_{[\Si]}$,
		\begin{align}
			\inf\dimh\sing([\Si],g)\ge \max\{d-5,d-c\},
		\end{align}provided both $d,c\ge 3.$
	\end{thm}
	Here by an embedded representative with integer multiplicities on its connected components, we mean that $$[\Si]=\sum_{j=1}^nk_j[N_j],$$ with each $k_j\in\Z_+$ and $\{N_1,\dots,N_n\}$ a collection of pairwise disjoint smooth submanifolds. The topology we use on the space of metrics is the classical weak $C^\infty$ topology \cite{MH,GG}, which is metrizable \cite{RHif}\cite[page 35]{MH} and is the same as the Whitney $C^\infty$ topology on compact manifolds \cite{MH}.

	In other words there exists an open set $\Om_{[\Si]}$ in the space of smooth Riemannian metrics, such that every area-minimizing integral current in $[\Si]$ with respect to an ambient metric $g\in \Om_{[\Si]}$  has a non-empty singular set of Hausdorff dimension at least $\max\{d-5,d-c\}$, provided both $d\ge 3,c\ge 3.$ The statement is not vacuous only when $d\ge 5,c\ge 3$ or $3\le c\le d\le 4,$ thus implying Theorem \ref{mts}. 
	
	Theorem \ref{mtc} also implies that when $d$ or $(d-c)$ is large, the singular sets of area-minimizing submanifolds are forced to be large with respect to an open set $\Om_{[\Si]}$ in the space of smooth Riemannian metrics as well.
	
	For the six families $c=1,$ $c=2,$ $d=1,$ $d=2,$ $d=3<c,$ or $d=4<c,$ that are not dealt with in our Theorems \ref{mts} and \ref{mtc}, the case of $d=1$ is smooth, as $1$-dimensional area-minimizing submanifolds are embedded geodesics \cite{HFts}. The case of surfaces, i.e., $d=2$, is proven to be generically smooth by White \cite{BWgt}. Since our Theorem \ref{mtc} includes all dimensions $d\ge 3,$ not considering the the value of $c,$ the theorem is sharp in terms of the dimension $d$ alone. The case of $c=1,$ i.e., hypersurfaces, is widely believed to be generically smooth as well, with the best result obtained by Chodosh-Mantoulidis-Schulze \cite{CMS1,CMS2}. For the rest of the cases that are not covered by our theorem, all the known examples are generically smooth as well. Thus, we conjecture that Theorems \ref{mts} give necessary and sufficient conditions on the dimension and codimension pairs for the generic smoothness of area-minimizing submanifolds:
	\begin{conj}
		Area-minimizing submanifolds are smooth in generic metrics if and only if $c=1,$ $c=2,$ $d=1,$ $d=2,$ $d=3<c,$ or $d=4<c.$ 
	\end{conj}
	For a sanity check of the above conjecture, Thom's result \cite{RT} shows that one can find integral homology classes with no smooth representatives if and only if $d\ge 7$ and $c\ge 3.$ Thus, homology classes in the six families $c=1,$ $c=2,$ $d=1,$ $d=2,$ $d=3<c,$ or $d=4<c$ always have smooth topological representatives. Interestingly, by Thom's classical work \cite{RT}, persistent singular sets in topological representatives start appearing from $d=7,$ while in our Theorems \ref{mts} and \ref{mtc}, persistent singular sets in area-minimizing representatives  start appearing from $d=3.$
	
	We also conjecture that the dimension lower bound on the singular sets in Theorem \ref{mtc} is an upper bound in generic metrics:
	\begin{conj}
		In generic metrics, area-minimizing submanifolds have a singular set of dimension at most $\max\{d-5,d-c\},$ provided $d\ge 5,c\ge 3,$ or $3\le c\le d\le 4,$ and this upper bound $\max\{d-5,d-c\}$ can always be achieved in an open subset of the set of generic metrics.  
	\end{conj}
	In other words in generic metrics, $\max\{d-5,d-c\}$ is the maximum dimension of singular sets, and this maximum dimension can be achieved quite often.
	
	Theorems \ref{mts} and \ref{mtc} also hold if we replace integral homology with mod $2$ homology or dropping the orientability of $M$ in Definition \ref{defnm}. We will discuss this in Section \ref{conc}.
	\subsection{Plan of our proof}\label{planpf}
	We now outline the plan of our proof for Theorems \ref{mts} and \ref{mtc}. 
	
	Theorem \ref{mts} is a direct corollary of Theorem \ref{mtc}, by adding to the conclusion homology classes with no smooth representatives at all due to topological reasons. 
	
	From now on, we will focus on proving Theorem \ref{mtc}.
	
	Recall that our goal is to find open subsets of the space of smooth Riemannian metrics in which the singular sets of area-minimizing submanifolds are of dimension at least  $\max\{d-5,d-c\}.$ We adopt the following definition about persistent singular sets.
	\begin{defn}\label{psing}
		We say a subset $\kk$ of the singular set of an area-minimizing submanifold $T$ in a Riemannian metric $g$ is a \textbf{persistent} \textbf{singular} \textbf{subset}, if there exists an open subset $\Om_T$ containing $g$ of the space of smooth Riemannian metrics, so that for all $h\in\Om_T$
		\begin{align*}
			\inf\dimh\sing([T],h)\ge \dimh \kk.
		\end{align*}
	\end{defn}
	Then our goal is to find an area-minimizing submanifold with persistent singular subsets of dimension $\max\{d-5,d-c\}$ in a smooth Riemannian metric.
	
	For simplicity of exposition, assume that $[\Si]$ can be represented by a connected smooth submanifold $N.$ The general case in Theorem \ref{mtc} of $[\Si]$ having a representative with several connected components of integer multiplicities needs a slightly modified argument but all the core ideas remain the same.
	
	The topological representative $N$ has no singular set to start with. Thus, to achieve our goal of having persistent singular sets, we must \textbf{add} \textbf{singular} \textbf{sets} \textbf{manually} to alter the topological representative $N$ while achieving the following four features:
	\begin{enumerate}
		\item  The added singular sets give an altered topological representative that stays in the same homology class $[\Si]$.\label{ft0}
		\item The added singular sets should have dimension $\max \{d-5,d-c\}.$\label{ft1}
		\item We can find a smooth Riemannian metric $h$, in which the altered topological representative is \textbf{area-minimizing}.\label{ft2}	
		\item The added singular sets have \textbf{persistent singular subsets} for the altered representative area-minimizing in $h$.\label{ft3}
	\end{enumerate}
	Features (\ref{ft0}) and (\ref{ft1}) are only about the altered topological representative, while in Features (\ref{ft2}) and (\ref{ft3}) the altered topological representative is area-minimizing in our newly found metric $h.$
	
	Features (\ref{ft0}), (\ref{ft1}) and (\ref{ft2}) are relatively straightforward to achieve. Feature (\ref{ft0}) uses a classical construction due to Smale \cite{NS}, from which Feature (\ref{ft1}) follows naturally. Feature (\ref{ft2}) uses Zhang's work on making topological representatives area-minimizing in \cite{YZa,YZj,YZt}, while Feature (\ref{ft3}) is the hardest part and will take the most efforts.
	
	The core of our paper will center around constructing and proving persistent singular sets, i.e., achieving Feature (\ref{ft3}). We rely on the following two central facts, the proof of which will take up the majority of our manuscript:
	\begin{fact}\label{fctgl}
		Orthogonal \textbf{transverse} intersections in immersed area-minimizing submanifolds are persistent singular sets, provided the dimension and codimension of the immersions with respect to the ambient space are both \textbf{at least} $3.$ 
	\end{fact}
	Frank Morgan hinted the above fact in \cite{FMtd}. In the case of orthogonal transverse intersections of $3$-dimensional subspaces of $\R^6,$ Gary Lawlor proved the above fact for metric perturbations among Euclidean metrics in \cite[Section 6.4]{GL}. Our proof of the above Fact \ref{fctgl} is essentially a vastly upgraded version of Lawlor's proof. 
	
	Roughly speaking, Lawlor found a sufficient condition for whether intersection singular sets are persistent, which only starts to work from dimension and codimension at least $3$. In Euclidean metrics, the sufficient condition reduces to the existence of solutions to an ODE. The key to our proof is finding the right correction factor for the ODE  to work in general smooth Riemannian metrics.
	
	The integer $(d-c)$ is the dimension of the intersection set of $d$-dimensional transverse immersions of codimension $c$ with respect to the ambient space. This explains the factor $(d-c)$ in our lower bound $\max\{d-5,d-c\}$ on the singular sets of area-minimizing submanifolds in Theorem \ref{mtc}.
	
	Another fact we use is the following. Embed $\cpt$ into the unit sphere $S^7$ of $\R^8$ via the Veronese embedding using Laplacian eigenfunctions as coordinates. Define the cone $C(\cpt)$ to be the union of all rays starting from the origin of $\R^8$  and passing through a point in the Veronese embedding of $\cpt.$ The cone $C(\cpt)$ is area-minimizing (Lemma \ref{fctcpt}), and the product $C(\cpt)\times S^{d-5}$ is area-minimizing in the product $\R^8\times S^{d-5}$ with the product metric \cite[Chapter II, Proposition 7.10]{HL}, where $S^{d-5}$ is the standard $(d-5)$-dimensional sphere. Then we have the following.
	\begin{fact}\label{fctcb}
		The singular set $\{0\}\times S^{d-5}$ of the $d$-dimensional area-minimizing submanifold $C(\cpt)\times S^{d-5}$ is a persistent singular set.
	\end{fact}
	To see the ideas behind the above fact, let us think of the simplest case where $d=5,$ i.e., we have two disjoint copies of $\R^8,$ each copy with an area-minimizing $C(\cpt)$ inside. Note that topologically $\cpt$ as a manifold cannot be the boundary of any $5$-dimensional manifold \cite[page 203]{MScs}. Thus any reasonable perturbation of the Veronese embedding of $\cpt$ will not bound a smooth area-minimizing submanifold. The above Fact \ref{fctcb} is heuristically due to this non-bounding property of $\cpt$. The analogues of the above Fact \ref{fctcb} were first noticed by \cite{HP}.
	
	Note that the singular set in Fact \ref{fctcb} has dimension $(d-5)$. This explains the factor $(d-5)$ in our lower bound $\max\{d-5,d-c\}$ on the singular sets of area-minimizing submanifolds in Theorem \ref{mtc}.
	
	With Facts \ref{fctgl} and \ref{fctcb} in hand, our plan is to manually introduce orthogonal transverse immersions or $C(\cpt)\times S^{d-5}$ into our topological representative of $[\Si]$, while achieving Features \ref{ft0}, \ref{ft1} and \ref{ft2}. Then Feature \ref{ft3} is a direct consequence of Facts \ref{fctgl} and \ref{fctcb}, which finishes the proof.
	\subsection{Overview of the paper}
	Our paper will be structured as follows. 
	In Section \ref{bsdefn}, we will give some basic definitions.
	
	In Section \ref{preprep}, we will prepare a topological representative of $[\Si]$. In Section \ref{modsing}, we will alter the topological representative of $[\Si]$ to add our desired singular sets, while establishing Features \ref{ft0} and \ref{ft1} of our plan. In Section \ref{tbssc}, we establish some technical topological lemmas that are needed for Section \ref{zhang}.
	
	In Section \ref{zhang}, we will find a smooth Riemannian metric in which the altered representative produced in Section \ref{modsing} is area-minimizing, thus achieving Feature \ref{ft2}.
	
	In Section \ref{lawlor}, we will prove Fact \ref{fctgl} in the case of our altered topological representative.
	
	In Section \ref{seccpt}, we will prove Fact \ref{fctcb} in the case of our altered topological representative.
	
	In Section \ref{conc}, we will finish our proof by wrapping up the conclusions from the previous sections. Then we will give some concluding remarks.
	\section*{Acknowledgements}
	The author acknowledges the support
	of the NSF through the grant FRG-1854147. The author is deeply grateful to his advisor, Professor Camillo De Lellis, for presenting this problem, directing attention to the reference \cite{GL}, recommending the use of Zhang's gluing constructions and proofreading this manuscript. His unwavering support has been invaluable. 
	Sincere thanks go to Professor Hubert Bray, who has generously shared with the author invaluable advice about virtually every aspect of this world. The author would also like to thank Professor Frank Morgan and Professor Simon Brendle, for their constant encouragements and supports. Another thank you goes to Professor Yang Li and Professor Song Sun, who raised to the author the issue of perturbing away non-isolated singularities. Sincere tribute is paid to Professor Gary Lawlor and Professor Yongsheng Zhang, whose pioneering work on constructing calibrations is the core of this paper. The author would also like to thank the following colleagues for their interests in this work: Professor\textbf{s} 
	Nicolau Aiex, Otis Chodosh, Tam\'{a}s Darvas, Herman Gluck, Or Hershkovits, Robert Kusner, Yi Lai, 
Davi Maximo,	Peter McGrath,  Andr\'{e} Neves, Antoine Song, Daniel Stern, Ao Sun, and Shing-Tung Yau.
	
	Special thanks are extended to the referees, whose careful reading of this manuscript and uncountably many helpful suggestions are what brought this manuscript to the current state. 
	\section{Basic definitions}\label{bsdefn}
	We will give some basic definitions in this section. The contents are more or less standard and the experienced reader can skip this section.
	\subsection{Manifolds}\label{bdmn}
	As we deal with integral homology in this manuscript, unless otherwise stated, every manifold and submanifold we mention in this manuscript is assumed to be orientable.  
	
	In general,  the symbol $M$ will be reserved for a $(d+c)$-dimensional closed ambient manifold and the symbol $[\Si]$ will be reserved for a $d$-dimensional homology class. When we do not mention the ambient manifold explicitly, it is understood that the ambient manifold is $M.$ 
	
	We will often speak of smooth neighborhoods or smooth open sets, by which we mean an open set with smooth boundary.
	
	We will also make frequent use of transversality \cite{MH}. By this we mean that whenever we have a  pair of smooth submanifolds of possibly different dimension, we can assume without loss of generality that the pair is transverse to each other and thus intersects along a submanifold.
	\subsection{Riemannian geometry}
	When a submanifold $N$ is equipped with the induced Riemannian metric, we will use $\fd_N(p,q)$ to denote the Riemannian distance between points $(p,q)$ on $N$. And $\fd_M(p,K)$ will be used to denote the Riemannian distance of $p$ to a closed subset $K$ on $M.$
	
	We will reserve the symbol $\no{\cdot}$ to denote the Riemannian length of vectors, forms, etc. Whenever we use the symbol $|\cdot|,$ it is understood that we are either taking absolute value or taking the square root of the sum of squares of components in coordinates, which thus are different from $\no{\cdot}$ on general manifolds.
		\subsection{Smooth functions with controlled zero set}
	Let $N$ be a smooth orientable Riemannian manifold. Let $K$ be a compact closed subset of $N.$ We need the following lemma.
	\begin{lem}\label{lemzs}
		There exists a smooth function $f$ on $N$ such that the zero set of $f$ is $K, $ i.e., $$f\m(0)=K,$$ and $f$ vanishes to infinite order at $K.$
	\end{lem}
	\begin{proof}This is a well-known fact. We only give a sketch here.
		By Whitney embedding theorem in \cite{CW} and the fact that smooth function restricts to smooth function on submanifolds, we reduce to the case of finding a smooth function vanishing with zero set $K$ and infinite order vanishing at $K$ in Euclidean space. By \cite[3.1.13]{HF}, there exists a good cover of the complement of $K$ with uniform bound on number of intersections. Then for each ball in the covering, we can take a bump function supported on the ball and non-vanishing in the interior of the ball. Now take a sum of these bump functions with small constants in front of each term, so that the sum of $C^k$ norm is absolutely convergent for all $k$. Thus the infinite sum is smooth, and by construction has zero set precisely $K$. After restricting to an embedding of $N, $ we can adjust the $C^k$ norm of the infinite sum by multiplying with small constants.
	\end{proof} 
	\subsection{Unique continuation of minimal submanifolds}
	In this manuscript, we will frequently use the following classical result:
	\begin{lem}\label{suct}(Strong unique continuation) If two properly embedded $d$-dimensional minimal submanifolds $F,G$ in a $(d+c)$-dimensional Riemannian manifold are tangent to each other of infinite order at a point $p\in F\cap G$, then in a neighborhood of $p$ they coincide.
	\end{lem}
	\begin{proof}
		This is a classical result, but the author could not find a reference for general codimension in arbitrary ambient metrics. We will only provide a sketch of the proof. Set up a Fermi coordinate \cite[Chapter II]{AG} chart $(x_1,x_2,\dots,x_d,y_1,\dots,y_c)$ centered at $p$ and adapted to $G,$ i.e., $G$ being the $x_1x_2\dots x_d$-plane in the coordinate chart. By restricting to a small neighborhood of $p,$ we can assume that $F$ can be written as the graph over the $x_1\dots x_d$-plane of a smooth function $\ov{F}:\R^d\to \R^c$. Our assumption of $F$ and $G$ tangent to each other of infinite order at $p$ means that $\ov{F}$ vanishes to infinite order at $\{0\}^d.$
		
		Write $r=\sqrt{\sum_j x_j^2}.$ Direct calculation shows 
		\begin{align}\label{keq}
			|\De_G \ov{F}|=O\left(|\ov{F}|+r|\ed \ov{F}|\right).
		\end{align}Here $\De_G$ is the Laplacian on $G.$ Now a classical argument using the monotonicity formula of frequency functions finishes the proof. For precise reference, we can use \cite[Theorem 1.8, condition (1.4) with $f=r^2$]{JK} to deduce that $\ov{F}$ vanishes around $p$, and we are done. Note that  \cite[Theorem 1.8]{JK} says the required condition is \cite[Condition (1.3)]{JK}. Unfortunately that is a typo in \cite{JK}, as all the proof in \cite{JK} cites \cite[Condition (1.4)]{JK}.
	\end{proof}
	\subsection{Definition of area-minimizing.}
	The right category to discuss area-minimizing representatives is the category of integral currents. We will use Federer's definitive monograph \cite{HF} and Simon's classical lecture notes \cite{LS} as basic references. 
	
	For our purposes, the reader can just regard integral currents as integer coefficient chains of simplicial complexes in algebraic topology. Note that integral currents are allowed to have finite area boundaries.
	
	We say an integral current $T$ is area-minimizing, if $T$ has the least area among all integral currents homologous to $T$:
	\begin{defn}\label{defnam}
		A $d$-dimensional integral current $T$ is area-minimizing if 
		\begin{align*}
			\ms(T)\le \ms(T+\pd V),
		\end{align*}for all  $(d+1)$-dimensional integral currents $V$. 
	\end{defn}
	Here $\ms$ is the mass of the current, which in our context, is just the area of the underlying set with multiplicity included. For instance, $\ms(\pm 2 S^1)=2\ms(S^1)=4\pi$ for the unit circle $S^1$ in $\R^2.$
	\subsection{Definition of singular sets}
	\begin{defn}\label{defnsm}
		We say an integral current $T$ is smooth at a point $p$ in the support of $T$ if there exists an open set $U$ containing $p$ on $M$, such that  $T$ restricted to $U$ equals an integer multiple of a smooth submanifold $N.$ The definitions of regular sets and singular sets of $T$ are as follows:
		\begin{itemize}
			\item The singular set of $T,$ $\sing T,$ is defined as set of points in the support of $T$ where $T$ is not smooth.		
			\item 
			The regular set of $T,$ $\operatorname{Reg}T$, is defined as the set of the points in the support of $T$ where $T$ is smooth. 
		\end{itemize}
	\end{defn} 
	Here support means the underlying set of an integral current. From now on, we will also use the symbol $\supp T$ to mean the support of $T.$
	\subsection{Multiplicity $1$ integral currents}
	Many powerful theorems from geometric measure theory deal with a special class of integral currents called multiplicity $1$ currents. We give a formal definition as follows.
	\begin{defn}\label{defnmo}
		We say an integral current $T$ is of multiplicity $1,$ if $T$ restricted to a neighborhood of each smooth point $p$ equals a smooth submanifold $N$ and the regular set of $T$ is dense in $\supp T.$
	\end{defn}
	For example, a figure $8$ in $\R^2$ is a multiplicity $1$ integral current.
	\subsection{Connected sum of multiplicity $1$ integral currents}
	Many times in this manuscript, we need to use  a connected sum to connect a representative of $[\Si]$ with multiple components. We give a formal construction as follows.
	
\begin{assump}\label{assumpcs}
Assume that $T$ and $V$ are two $d$-dimensional multiplicity $1$ integral currents of dimension $d$ and codimension ${c}$ \textbf{both at least} $\mathbf{2}$, and both $\reg T$ and $\reg V$ are  non-empty and not equal to each other.
\end{assump}		
	Our goal is to construct a connected sum $T\# V$.
	\begin{lem}\label{lemcs}
Under Assumptions \ref{assumpcs}, there is a multiplicity $1$ integral current $T\# V$, called the connected sum of $T$ and $V,$ so that 
		\begin{itemize}
			\item $\sing (T\# V)=\sing (T+V)$.
			\item Restricted to a neighborhood of $\sing (T\# V),$ the integral current $T\#V$ equals $T+V$.
			\item If $\operatorname{Reg}T$ and $\operatorname{Reg}V$ are both connected, then $\operatorname{Reg}T\# V$ is also connected.
			\item $[T\# V]=[T]+[V]$ as homology classes.
		\end{itemize}
	\end{lem}
	For instance, the connected sum of two disjoint embedded tori in $\R^4$ is a genus $2$ surface. Here the symbol $+$ in $T+V$ means the sum as integral currents. 
	\begin{proof}
		In case $T$ and $V$ are disjoint submanifolds, the construction is classical and can be found in \cite{CL}.

		In our case, $T$ and $V$ are general integral currents that  may intersect each other and may have large singular sets. However, the idea is similar. We need to remove a $d$-dimensional ball from $T$ and $V,$ respectively, then use a tube diffeomorphic to $S^{d-1}\times[0,1]$ to connect $T$ to $V.$ 	The precise construction is as follows.
		
		Take a smooth point $p$ of $T$, and a smooth point $q$ of $V,$ with $p\not=q.$ Without loss of generality, we can assume that there is an open set $U$ containing both $p$ and $q$, such that we have a coordinate $$(x_1,\dots,x_d,y_1,\dots,y_c),$$
		on $U$, with
		\begin{align*}
			p=(\overbrace{0,\dots,0}^{d},\overbrace{0,0,\dots,0}^{c}),q=(\overbrace{0,\dots,0}^{d},1,\overbrace{0,\dots,0}^{c-1})
		\end{align*}Moreover,  the integral current $T$ restricted to $U$ is parameterized by
		\begin{align*}
			(x_1,\dots,x_d,0,0,\dots,0),	\end{align*}with each $x_j\in\R,$ and the integral current $V$ restricted to $U$ is parameterized by
		\begin{align}\label{xdminus}
			(-x_1,x_2,\dots,x_d,1,0,\dots,0),
		\end{align}with each $x_j\in\R.$ The minus sign is to give opposite orientations in order to carry out the connected sum.
		
To construct the desired coordinate chart $U$ as above, we need to use transversality. Though in general $T$ and $V$ per se may have complicated singular sets, by Definition \ref{defnmo} a small perturbation of a smooth curve $\ga$ from $p$ to $q$ will only meet $T$ and $V$ on their regular sets. By our assumption of codimension $c\ge 2,$ we have $d+1<d+c=\dim M,$ so transversality means that $\ga$ intersect $T$ and $V$ only at end points, and $\ga$ is not tangent to $\reg T$ and $\reg V.$ Now take an arbitrary Riemannian metric on $M,$ and extend the unit speed tangent vectors to $\ga$ to a smooth vector field $X_\ga$ on $M.$ Then apply the flow generated by $X_\ga$ we can bring $p$ and $q$ to the same coordinate neighborhood. It is then straightforward to verify that we can obtain the above coordinate chart $U$ as desired.
		
		Now remove the $d$-dimensional ball centered at $p$ $$B_1^d(0)\times\{0\}\times\{0\}^{c-1},$$ from $T$ and remove the $d$-dimensional ball centered at $q$ $$B_1^d(0)\times\{1\}\times\{0\}^{c-1},$$ from $V$, and  glue the tube
		\begin{align*}
			S_1^{d-1}(0)\times [0,1]\times \{0\}^{c-1},
		\end{align*}smoothly onto $T+V$ with the two $d$-dimensional balls removed,
		we obtain $T\# V.$ Here by gluing smoothly we mean smoothing the corner introduced by the boundary of the tube. Since our constructions only alter $\reg T$ and $\reg V,$ we deduce that the first two bullets in Lemma \ref{lemcs} hold. Since $T$ and $V$ are of dimension at least $2$, we deduce the third bullet in Lemma \ref{lemcs}. For the last bullet in Lemma \ref{lemcs}, since
		\begin{align*}
			&	S_1^{d-1}(0)\times [0,1]\times \{0\}^{c-1}\\=&\pd\bigg(B_1^d(0)\times [0,1]\times\{0\}^{c-1}\bigg)	+B_1^d(0)\times\{0\}\times\{0\}^{c-1}-B_1^d(0)\times\{1\}\times\{0\}^{c-1},
		\end{align*}we deduce that $[T\# V]=[T]+[V]$. The minus sign in $-B_1^d(0)\times\{1\}\times\{0\}^{c-1}$ explains the minus sign in parameterization (\ref{xdminus}) of $V$.
	\end{proof}
	\subsection{Calibrations}
	The central object of our manuscript is area-minimizing integral currents. Thus, we must have a way to prove that an integral current is area-minimizing. The way we prove area-minimizing is via calibrations. 
	The notion of calibrations started with the classical paper of Harvey-Lawson \cite{HL}.
	
	Recall that the comass of a  $d$-dimensional differential form  $\phi$ is the maximum of $\phi$ evaluated on simple unit $d$-vectors in the tangent space to $M$ among all  points \cite[Section 1.8]{HF}. In other words
	\begin{align*}
		\cms_g\phi=\max_{p\in M}\max_{\substack{P\s T_pM\\ \dim P=d}}\phi(P).
	\end{align*}
	Here we use the symbol $P$ to denote both a $d$-dimensional plane $P$ and the unit simple $d$-vector representing $P.$
	
	For instance, if $\phi$ is a simple form, then the comass of $\phi$ is equal to the maximum Riemannian length of $\phi$ over all points. 
	
	For the reader's convenience, we record some facts about comass.
	\begin{fact}\label{cmsvec}
		Let $g,h$ be positive definite quadratic forms on a finite dimensional vector space, with $\psi$ a constant differential form.\begin{enumerate}
			\item $\cms_h\psi\le\cms_g\psi$ if $h\ge g$ as quadratic forms.\label{cms1}
			\item $\cms_{\lam^2 h}\psi=|\lam|^{-\deg\psi}\cms_h\psi$.\label{cms2}
			\item $\cms_g \psi=\cms_g\ast\psi$ and $\psi(P)=\cms_g(\psi)\iff \ast\psi(\ast P)=\cms_g(\ast\psi).$\label{cms3}
		\end{enumerate}
	\end{fact}Here the parameter $\lam\in \R^+$, the symbol
	$\cms_g\psi$ means the comass of $\psi$ with respect to $g$ and the symbol $\ast$ is the Hodge star. All the bullets above are direct corollaries of definition of comass. The first two bullets are proved in \cite[Lemma 2.1.9, 2.1.20]{YZt} and the last bullet is proved in \cite[Remark at the bottom of p.79]{HL}. One can easily adapt Fact \ref{cmsvec} to general Riemannian manifolds.
	
	Now we are ready to give the definition of calibrations.
	\begin{defn}(Definition of calibrations)
		\begin{itemize}
			\item 	We say a closed $d$-dimensional Hausdorff measurable $d$-form $\phi$ on a (possibly open) ambient manifold is a calibration if its comass is at most $1.$ 
			\item 
			We say a $d$-dimensional integral current $T$ is calibrated by $\phi,$ if $d$-dimensional Hausdorff measure almost everywhere $\phi$ restricted to the tangent space of $T$ equals the volume form of $T.$
		\end{itemize}
	\end{defn}
	In general, if the calibration form $\phi$ has sufficient regularity, then any integral current $T$ calibrated by $\phi$ will be area-minimizing. The regularity assumption is as follows.
	\begin{defn}\label{defnlip}(Regularity for calibrations to guarantee area-minimizing)
		\begin{itemize}
			\item	We say a degree $d$ closed measurable form $\phi$ has Lipschitz anti-derivatives, if for every point $p$ in the domain of $\phi,$ there exists a neighborhood $U$ of $p$ such that we have $\phi=d\omega$ in $U$ for a  Lipschitz $(d-1)$-form $\omega$ in $U$. We call $\omega$ a Lipschitz anti-derivative of $\phi$ at $p.$ 
			\item We say a $d$-form $\phi$ with Lipschitz anti-derivative is almost continuous on an integral current $T,$ if for every point $p$ in the domain of $\phi,$ there is a Lipschitz anti-derivative $\om$ of $\phi$ at $p,$ such that $\om$ is continuously differentiable on a subset of $\supp T\cap \operatorname{domain} \om$ of full $d$-dimensional Hausdorff measure. 
		\end{itemize}
	\end{defn}
	\begin{lem}\label{fcal}
		If a $d$-dimensional integral current $T$ on a compact ambient manifold is calibrated by a degree $d$ calibration form $\phi$ with Lipschitz anti-derivative, then $T$ is area-minimizing (Definition \ref{defnam}) provided that $\phi$ is almost continuous on $T.$ Furthermore, any integral current homologous to $T$ and having the same area as $T$ must also be calibrated by $\phi.$
		\end{lem}
	\begin{proof}
		This proof is a straightforward adaptation of \cite[Theorem A8]{GL} from Euclidean space to general Riemannian manifolds. We will only give a sketch of the argument.
		
		If $\phi$ is smooth, then this is the classical fundamental theorem of calibrated geometry proved by Harvey-Lawson \cite[Chapter II, Theorem 4.2]{HL}. The proof goes as follows. The condition of comass at most $1$ implies that $R(\phi)\le\ms(R)$ for all $d$-dimensional integral currents $R.$ Here $R(\phi)$ means integrating $\phi$ on the integral current $R.$ The condition of $T$ calibrated by $\phi$ means that $T(\phi)=\ms(T).$ Stokes Theorem and closedness of $\phi$ imply that $T(\phi)=(T+\pd V)(\phi)$ for any $(d+1)$-dimensional integral current $V.$ Thus, for an integral current $T$ calibrated by $\phi,$ we deduce that\begin{align}\label{ftmc}
			\ms(T)=T(\phi)=(T+\pd V)(\phi)\le\ms(T+\pd V).
		\end{align}In other words $T$ is area-minimizing. Furthermore, equality in (\ref{ftmc}) holds if and only if $T+\pd V$ is also calibrated by $\phi.$
		
		When $\phi$ is not smooth, we cannot integrate $\phi$ directly on currents. However, when the ambient space is the Euclidean space we can integrate mollifications $\phi_\e$ of $\phi$ on $T.$ The inequality (\ref{ftmc}) becomes
		\begin{align}\label{ftmce}
			\ms(T)=(1+O(\e))T(\phi_\e)=(1+O(\e))(T+\pd V)(\phi_\e)\le (1+O(\e))\ms(T+\pd V).
		\end{align}Taking $\e\to 0$ shows that $\ms(T)\le\ms(T+\pd V)$. This is precisely Lawlor's proof in \cite[Theorem A8]{GL}. Again as $\e\to0,$ we can deduce that $\ms(T)=\ms(T+\pd V) $ if and only if $T+\pd V$ is also calibrated by $\phi.$
		
		For our proof, we need to work with general manifolds, which does not have the structure of Euclidean spaces. To overcome this, we need to embed our ambient space into Euclidean space to obtain mollifications of $\phi$. Precisely speaking, isometrically embed the compact ambient manifold into Euclidean space using Nash embedding. Use $\pi$ to denote the nearest distance projection onto the embedded image of the ambient manifold. Then $\pi\du(\phi)$ is a closed form defined in a tubular neighborhood of the embedded image of the ambient manifold. It is straightforward to check that mollifications of $\pi\du(\phi)$ can serve the same role as mollifications of $\phi$ in estimates (\ref{ftmce}). Thus, we can still establish an estimate of (\ref{ftmce}). The claim about $\ms(T)=\ms(T+\pd V)$ implying $T+\pd V$ calibrated by $\phi$ follows directly from this. \end{proof}
	\section{Preparing a representative of $[\Si]$}\label{preprep}Let us first collect the several assumptions in Theorem \ref{mtc}.
\begin{assump}\label{assumpdc3}
	Based on Definition \ref{defnm}, assume both $d,c\ge 3$ and that the non-zero integral homology class $[\Si]$ has a smooth representative with integer multiplicities on connected components.
\end{assump}
In other words $[\Si]$ has an integral current representative as follows:
	\begin{align}\label{sigs}
		\Si=\overbrace{N_1+\dots+N_1}^{k_1}+\overbrace{N_2+\dots+N_2}^{k_2}+\dots+\overbrace{N_n+\dots+N_n}^{k_n},
	\end{align}where $\{N_1,\dots, N_n\}$ is a collection of pairwise disjoint $d$-dimensional smooth submanifolds and $k_1,\dots,k_n\in \Z_+.$
	
	However, it will turn out later that the multiplicities $k_1,\dots,k_n$ will prevent us from using many powerful theorems from geometric measure theory. Thus, we need a representative that is of multiplicity $1$ (Definition \ref{defnmo}). 
	
	Also, it will become clear later that to make topological representatives area-minimizing, we need the topological representatives to be connected on the regular set.
	
	In this section we will solve these two issues.
	\subsection{Reduction to a multiplicity $1$ current}
	To destroy the multiplicities of $\Si$ in (\ref{sigs}), first put all the summands in (\ref{sigs}) into general position:
	\begin{align}\label{sigdis}
		\Si=(N_1^1+\dots+N_1^{k_1})+\dots+(N_n^1+\dots +N_n^{k_n}).
	\end{align}
	Here each $N_j^k$ is a perturbation of $N_j$ by a generic ambient diffeomorphism into general position. In other words $\Si$ is an immersion with normal crossings. 
	\begin{defn}\label{defnimnc}
		An immersion with normal crossings is defined to be an immersion that is a stable mapping \cite[Definition 1.1]{GG}.
	\end{defn} A mapping is called stable if in the space of mappings there exists an open subset containing the mapping so that all of the mappings in this open subset are equivalent up to an ambient diffeomorphism of the target manifold. Here the topology on the space of mappings is the Whitney $C^\infty$ topology \cite{GG}, which is the same as the metrizable weak topology on compact manifolds \cite{MH}.
	
	For example, an immersion with transverse intersections is an immersion with normal crossings. Also, the union of $xy$, $xz$, and $yz$-planes in $\R^3$ with coordinate $(x,y,z)$ is an immersion with normal crossings. A detailed exposition of immersions with normal crossings is in the classical work \cite[Section III.3]{GG}. However, the reader can just regard this notion as a black box, as we only need three consequences of the property of being an immersion with normal crossings.
	\begin{fact}\label{inct}
		Immersed submanifolds $N$ with normal crossings have tubular neighborhoods that deformation retract onto $N$.
	\end{fact}
	\begin{fact}\label{incm}
		Immersed submanifolds $N$ with normal crossings are multiplicity $1$ currents.
	\end{fact}
	\begin{fact}\label{imcal}
		There exists a smooth metric $h_N$ and a smooth calibration form $\phi_N$ on a neighborhood $U(N\cap N)$ of the self-intersection set $N\cap N$ of an immersed submanifold $N$ with normal crossings  such that $\phi_N$ calibrates $N$ restricted to $U(N\cap N),$ provided the dimension and codimension of $N$ are both at least $3.$
	\end{fact}
	We will prove a much stronger version of Fact \ref{inct} later in the Section \ref{tbssc}. Fact \ref{incm} follows directly from definition of immersions with normal crossings in \cite[Section III.3]{GG}. 
	
	Fact \ref{imcal} is a byproduct of the proof of \cite[Theorem 3.20]{YZa}. \cite[Theorem 3.20]{YZa} states that for embedded submanifolds in general position and linearly independent in real homology, we can find a smooth differential form in a Riemannian metric that calibrates them all at the same time. However, in the proof of \cite[Theorem 3.20]{YZa}, a calibration form is first constructed around the intersection set of the embedded submanifolds, and does not use the linear independence condition, which is only used to extend the calibration to be defined everywhere. Thus, the proof of \cite[Theorem 3.20]{YZa} applies  to prove Fact \ref{imcal}.
	\subsection{Connecting the regular set}
	To make the regular set of $\Si$ in (\ref{sigdis}) connected, 
	we need to apply  the connected sum of multiplicity $1$ currents defined in Lemma \ref{lemcs} to $\Si$. In other words we obtain a representative $N$ of $\Si$, where
	\begin{align*}
		N=(N_1^1\# \cdots\# N_1^{k_1})\#\cdots\#(N_n^1\# \cdots\# N_n^{k_n}).
	\end{align*}
	By Lemma \ref{lemcs}, $N$ is an immersed connected submanifold with normal crossings, representing the homology class $[\Si]$.
	
	To sum it up, the assumption on $[\Si]$ in the statement of Theorem \ref{mtc} can be reduced to the following assumption.
	\begin{assump}\label{simpa}
		The integral homology class $[\Si]$ is represented by a smoothly immersed connected submanifold $N$ with normal crossings (Definition \ref{defnimnc}).
	\end{assump}\raggedbottom
	\section{Construction of singular sets }\label{modsing}
With Assumption \ref{simpa} in mind, as outlined in the plan of our proof (Section \ref{planpf}),  we need to add singular sets manually to the topological representative $N$ of $[\Si]$. The way to do this is to construct a homologically trivial cycle $\si_s(C)$ on $M$  in Section \ref{secssc} with our desired singular set and using Lemma \ref{lemcs} to do a connected sum $N\#\si_s(C).$ This section will be devoted to the construction of the homologically trivial cycle $\si_s(C).$	
	\subsection{The $s$-th spherical link of a cone $C$}\label{secssc}
	Let $s$ denote a nonnegative integer at most $(d-3)$ that we will set later to be the dimension of the singular set we want to add. Let $V$ be a $(d-s-1)$-dimensional not necessarily connected closed smooth submanifold of the unit sphere $S^{d-s+c-1}$of $\R^{d-s+c}$. 
	\begin{defn}\label{defncv}
		The cone $C(V)$  over $V$ in $\R^{d-s+c}$ is defined to be the union of all rays from the origin passing through $V.$
	\end{defn}
	Here by a ray passing through a point $p$ we mean the half line parameterized by $tp$, with $t\in \R_{\ge 0}.$
	
	It is straightforward to verify that $C(V)$ is a boundaryless multiplicity $1$ integral current. When $V$ is clear from the context, we will drop the symbol $V$ and write $C$ only. In case $C$ is not a flat subspace of $\R^{d-s+c}$, it is straightforward to see that
	\begin{align*}
		\sing C=\{0\}^{d-s+c}.
	\end{align*}
	From now on, we assume that $C$ is not a flat subspace of $\R^{d-s+c}$. 
	\begin{defn}
		Define the $s$-th spherical link of $C$ to be				\begin{align}\label{prodl}
			\si^s(C)=(C\times \R^{s+1})\cap{ S^{d+c}}.
		\end{align}
		The symbol $S^{d+c}$ denotes the $(d+c)$-dimensional unit sphere in $\R^{d+c+1}$ and $C\times \R^{s+1}$ is embedded naturally in the product $\R^{d-s+c}\times \R^{s+1}\cong \R^{d+c+1}.$ 		
	\end{defn}
	It is straightforward to verify that $\si^s(C)$ is a $d$-dimensional multiplicity $1$ integral current on the $(d+c)$-dimensional sphere $S^{d+c}$.
	\begin{figure}
		\centering
		\def\svgwidth{0.9\paperwidth}
\begingroup%
  \makeatletter%
  \providecommand\color[2][]{%
    \errmessage{(Inkscape) Color is used for the text in Inkscape, but the package 'color.sty' is not loaded}%
    \renewcommand\color[2][]{}%
  }%
  \providecommand\transparent[1]{%
    \errmessage{(Inkscape) Transparency is used (non-zero) for the text in Inkscape, but the package 'transparent.sty' is not loaded}%
    \renewcommand\transparent[1]{}%
  }%
  \providecommand\rotatebox[2]{#2}%
  \newcommand*\fsize{\dimexpr\f@size pt\relax}%
  \newcommand*\lineheight[1]{\fontsize{\fsize}{#1\fsize}\selectfont}%
  \ifx\svgwidth\undefined%
    \setlength{\unitlength}{2160bp}%
    \ifx\svgscale\undefined%
      \relax%
    \else%
      \setlength{\unitlength}{\unitlength * \real{\svgscale}}%
    \fi%
  \else%
    \setlength{\unitlength}{\svgwidth}%
  \fi%
  \global\let\svgwidth\undefined%
  \global\let\svgscale\undefined%
  \makeatother%
  \begin{picture}(1,0.5)%
    \lineheight{1}%
    \setlength\tabcolsep{0pt}%
    \put(0,0){\includegraphics[width=\unitlength,page=1]{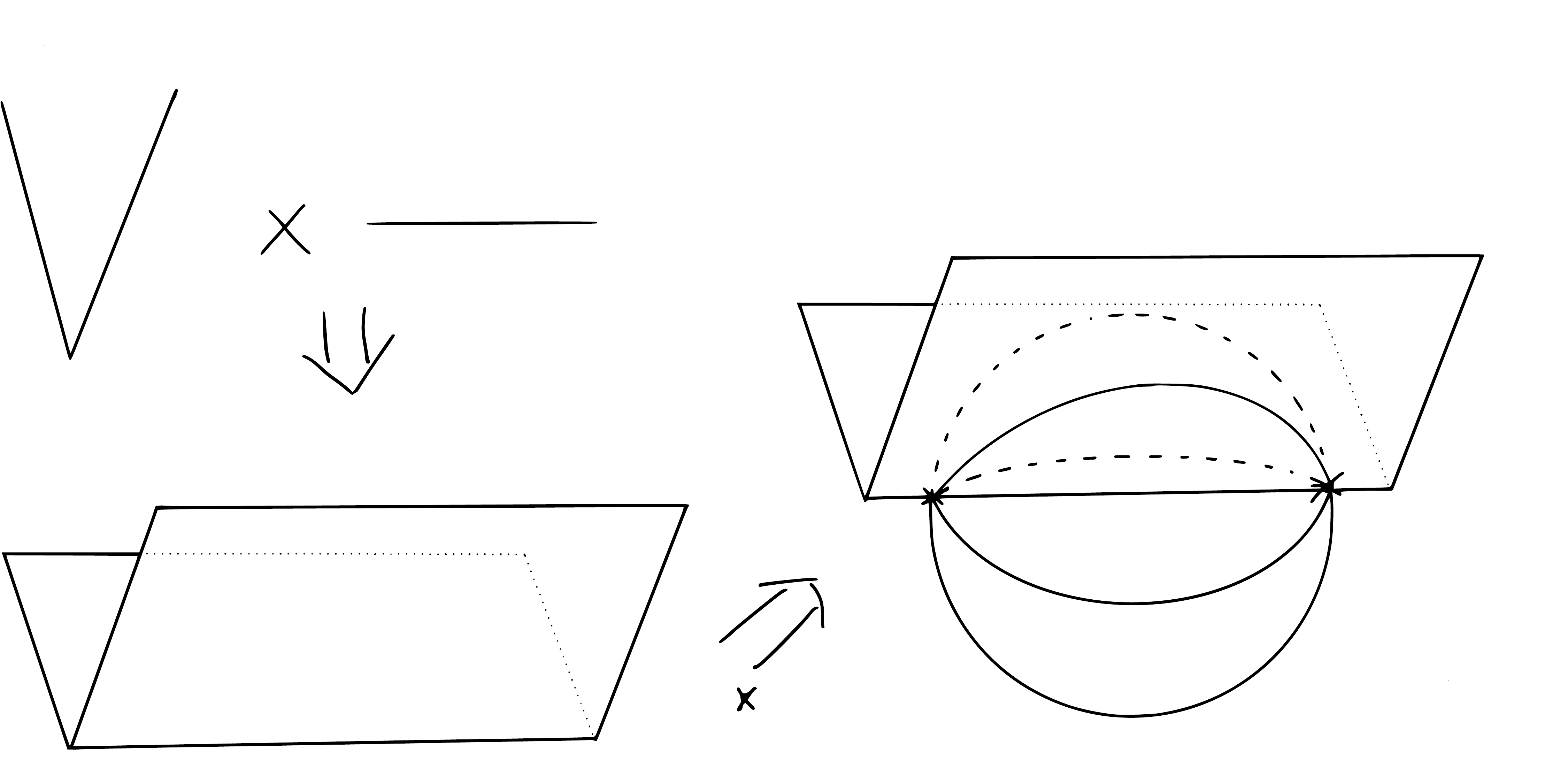}}%
    \put(0.00496863,0.45693687){\color[rgb]{0,0,0}\makebox(0,0)[lt]{\lineheight{1.25}\smash{\begin{tabular}[t]{l}$C^{d-s}\s \R^{d-s+c}$\end{tabular}}}}%
    \put(0.27824023,0.38284622){\color[rgb]{0,0,0}\makebox(0,0)[lt]{\lineheight{1.25}\smash{\begin{tabular}[t]{l}$\R^{s+1}$\end{tabular}}}}%
    \put(0.10669392,0.20236793){\color[rgb]{0,0,0}\makebox(0,0)[lt]{\lineheight{1.25}\smash{\begin{tabular}[t]{l}$C^{d-s}\times \R^{s+1}\s \R^{d+c+1}$\end{tabular}}}}%
    \put(0.42363669,0.02495004){\color[rgb]{0,0,0}\makebox(0,0)[lt]{\lineheight{1.25}\smash{\begin{tabular}[t]{l}$\sing \ssc =\{0\}^{d-s+c}\times S^s$\end{tabular}}}}%
    \put(0.56380286,0.3628129){\color[rgb]{0,0,0}\makebox(0,0)[lt]{\lineheight{1.25}\smash{\begin{tabular}[t]{l}$\ssc=(C^{d-s}\times \R^{s+1})\cap S^{d+c}$\end{tabular}}}}%
    \put(0.67572663,0.08771092){\color[rgb]{0,0,0}\makebox(0,0)[lt]{\lineheight{1.25}\smash{\begin{tabular}[t]{l}$S^{d+c}\s \R^{d+c+1}$\end{tabular}}}}%
  \end{picture}%
\endgroup%

				\caption{Construction of $\ssc$}
		\label{fig:ssc}
	\end{figure}
		\begin{lem}\label{singssc}
		There exists a neighborhood $U(\{0\}^{d-s+c}\times S^s)$ of $\{0\}^{d-s+c}\times S^s$ in $S^{d+c}$ diffeomorphic to the product of the unit ball $B_1^{d-s+c}$ in $\R^{d-s+c}$ with  $S^s$ by a diffeomorphism $\Ga$, i.e.,
		\begin{align*}
			U(\{0\}^{d-s+c}\times S^s)\overset{\Ga}{\cong}B_1^{d-s+c}\times S^s,
		\end{align*}such that $\ssc$ restricted to $U(\{0\}^{d-s+c}\times S^s)$ equals 
		\begin{align*}
			\big(C|_{B_1^{d-s+c}}\big)\times S^s,
		\end{align*}via $\Ga.$
		Moreover, the singular set of $\ssc$ is $\{0\}^{d-s+c}\times S^s$.
	\end{lem}
	Here $S^s$ denotes the unit sphere in $\R^{s+1} $ and we regard $S^0$ as two disjoint points. We use the symbol $T|_K$ to denote the restriction of an integral current to a set $K.$ 
	\begin{proof}
		Let us construct explicitly the diffeomorphism $\Ga$ and its inverse $\Ga\m.$
		
		Adopt a coordinate system
		\begin{align*}
			x=	(x_1,\dots,x_{d-s+c}),
		\end{align*}on $\R^{d-s+c}$, and a coordinate system
		\begin{align*}
			y=(y_1,\dots,y_{s+1}),
		\end{align*}on $\R^{s+1}.$ The product structure $\R^{d+c+1}\cong \R^{d-s+c}\times \R^{s+1}$ gives a coordinate system
$		(x,y)
$		on $\R^{d+c+1}.$
		
		Define a map $\Ga:\R^{d+c+1}\to \R^{d+c+1}$ by 
		\begin{align*}
			\Ga(x,y)=\bigg(\frac{x}{|y|},\frac{y}{|y|}\bigg).
		\end{align*}
		\begin{claim}\label{gabi}
			We claim that $\Ga$ is a smooth diffeomorphism from $S^{d+c}\setminus\{|y|=0\}$ onto $\R^{d-s+c}\times S^s.$
		\end{claim}
		To prove the above claim, define a smooth map $\Ga':\R^{d+c+1}\to\R^{d+c+1},$
		\begin{align*}
			\Ga'(x,y)=\bigg(\frac{x}{\sqrt{1+|x|^2}},\frac{y}{\sqrt{1+|x|^2}}\bigg).
		\end{align*}
		The map	$\Ga'$ is smooth everywhere and we wish to show that
		\begin{align*}
			\left(\Ga|_{S^{d+c}\setminus\{|y|=0\}}\right)\m=\Ga'|_{\R^{d-s+c}\times S^s}.
		\end{align*}Direct calculation shows that the map $\Ga'$ restricted to $\R^{d-s+c}\times S^s,$ i.e., where $|y|=1,$ satisfies
		\begin{align*}
			y\big(\Ga'|_{\R^{d-s+c}\times S^s}(x,y)\big)\not=&\,0,\\\big|\Ga'|_{\R^{d-s+c}\times S^s}(x,y)\big|=&\,1,\\\Ga|_{S^{d+c}\setminus\{|y|=0\}}\circ\Ga'|_{\R^{d-s+c}\times S^s}(x,y)=&\,(x,y).
		\end{align*}
		Here $y\big(\Ga'|_{\R^{d-s+c}\times S^s}(x,y)\big)$ means the $y$ coordinate of $\Ga'|_{\R^{d-s+c}\times S^s}(x,y).$ In other words $\Ga'|_{\R^{d-s+c}\times S^s}$ maps $\R^{d-s+c}\times S^s$ to $S^{d+c}\setminus\{|y|=0\}$
		and $\Ga'|_{\R^{d-s+c}\times S^s}$ is a right inverse of $\Ga|_{S^{d+c}\setminus\{|y|=0\}}.$ Thus $\Ga'|_{\R^{d-s+c}\times S^s}$ must  be injective and $\Ga|_{S^{d+c}\setminus\{|y|=0\}}$ must be surjective onto $\R^{d-s+c}\times S^s$. 
		
		On the other hand, direct calculations also show that $\Ga$ restricted to $S^{d+c}\setminus\{|y|=0\},$ i.e., where
		$|x|^2+|y|^2=1,|y|\not=0,$ is a smooth map that satisfies
		\begin{align*}
			\Big|y\big(\Ga|_{S^{d+c}\setminus\{|y|=0\}}(x,y)\big)\Big|=&\,1,\\\Ga'|_{\R^{d-s+c}\times S^s}\circ\Ga|_{S^{d+c}\setminus\{|y|=0\}}(x,y)=&\,(x,y).
		\end{align*}
		In other words $\Ga$ maps $S^{d+c}\setminus\{|y|=0\}$ to $\R^{d-s+c}\times S^s$ and $\Ga'|_{\R^{d-s+c}\times S^s}$ is also a left inverse of $\Ga|_{S^{d+c}\setminus\{|y|=0\}}.$ Thus $\Ga'|_{\R^{d-s+c}\times S^s}$ must be surjective onto $S^{d+c}\setminus\{|y|=0\}$ and $\Ga|_{S^{d+c}\setminus\{|y|=0\}}$ must be injective.
		
		To sum it up, both $\Ga|_{S^{d+c}\setminus\{|y|=0\}}$ and $\Ga'|_{\R^{d-s+c}\times S^s}$ are smooth bijective maps and they are inverse maps of each other. This finishes the proof of Claim \ref{gabi}.
		
		Now we are ready to find the neighborhood $U(\{0\}^{d-s+c}\times S^s)$ in the statement of Lemma \ref{singssc}. Since $\{0\}^{d-s+c}\times S^s$ is an orbit of the standard Lie group action $$\operatorname{O}(d-s+c)\times\operatorname{O}(s+1),$$ on $\R^{d-s+c}\times \R^{s+1},$ we deduce that the tubular neighborhoods $U_r$ of radius $r$ around $\{0\}^{d-s+c}\times S^s$ in $S^{d+c}$ are also invariant under the action of $O(d-s+c)\times O(s+1)$ as well. Using this, we can calculate that
		\begin{align*}
			U_r=\{(x,y)\,|\,|x|^2+|y|^2=1,|x|\le\sin r,(x,y)\in \R^{d-s+c}\times \R^{s+1}\}.
		\end{align*}
		The same calculations as in the proof of Claim \ref{gabi} show that $U_{\frac{\pi}{4}}$ is mapped by $\Ga|_{S^{d+c}\setminus\{|y|=0\}}$ bijectively to
		$		B_{1}^{d-s+c}\times S^s.$ 
		\begin{defn}
			Define $$U(\{0\}^{d-s+c}\times S^s)=U_{\frac{\pi}{4}}.$$
		\end{defn}
		To simplify our notations, in this subsubsection we will still use $U_{\frac{\pi}{4}}$ instead of $U(\{0\}^{d-s+c}\times S^s).$ 
		\begin{claim}\label{mapc}
			\begin{itemize}
				\item The map $\Ga|_{S^{d+c}\setminus\{|y|=0\}}$ sends $\si^s(C)|_{ U_{\frac{\pi}{4}}}$ to $C|_{B_1^{d-s+c}}\times S^s.$
				\item The map $\Ga'|_{\R^{d-s+c}\times S^s}$ sends $C|_{B_1^{d-s+c}}\times S^s$ to $\ssc|_{ U_{\frac{\pi}{4}}}.$
			\end{itemize}
		\end{claim}
		To verify the first bullet of Claim \ref{mapc},  we can parameterize $\ssc|_{ U_{\frac{\pi}{4}}}$ by $(v,y)$ with	\begin{align*}
			&v\in C,\\&|v|\le \sin\frac{\pi}{4},\\&|v|^2+|y|^2= 1.
		\end{align*} Thus, on $\ssc|_{ U_{\frac{\pi}{4}}}$, we have
		\begin{align*}
			\Ga\big|_{\ssc|_{ U_{\frac{\pi}{4}}}}(v,y)=\bigg(\frac{v}{|y|},\frac{y}{|y|}\bigg)\in C|_{ B_1^{d-s+c}}\times S^s,
		\end{align*}where we have used the fact that $\lam v\in C$ for any $\lam\in\R_{\ge 0}.$ This finishes the first bullet of Claim \ref{mapc}
		
		For the second bullet of Claim \ref{mapc}, we can parameterize $C|_{ B_1^{d-s+c}}\times S^s$ by $(v,y)$ with
		\begin{align*}
			&v\in C,\\&|v|\le 1,\\&|y|=1.
		\end{align*}
		Thus, on $C|_{ B_1^{d-s+c}}\times S^s$, we have
		\begin{align*}
			&\Ga'\big|_{C|_{ B_1^{d-s+c}\times S^s}}(v,y)=\bigg(\frac{v}{\sqrt{1+|v|^2}},\frac{y}{\sqrt{1+|v|^2}}\bigg)\\\in &(C\times \R^{s+1})\cap  U_{\frac{\pi}{4}}=\ssc|_{ U_{\frac{\pi}{4}}}.
		\end{align*}
		Again we have used the fact that $\lam v\in C$ for any $\lam \in \R_{\ge0}$. This finishes the second bullet of Claim \ref{mapc}.
		
		Since $\Ga$ and $\Ga'$ are inverses of each other, by Claim \ref{mapc} we deduce that  $\ssc$ restricted to $U(\{0\}^{d-s+c}\times S^s)$ is equal to
		\begin{align*}
			C|_{B_1^{d-s+c}}\times S^s,
		\end{align*}via $\Ga.$
		
		Now to finish the proof, we need to determine $\sing\ssc.$ Note that $C\times \R^{s+1}$ is invariant under uniform scalings of $\R^{d+c+1},$ i.e., scalar multiples of the identity matrix.  Thus, we deduce that the regular set $\reg (C\times \R^{s+1})$ intersects the unit sphere $S^{d+c}$ in $\R^{d+c+1}$ orthogonally. Orthogonality here gives transversality of $\reg (C\times \R^{s+1})$ with $S^{d+c}$, which implies that 
		\begin{align*}
			\big(		\reg (C\times \R^{{s+1}})\big)\cap S^{d+c}\subset\reg\ssc.
		\end{align*}
		By Definition \ref{defnsm}, this implies that
		\begin{align*}
			\big(		\sing (C\times \R^{{s+1}})\big)\cap S^{d+c}\supset\sing\ssc.	
		\end{align*}The singular set of $C\times\R^{s+1}$ is $\{0\}^{d-s+c}\times\R^{s+1}.$ Consequently, we deduce that
		\begin{align*}
			\sing \ssc\s\big(\{0\}^{d-s+c}\times\R^{s+1}\big)\cap S^{d+c}=\{0\}^{d-s+c}\times S^s.
		\end{align*}
		Now recall that $\ssc$ restricted to $U_{\frac{\pi}{4}}$  equals  $(C|_{B_1^{d-s+c}})\times S^s$ via $\Ga$.  The singular set of $C|_{B_1^{d-s+c}}\times S^s$ is  $\{0\}^{d-s+c}\times S^s.$ The inverse image of $\{0\}^{d-s+c}\times S^s$, i.e., $\Ga'(\{0\}^{d-s+c}\times S^s)$ is precisely $\{0\}^{d-s+c}\times S^s.$ We are done.
		(We emphasize that the two $\{0\}^{d-s+c}\times S^s$ are in different coordinates but are the same set via $\Ga.$)
	\end{proof}
	\subsection{Adding the singular sets to make an altered representative}\label{altered}
	Now we are ready to make the altered representative of $[\Si]$ with singular sets. 
	\begin{assump}\label{assumpc}
		Assume $C$ is an area-minimizing cone in $\R^{d-s+c}$ calibrated by a smooth calibration form $\phi_C$.
	\end{assump}
	First recall that $[\Si]$ can be represented by an immersed connected submanifold $N$ with normal crossings (Assumption \ref{simpa}). 
	\begin{figure}[h]
		\centering
		\def\svgwidth{0.9\paperwidth}
		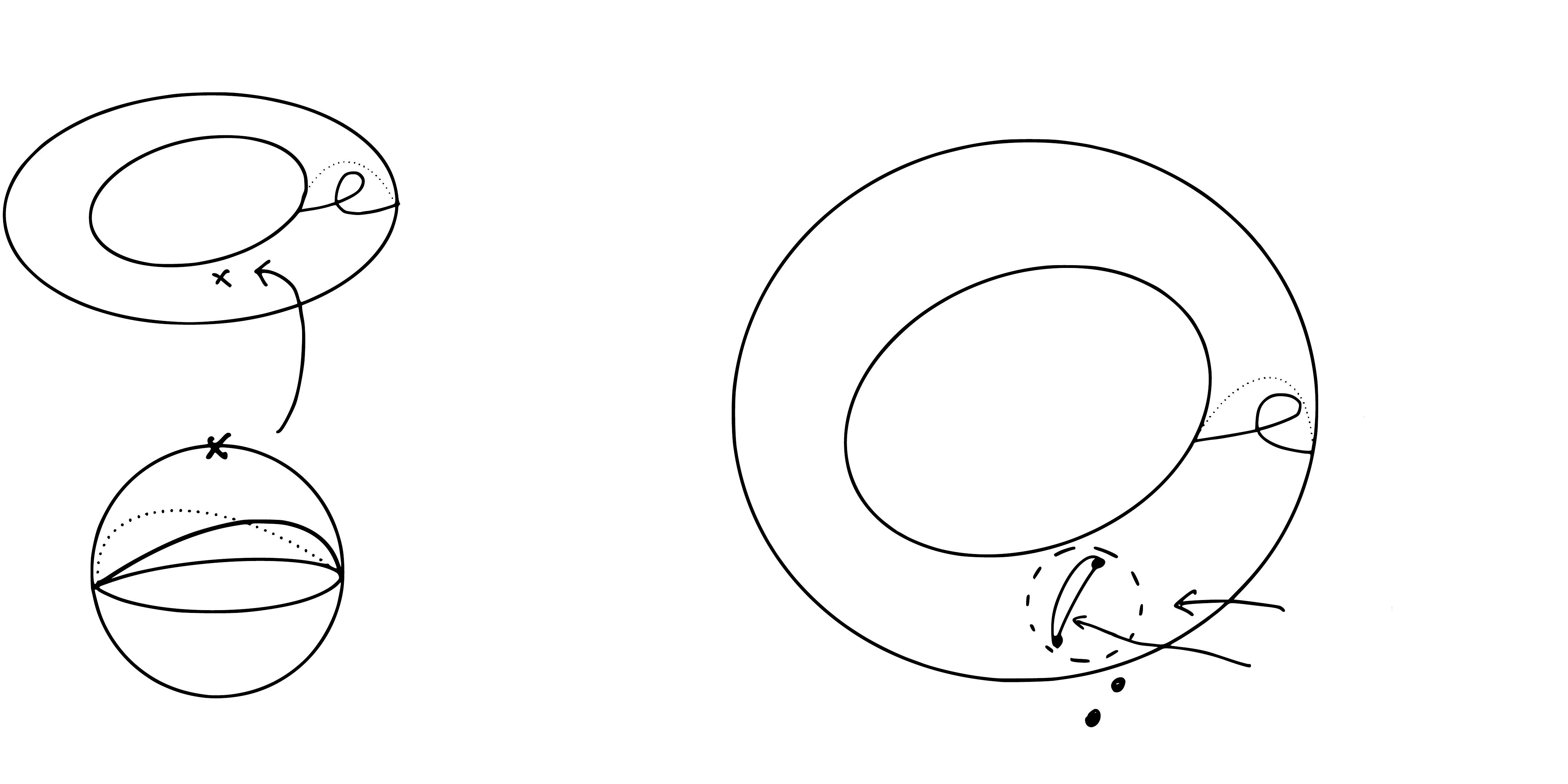
		\caption{Realizing $\ssc$ on $M$}
		\label{fig:ssccycle}
	\end{figure}
	Now, let us realize $\ssc$ as  a homologically trivial cycle on $M.$ 
	Pick any point $p$ on $S^{d+c}$ with $p\not\in \supp\ssc.$ Then $S^{d+c}\setminus\{p\}$ is diffeomorphic to an open $(d+c)$-dimensional ball. Pick any point $q$ on $M,$ with $q\not\in\supp N$ , then we can embed $S^{d+c}\setminus\{p\}$ onto an open ball centered at $q$ on $M$, not intersecting $\supp N.$ This way we can regard $\ssc$ as a homologically trivial cycle on $M,$ in the embedding of $S^{d+c}\setminus\{p\}.$ 
	
	Recall Fact \ref{imcal}. Now we can use the connected sum defined in Lemma \ref{lemcs} to construct our desired representative $N\#\ssc$ of $[\Si]$ with the following properties:
	\begin{assump}\label{assumpns}
		$N\#\ssc$ is a multiplicity $1$ integral current with \textbf{connected} regular set, and 
		\begin{align*}
			N\#\ssc|_{U(N\cap N)\cup U(\sing\ssc)}&=N|_{U(N\cap N)
			}+\ssc|_{
			U(\sing\ssc)},\\
			\sing\left( N\#\ssc\right)&=(N\cap N)\cup\sing \ssc\textnormal{ as a disjoint union}.
		\end{align*}
	\end{assump}
	Here $U(\sing\ssc)$ is just the set $U(\{0\}^{d-s+c}\times S^s)$ in Lemma \ref{singssc}. The change of notation is for coordinate independent expressions.
	
	It may happen that $\reg \ssc$ has several connected components. We can still use the connected sum in Lemma \ref{lemcs} to connect each connected component of $\reg\ssc$ to achieve a connected regular set.
	\section{Regular neighborhoods of $N\#\si^s(C)$}\label{tbssc}
	In the next section, we will need the existence of regular neighborhoods of $N\#\ssc$ of arbitrarily small radii. A regular neighborhood is a classical concept in piecewise linear topology, and basically is the simplicial version of tubular neighborhoods around simplicial subcomplexes. Heuristically speaking, we need such regular neighborhoods to have a well-controlled local geometry around $N\#\ssc,$ which will enable us to make $N\#\ssc$ area-minimizing in a smooth metric.
	
	Finding such arbitrarily small radii regular neighborhoods around smooth submanifolds is a classical fact, i.e., just taking tubular neighborhoods of arbitrarily small radii. Unfortunately, $N\#\ssc$ has fairly complicated singular set, so we need to go through a more complicated route, by realizing $N\#\ssc$ as a subcomplex of a triangulation and then using the classical regular neighborhood theorem from piecewise linear topology. This section is technical, and the reader can just keep the conclusion of Lemma \ref{tube} below in mind while skipping the proof during the first read through.
	
	We plan to achieve our goal of constructing regular neighborhoods in four steps.
	\begin{enumerate}
		\item Prove that $\supp N\#\ssc$ is a Whitney stratified set. \label{stepw1}
		\item Whitney stratified sets are triangulable, so $\supp N\# \ssc$ can be realized as a subcomplex of a triangulation of $M.$\label{stepw2}
		\item By the classical simplicial regular neighborhood theorem, the $n$-th barycentric subdivision of the triangulation of $M$  with large enough $n$ will yield regular neighborhoods $U_{\e}$ that deformation retract onto $\supp N\#\ssc$ with $U_\e$ arbitrarily close in distance to $\supp N\#\ssc.$\label{stepw3}
		\item Use Hirsch's smooth regular neighborhood theorem to upgrade the simplicial regular neighborhoods in the previous step to smooth regular neighborhoods.\label{stepw4}
	\end{enumerate}
	Now equip $M$ with an arbitrary ambient metric $h$. We are ready to state the main lemma in this section.
	\begin{lem}\label{tube}
		For any $\e>0,$ there exists a smooth open set $U_\e(N\#\ssc)$ containing $\supp N\#\ssc$, such that
		\begin{enumerate}
			\item  	 $U_\e(N\#\ssc)$ deformation retracts onto $\supp N\#\ssc$ by a map $\pi_\e$,\label{homue0}
			\item $U_\e(N\#\ssc)$ is arbitrarily close to $\supp N\#\ssc,$ i.e.,\label{homue1}
			\begin{align*}
				U_\e(N\#\ssc)\s \{p\,|\,\operatorname{dist}_h(p,\supp(N\# \ssc))<\e\}.
			\end{align*}
			\item The $d$-th homology group of $U_\e(N\#\ssc)$ is generated by $N\#\ssc$, i.e.,\label{homue}
			\begin{align}\label{homns}
				H_d(U_\e(N\#\ssc),\Z)=\Z[N\#\ssc].
			\end{align}
		\end{enumerate}
	\end{lem}
	This section will be devoted to proving Lemma \ref{tube}.
	The proof will be carried out in the order of Step (\ref{stepw1}) and Step (\ref{stepw2}) together, Step (\ref{stepw3}) and Step (\ref{stepw4}) together. Lemma \ref{tube} will be obtained as a direct consequence of the above four steps.	
	\subsection{$\supp N\#\ssc$ is a Whitney stratified set and triangulable}
	Our plan in this subsection is to stratify $\supp N\#\ssc$ into disjoint unions of smooth submanifolds of different dimensions, and then prove that the stratification is a Whitney stratification. The Whitney stratification allows us to give a triangulation of $M$ with $N\#\ssc$ as a subcomplex. 
	
	An excellent introduction to Whitney's stratification is available in \cite{DTs}, and the classical reference for triangulation of Whitney stratified sets is \cite{MG}. However, since we only use Whitney's stratification to get a triangulation, the reader can just regard Whitney's stratification as a black box. Intuitively, it is clear from the constructions alone that $\supp N\#\ssc$ can be triangulated.
	
	A Whitney stratification of a closed set $K$ is a finite filtration of closed sets $K_0\s K_1\s K_2\s\dots\s K_J=K,$ such that (the possibly empty) $K_j\setminus K_{j-1}$ is a locally finite disjoint union of $j$-dimensional submanifolds, called $j$-th strata of $K,$ and each pair of strata satisfies Whitney's condition (b). Whitney's condition (b) on an ambient manifold $M$	(\cite[Section 1]{DTg}) is defined as follows:
	\begin{defn}\label{whitb}
		For disjoint submanifolds $X,Y$ on $M$, let a point $y$ be a point in $Y\cap\ov{X}.$ We say $X$ is (b) regular over $Y$ at $y$, if given sequences $\{x_i\}$ on $X$, $\{y_i\}$ on $Y$, both converging to $y,$ such that the tangent space $T_{x_i}X$ converges to a plane $P$ in a coordinate chart, and the unit vector in the direction of $x_i$ to $y_i$ converges to $v,$ then $v\in P.$
	\end{defn}  
	The above definition of Whitney's stratification is taken from \cite[Section 2]{MG}. In some literature, authors also require strata in a Whitney stratification to satisfy Whitney's condition (a), which is long known to be implied by Whitney's condition (b), e.g., \cite[p. 454 last paragraph]{DTg}. 
	\subsubsection{Strata of $\supp N\#\ssc$}
	First we want to define a filtration $$K_0\subset K_1\subset K_2\subset\dots\subset K_d$$ for $\supp N\#\ssc$ with $K_j\setminus K_{j-1}$ a $j$-dimensional submanifold.
	
	To do this, let us discuss the singularity structure of $\supp N\#\ssc.$ By Assumption \ref{assumpns}, $\sing\left( N\#\ssc\right)$ can be decomposed into the union of $N\cap N,$ the self-intersections set of $N$, and $\sing\ssc,$ i.e.,
	\begin{align*}
		\sing\left( N\#\ssc\right)=(N\cap N)\cup \sing\ssc
	\end{align*} 
	
	Recall Lemma \ref{singssc}, $\sing\ssc$ is an embedded standard $s$-dimensional sphere $S^s$ with $s\ge 0.$ Restricted to a neighborhood $U(\sing\ssc)$ around $\sing\ssc$, $\ssc$ equals $C|_{B_1^{d-s+c}}\times S^s$ via a diffeomorphism $\Ga.$  Thus,  in our construction, we want to make sure $\sing\ssc$ is contained in  $$K_s\s K_{s+1}\s\dots\s K_d$$  and contributes to only the strata $K_s\setminus K_{s-1},$ where we regard $K_{-1}$ as an empty set.
	
	The self-intersection set $N\cap N$ is a bit more complicated. We can decompose $N\cap N$ into disjoint subsets with the same intersection multiplicity, i.e.,
	\begin{align*}
		N\cap N=(N\cap N)_2\cup (N\cap N)_3\cup\dots,
	\end{align*}where $(N\cap N)_j$ denotes the multiplicity $j$ intersection set, i.e. the self intersection points coming from the intersection of $j$ embedded pieces.
	
	In the example of the union of $xy,yz$ and $zx$-axes planes in $xyz$-space, the multiplicity $2$ intersection set $(N\cap N)_2$ consists of the $x,y$ and $z$-axis, all minus the origin. The multiplicity $3$ intersection set $(N\cap N)_3$ consists of the origin.
	
	Since immersions with normal crossings are stable mappings, it is not hard to show that at any point $p\in N\cap N,$ there is a coordinate system on $M$ centered at $p,$ in which $N\cap N$ becomes the union of at most $$\left\lfloor\frac{d}{c}\right\rfloor+1$$ $d-$dimensional coordinate axis planes  $P_1,P_2,\dots, P_k$ with $k\le\left\lfloor\frac{d}{c}\right\rfloor+1,$ such that the each pair of planes $P_i$ and $P_j$ is transverse for $i\not =j$, each plane $P_i$ is transverse to another pair's transverse intersection $P_j\cap P_l,$ which is of dimension $(d-c)$, and so on. Thus, the self-intersection set of multiplicity $j$, i.e., $(N\cap N)_j$, is an open submanifold of dimension $\left(d-(j-1)c\right),$ i.e.,
	\begin{align*}
		\dim (N\cap N)_j=d-(j-1)c.
	\end{align*} The above paragraph is taken from \cite[Section III.3]{GG} and \cite[Main Results, p. vii]{RHm}.
	
	Thus, we want to make sure that in our construction $(N\cap N)_j$ contributes only to the strata $K_{d-(j-1)c}\setminus K_{d-(j-1)c-1}.$
	
	Now we are ready to describe the filtration $K_0\s\dots\s K_d$ of $\supp N\#\ssc.$ By the above discussion, we have a decomposition of $\supp N\#\ssc$ into disjoint submanifolds of different dimensions as follows:
	\begin{align}\label{sscdec}
		\supp N\#\ssc=\reg \left(N\#\ssc\right)\cup \sing\ssc\cup (N\cap N)_2\cup \dots\cup (N\cap N)_{\left\lfloor\frac{d}{c}\right\rfloor+1},
	\end{align}
	where each factor in the right hand side is a (possibly empty) smooth submanifold of $M.$  Here we regard $(N\cap N)_{\left\lfloor\frac{d}{c}\right\rfloor+1}$ as an empty set when $N\cap N$ has no points of multiplicity $\left(\left\lfloor\frac{d}{c}\right\rfloor+1\right),$ which is equivalent to $\frac{d}{c}$ being an integer. 
	
	If either $s=0$ or $\frac{d}{c}$ is an integer, set $K_0=\sing\ssc\cup (N\cap N)_{\left\lfloor\frac{d}{c}\right\rfloor+1}.$ Otherwise, set $K_0=K_{-1}=\es.$ Suppose we have constructed a $K_k.$ If some factor in the right hand side of equation (\ref{sscdec}) has dimension $(k+1)$, then define $K_{k+1}$ as the union of $K_k$ with all the factors in the right hand side of (\ref{sscdec}) that has dimension $(k+1).$ Otherwise, set $K_{k+1}=K_k.$ Thus, inductively, we can construct a filtration 
	\begin{align}\label{filterk}
		K_0\s\dots\s K_d
	\end{align} of $\supp N\#\ssc,$ so that $K_j\setminus K_{j-1}$ is a (possibly empty)  smooth submanifold of dimension $j.$
	\subsubsection{Verifying Whitney's condition (b)}
	Now we have to verify that each pair of strata satisfy Whitney's condition (b) (Definition \ref{whitb}). 
	
	First let us decide which strata $K_j\setminus K_{j-1}$ can limit to another strata $K_l\setminus K_{l-1}$.
	
	Recall that all strata are the factors of the same dimension in the right hand side of (\ref{sscdec}). Let us consider factor by factor.
	
	Since the regular set $\reg \left(N\#\ssc\right)$ is dense, $\reg \left(N\#\ssc\right)$ can limit to any other strata. Let $\{x_j\}$ be a sequence of points in $\reg \left(N\#\ssc\right)$ that converges to a point $x.$ 
	
	If $x\in \sing\ssc,$ then let $\{y_i\}$ be a sequence of points in $\sing\ssc$ that also converges to $x.$ By Lemma \ref{singssc} and Assumption \ref{assumpns}, restricted to a neighborhood of $\sing\left( N\#\ssc\right)$, $\ssc$ equals $C|_{B_1^{d-s+c}}\times S^s\s B_1^{d-s+c}\times S^s.$ Thus, by choosing a coordinate system on $S^s,$ we can assume that there is a coordinate chart on $M$ containing $x,$ such that $N\#\ssc$ restricted to this chart equals\begin{align*}
		C\times \R^s\s \R^{d-s+c}\times \R^s
	\end{align*}
	Then the $\{y_j\}$ lies in $\{0\}^{d-s+c}\times \R^s$, and $\{x_j\}$ lies in $C\times \R^s\setminus(\{0\}^{d-s+c}\times \R^s).$

	We claim that line segment $\overrightarrow{x_iy_i}$ lies entirely on $C\times \R^s,$ so by definition of tangent space as the set of tangent vectors to curves, $\overrightarrow{x_iy_i}\in T_{x_i}\reg \left(N\#\ssc\right)$. Thus, Whitney's condition (b) is satisfied for the strata $\reg \left(N\#\ssc\right)$ over $\sing\ssc$. To see this, note that any point $q$ on $\overrightarrow{x_iy_i}$ equals\begin{align*}
		q=t\overrightarrow{0x_i}+(1-t)\overrightarrow{0y_i},
	\end{align*}for $t\in[0,1]$. Since $C\times \R^s$ is invariant under uniform scalings, $tx_i$ also lies in $C\times \R^s$. Since $C\times \R^s$ is invariant under translations generated by vectors in $\{0\}^{d-s+c}\times\R^s$, and since $(1-t)y_i\in \{0\}^{d-s+c}\times\R^s$, we see that $q\in C\times \R^s$ as claimed.

	If $x\in (N\cap N)_j,$ adopt a coordinate chart so that $N\#\ssc$ becomes the union of coordinate planes $P_1,\dots,P_j$. Let $\{y_i\}$ be a sequence in $(N\cap N)_j$ with limit $x.$ Then $\overrightarrow{x_iy_i}$ is always contained in some $P_l$ among $P_1,\dots, P_j.$ By definition of tangent space, again we deduce that, $\overrightarrow{x_iy_i}\in T_{x_i}\reg \left(N\#\ssc\right).$ Thus, Whitney's condition (b) is satisfied for the strata $\reg \left(N\#\ssc\right)$ over $(N\cap N)_j$.
	
	The strata $\sing\ssc$ is closed and thus not limiting to any other strata. The strata $(N\cap N)_j$ can only limit to $(N\cap N)_k$ for $k>j.$ Again, by adopting a coordinate chart where $N\#\ssc$ becomes a union of coordinate planes, it is straightforward to verify Whitney's condition (b) using similar ideas as the previous paragraph.
	
	To sum it up, $\supp N\#\ssc$ is a Whitney's stratified set.
	\subsubsection{Realizing $N\#\ssc$ as a subcomplex of a triangulation of $M$}Since we want to use the fact that Whitney's stratified set can be triangulated, and since we want to realize $N\#\ssc$ as a subcomplex of $M,$ we need to include $M$ into the Whitney stratification as well. To do this, set $$\supp N\#\ssc=K_d=K_{d+1}=\dots=K_{d+c-1},K_{d+c}=M.$$ Since we only added the ambient manifold $M$ to the filtration (\ref{filterk}), it is straightforward to verify that
	\begin{align*}
		K_0\s K_1\s\cdots\s K_{d+c},
	\end{align*} gives a Whitney stratification of $M.$
	
	It is a classical fact that Whitney stratified sets admit triangulations. By \cite[Theorem p. 196]{MG}, $M$ admits a finite triangulation into a simplicial complex, such that $N\#\ssc$ is a subcomplex. \cite[Theorem p. 196]{MG} is stated for Thom-Mather stratified sets, of which Whitney's stratification is a special case as mentioned in \cite[Section 1, first paragraph]{MG}. To sum it up, we have proven that
	\begin{consq}\label{subcomplex}
		There is a finite triangulation of $M$ into a simplicial complex, such	$N\#\ssc$ is a subcomplex, and $\sing\left( N\#\ssc\right)$ is a subcomplex of $N\#\ssc$.
	\end{consq}
	We have finished Step (\ref{stepw1}) and (\ref{stepw2}) of our plan of this section. Now we are ready to deal with the next two steps.
	\subsection{Smooth regular neighborhoods of $N\#\ssc$ of arbitrarily small radii}
	In this section we will prove the first two bullets of Lemma \ref{tube}, i.e., constructing regular neighborhoods of arbitrarily small radii.
	\subsubsection{Simplicial regular neighborhoods of arbitrarily small radii}
	We need the following classical regular neighborhood theorem from piecewise linear topology:
	\begin{fact}\label{regnthm}
		(\cite[Theorem 2.11]{JH})	The derived neighborhoods $O_n$ of a subcomplex $K$ in the $n$-th barycentric subdivision are regular neighborhoods of $K$, provided $n\ge 2,$
		and each $O_n$ deformation retracts onto $K.$ \end{fact}
	Here the derived neighborhood means the union of all simplices that intersect our subcomplex in the $n$-th barycentric subdivision, or in simplicial jargon, union of all stars of $K$. The deformation retract of $O_n$ onto $K$ can be constructed by collapsing each simplex intersecting our subcomplex. 
	
	We also need another classical fact about diameters of simplices in barycentric subdivisions:
	\begin{fact}\label{radbd}
		The diameter of $j$-dimensional simplices resulting from the $n$-th barycentric division is at most $\big(\frac{j}{j+1}\big)^n$ times the diameter of the original $j$-dimensional simplex.
	\end{fact}
	This is proved in \cite[p. 120]{AH}.
	
	By Consequence \ref{subcomplex}, $N\#\ssc$ can be realized as a subcomplex of a finite triangulation of $M.$ By Fact \ref{regnthm}, the $n$-th derived neighborhood $O_n$ of the simplicial realization of $N\#\ssc$ in $M$ deformation retracts onto $N\#\ssc.$ By  Fact \ref{radbd}, for any $\e>0,$ there exists $n\in\N$ large enough, such that
	\begin{align}\label{simur}
		O_n\s	\{p\,|\,\operatorname{dist}(p,\supp N\#\ssc)<\e\}.
	\end{align}
	Though Fact \ref{radbd} is stated only for simplices in Euclidean space, we only have finitely many simplices in $M,$ and each simplex in $M$ is bi-Lipschitz equivalent to a simplex in $\R^{d+c}$. Thus, we can still conclude that the diameter on $M$ of simplices in the $n$-th barycentric subdivision will converge to $0$ as $n\to\infty.$
	\subsubsection{Upgrade to smooth regular neighborhoods}
	We have obtained simplicial regular neighborhoods of $N\#\ssc$ in (\ref{simur}). Now we want to upgrade it into a smooth one. This is precisely the conclusion in \cite[conlusion (1a') p. 525]{MHr}:
	\begin{fact}
		For a simplicial regular neighborhood
		$L$ of a subcomplex $K,$ and any open set $U$ containing $L,$ there exists a smooth regular neighborhood of $K$ in $U$ that deformation retracts onto $K$.\end{fact}
	Now applying the above fact to \begin{align*}
		&K=N\#\ssc,\\& L=O_n,\\&U=\{p|\operatorname{dist}(p,\supp N\#\ssc)<\e\},
	\end{align*} we arrive at a smooth open set $U_\e(N\#\ssc)$ containing $\supp N\#\ssc$ that deformation retracts onto $\supp N\#\ssc$ by a map $\pi_\e$ and 
	\begin{align*}
		U_\e(N\#\ssc)\s \{p|\operatorname{dist}(p,\supp N\#\ssc)<\e\}.
	\end{align*}This finishes the proof of the first two bullets in Lemma \ref{tube}.
	\subsection{Homology of $U_\e(N\#\ssc)$}
	Now we are left to prove the last bullet in Lemma \ref{tube}. 
	Since $U_\e(N\#\ssc)$ deformation retracts onto $\supp N\#\ssc,$ to show that $H_d(U_\e(N\#\ssc))=\Z[N\#\ssc],$ it suffices to prove that
	\begin{align*}
		H_d(\supp N\#\ssc)=\Z[N\#\ssc].
	\end{align*}
	Now we need the following classical constancy theorem from geometric measure theory \cite[Theorem 4.9]{FMgmt}:
	\begin{fact}\label{const}
		A $d$-dimensional integral current $T$ supported on a $d$-dimensional connected not necessarily closed orientable manifold $N$ equals an integer multiple of $N,$ provided $\pd T $ is supported on $\pd N.$
	\end{fact}
	Let $W$ be a simplicial cycle supported on $\supp N\#\ssc$. Then $W$ induces an integral current supported on $\supp N\#\ssc$. On the connected open manifold $\reg \left(N\#\ssc\right)$, we can use the Fact \ref{const} to deduce that $W$ restricted to $\reg \left(N\#\ssc\right)$ is an integer multiple of $N\#\ssc$ restricted to $\reg \left(N\#\ssc\right).$ Now that $\reg \left(N\#\ssc\right)$ contains all open $d$-dimensional simplices in $\supp N\#\ssc$, we deduce that $W=k(N\#\ssc)$, with $k\in\Z.$ This shows that $H_d(\supp N\#\ssc)$ is generated by $[N\#\ssc].$ 
	
	To finish the proof of Lemma \ref{tube}, we need to show that $[N\#\ssc]$ is not a torsion class on $\supp N\#\ssc.$ Note that there is no $(d+1)$-dimensional simplex in the triangulation of $N\#\ssc.$ Thus, a $d$-cycle on $N\#\ssc$ is a boundary if and only if it equals to $0$ as a cycle. However, $kN\#\ssc$ for any $k\not=0$ is not equal to $0$ as a cycle. We are done.
	\subsection{Remarks on $U(N\cap N)$ defined in Fact \ref{imcal}}
	Recall Consequence \ref{subcomplex} that $\sing\left( N\#\ssc\right)$, and thus $N\cap N$ are also a subcomplex of $M.$ The proof of Lemma \ref{tube} also applies to $N\cap N$. Consequently, without loss of generality, we can assume that
	\begin{assump}\label{assumpnn}
		$U(N\cap N)$  deformation retracts onto $N\cap N$.
	\end{assump}
	\section{Making the altered representative area-minimizing}\label{zhang}
	This section will be devoted to the proof of the following lemma:
	\begin{lem}\label{nsmin}
		There is a smooth metric $h$ on $M,$ such that $ N\#\ssc$ is the unique area-minimizing representative of $[\Si]$ with respect to $h$, provided $2[\Si]\not=0.$ When $2[\Si]=0,$ $\pm N\#\ssc$ are the unique area-minimizing representatives of $[\Si]$ in $h.$  Furthermore, we have
		\begin{itemize}
			\item 
			In a neighborhood $U_\rr(\sing\ssc)$ of $\sing \ssc$ diffeomorphic to $B_\rr^{d-s+c}\times S^s,$ via the diffeomorphism $\Ga$ defined in Lemma \ref{singssc},  the metric $h$ is equal to the pullback under $\Ga$ of the standard metric on $B_{\rr}^{d-s+c}\times S^s,$ and $N\#\ssc$ restricted to $U_\rr(\sing \ssc)$ equals to $C|_{B_\rr^{d-s+c}}\times S^s$ via $\Ga.$
		\end{itemize}
	\end{lem} 
	Here $B_\rr^{d-s+c}\times S^s$ is the product of a radius $\rr$ ball $B_\rr^{d-s+c}$ inside $\R^{d-s+c}$ and the standard $s$-dimensional sphere $S^s$ with $\rr>0$, and
	\begin{align*}
		U_\rr(\sing\ssc)\overset{\Ga}{\cong}B_\rr^{d-s+c}\times S^s.
	\end{align*}
	
	The negative sign $-$ means orientation reversing. Thus, $N\#\ssc$ and $-N\#\ssc$ are indeed different as integral currents. 
	
	The proof relies on Zhang's gluing of calibrations in  \cite{YZa,YZj,YZt}. For the reader's convenience, we reproduce the entire argument here.
	
	Roughly speaking, the construction of the metric $h$ involves five steps:
	\begin{enumerate}
		\item The way we construct $N\#\ssc$ provides a Riemannian metric $h_\sing$ and a smooth calibration form $\phi_\sing$ that calibrates $N\#\ssc$ restricted to $U(N\cap N)\cup U(\sing\ssc).$ Then we handcraft a Riemannian metric $h_\reg$ and a calibration form $\phi_\reg$ that calibrates the regular set of $N\#\ssc$.\label{z1}
		\item Gluing  the metrics $h_\reg$ to $h_\sing$ and the calibrations $\phi_\reg$ to $\phi_\sing$, we can construct a smooth metric $\hgd$ and a smooth form $\phi$ calibrating $N\#\ssc$ in a smooth regular neighborhood $U_\e(N\#\ssc)$ constructed in Lemma \ref{tube}.\label{z2}
		\item Multiplying the metric by a bump function to make $N\#\ssc$ the unique area-minimizing representative of $[N\#\ssc]$ in $U_\e(N\#\ssc)$. \label{z3}
		\item Making the metric in the complement of $U_\e(N\#\ssc)$ very large, so that any area-minimizing integral current $T$  in $[\Si]$ will stay in $U_\e(N\#\ssc).$\label{z4}
		\item Though homologous objects restricted to subsets of ambient space are not necessarily homologous, the equation (\ref{homns})  implies that $[T]= k[N\#\ssc]$ with $0\not=k\in \Z$ for any integral current cycle $T$ supported in $U_\e(N\#\ssc)$. Now use the calibration $\phi$ to deduce that $$|k|\ms(N\#\ssc)=|kN\#\ssc(\phi)|=|T(\phi)|\le \ms(T).$$ Then Step (\ref{z3}) forces $T=\pm N\#\ssc$ provided $T$ is area-minimizing in $[\Si].$
		\label{z5}
	\end{enumerate}
	After doing some basic preparational works, our proof will be carried in the order of the above five steps.
	\subsection{Preparing the neighborhoods}
	First of all, we need to prepare the neighborhoods $U_\e(N\#\ssc),U(N\cap N), U(\sing\ssc)$ for our constructions.
	\begin{defn}\label{defnus}
		Define
		\begin{align}\label{defnuns}
			U\left(\sing\left( N\#\ssc\right)\right)=U(N\cap N)\cup U(\sing\ssc).
		\end{align}
	\end{defn}Note that the union above is a disjoint union.
	
	In order to make space for gluing calibrations in later subsections, we need to retract $U\left(\sing\left( N\#\ssc\right)\right)$ a bit into its interior. Since $U\left(\sing\left( N\#\ssc\right)\right)$ is an oriented smooth submanifold with boundary, a neighborhood of its boundary $\pd\big(U\left(\sing\left( N\#\ssc\right)\right)\big)$ in its interior is diffeomorphic to $$\pd\big(U\left(\sing\left( N\#\ssc\right)\right)\big)\times [0,1).$$ We can delete this $\pd\big(U\left(\sing\left( N\#\ssc\right)\right)\big)\times [0,1)$ 
	to obtain a slightly smaller open set $\hat{U}\left(\sing\left( N\#\ssc\right)\right).$ In other words we have the following fact
		\begin{figure}[h]
		\centering
		\def\svgwidth{0.9\paperwidth}
		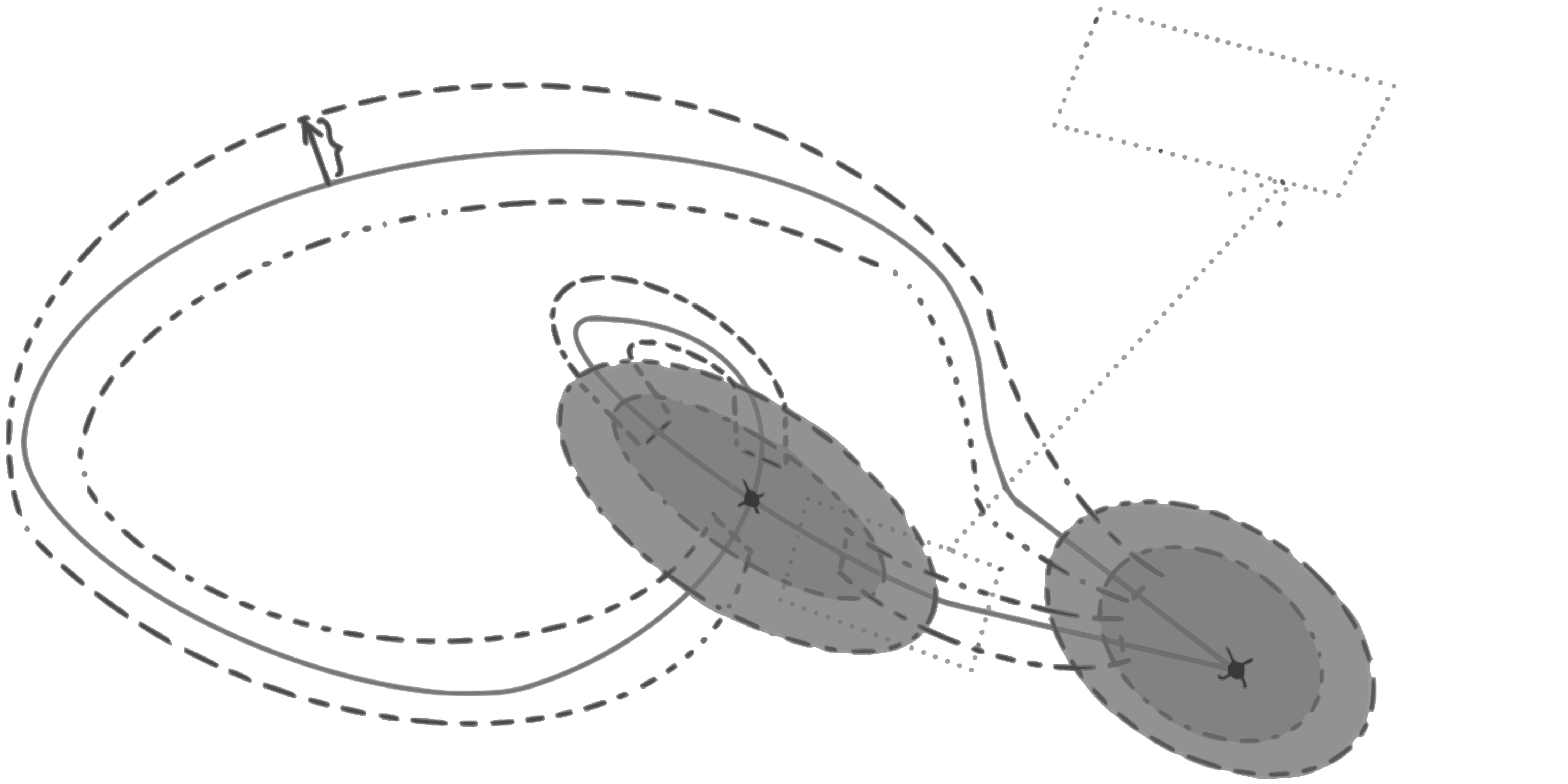
			\caption{Schematic illustration of the neighborhoods\\We emphasize that $\reg \left(N\#\ssc\right)$ should be connected! Impossible to draw on a $2$-d plane.}
	\label{fignbhd}
	\end{figure}
	\begin{fact}\label{uinner}
		The smooth open set $\hat{U}\left(\sing\left( N\#\ssc\right)\right)$ containing $\sing\left( N\#\ssc\right)$, is an open subset of $U\left(\sing\left( N\#\ssc\right)\right)$ and satisfies \begin{align*}
			U\left(\sing\left( N\#\ssc\right)\right)\setminus\hat{U}\left(\sing\left( N\#\ssc\right)\right)\overset{\Lambda}{\cong}\pd\big(U\left(\sing\left( N\#\ssc\right)\right)\big)\times (0,1],
		\end{align*}with $\Lambda$ a smooth diffeomorphism such that
		\begin{align*}
			\pd\big(U\left(\sing\left( N\#\ssc\right)\right)\big)&\overset{\Lambda}{\cong}\pd\big(U\left(\sing\left( N\#\ssc\right)\right)\big)\times\{0\},\\
			\pd\big(\hat{U}\left(\sing\left( N\#\ssc\right)\right)\big)&\overset{\Lambda}{\cong}\pd\big(U\left(\sing\left( N\#\ssc\right)\right)\big)\times\{1\}.
		\end{align*}
	\end{fact}
	The reader should now consult Figure \ref{fignbhd} for intuition.
	\subsection{Calibrating $N\#\ssc$ restricted to $U\left(\sing\left( N\#\ssc\right)\right)$}
	In this subsection, we will show how from our construction of $N\#\ssc$, we naturally get a smooth metric $h_\sing$ defined throughout $M$ such that $N\#\ssc$ is calibrated by a smooth form $\phi_\sing$ defined on $U\left(\sing\left( N\#\ssc\right)\right).$
	
	We need the following classical fact.
	\begin{fact}\label{calc}
		The closed smooth form $\pi_{\R^{d-s+c}}\du\phi_C\w\pi_{S^s}\du \operatorname{dvol}_{S^s}$ calibrates $C\times S^s$ in $\R^{d-s+c}\times S^s$ with the product Riemannian metric of standard metrics on the two factors.
	\end{fact}
	Here $\phi_C$ is the smooth form calibrating $C$, defined in Assumption \ref{assumpc}, and $\dvol_{S^s}$ is the volume form of $S^s$, $\pi_{\R^{d-s+c}}$ and $\pi_{S^s}$ are projections onto the $\R^{d-s+c}$ and $S^s$ factors, respectively, in the product $\R^{d-s+c}\times S^s.$ The above fact follows directly from \cite[Chapter II, Proposition 7.10]{HL}, which says that the wedge product of a calibration form $\psi$ and a unit simple  form remains a calibration in the product metrics, and calibrates precisely the product of calibrated submanifolds of $\psi$ and the submanifolds with tangent spaces dual to the unit simple form.
	
	Now, using Lemma \ref{singssc} we see that $N\#\ssc$ is calibrated in $U(\sing\ssc)$ by a smooth form $$\phi_\sing=\Ga\du(\pi_{\R^{d-s+c}}\du\phi_C\w\pi_{S^s}\du \operatorname{dvol}_{S^s})$$ in a smooth metric $h_\sing$ that is the pullback of the standard product metric on $\R^{d-s+c}\times S^s$ under $\Ga$. By Fact \ref{imcal}, $N\#\ssc$ is calibrated in $U(N\cap N)$ by a smooth form in a smooth metric. 
	
	Recall Definition \ref{defnuns}. Now we can sum up the construction in this section:
	\begin{fact}\label{calsing}
		We have a smooth Riemannian metric $h_\sing$ and a smooth form $\phi_\sing$ on $U\left(\sing\left( N\#\ssc\right)\right)$ such that $N\#\ssc$ is calibrated in $U\left(\sing\left( N\#\ssc\right)\right)$ by $\phi_\sing$ with respect to $h_\sing.$ 
	\end{fact}
	Now extend $h_\sing$ on $U\left(\sing\left( N\#\ssc\right)\right)$ to a smooth Riemannian metric $h_\sing$ defined throughout $M.$
	\subsection{Calibrating the regular set  of $N\#\ssc$}
	In this subsection, we want to calibrate the regular set of $N\#\ssc$ in a smooth metric. A natural idea pioneered in \cite{YZj} is to use the pullback of the volume form of $N\#\ssc$ under the normal bundle exponential map of $N\#\ssc$ as a calibration form in a renormalized metric.  Unfortunately, $\reg \left(N\#\ssc\right)$ is an open manifold with limit points in $\sing\left( N\#\ssc\right).$ Thus, we have to go through a more complicated route by restricting to a compact subset of $\reg \left(N\#\ssc\right)$ and taking into account that we have to make space for gluing of calibrations in the next subsection.
	
	Let us first discuss what compact subsets of $\reg \left(N\#\ssc\right)$ we want to deal with. Write
	\begin{align*}
		\widecheck{\reg}(N\#\ssc)=(\reg \left(N\#\ssc\right) )\cap \hat{U}\left(\sing\left( N\#\ssc\right)\right)\cp.
	\end{align*}
	By transversality, we can assume that $\widecheck{\reg}(N\#\ssc)$ is a smooth submanifold with boundary.
	
	Let $\exp^\perp$ be the exponential map from the normal bundle $\T^\perp\reg \left(N\#\ssc\right)$ of $\reg \left(N\#\ssc\right)$ to $M.$ Here $\T^\perp$ means the normal bundle of a submanifold with respect to the ambient manifold $M.$ Without loss of generality, we can assume that $\exp^\perp$ is a diffeomorphism when restricted to the following ball bundle \begin{align*}
		B^\perp_\de \widecheck{\reg}(N\#\ssc)=	\{(p,v)\,|\,p\in \widecheck{\reg}(N\#\ssc),\textnormal{ }v\in \T_p^\perp \reg \left(N\#\ssc\right),\textnormal{ }h_\sing(v,v)\le \de^2\},
	\end{align*} over the base $\widecheck{\reg}(N\#\ssc).$
	
	The image of $B^\perp_\de \widecheck{\reg}(N\#\ssc)$ under $\exp^\perp$ is a set of points with distance at most $\de$ to $\widecheck{\reg}(N\#\ssc)$. Write
	\begin{align*}
		B_\de(\widecheck{\reg}(N\#\ssc))=\exp^\perp B^\perp_\de \widecheck{\reg}(N\#\ssc).
	\end{align*}
	Without loss of generality, we can assume that
	\begin{align}\label{inc}
		\hat{U}\left(\sing\left( N\#\ssc\right)\right)\cp\cap U_\e(N\#\ssc)\s 	B_\de(\widecheck{\reg}(N\#\ssc)).
	\end{align}
	See Figure \ref{fignbhd} to get an intuition. This seemingly odd (\ref{inc}) is to prepare for the gluing of calibrations in the next subsection.
	
	Let $\widecheck{\pi}$ denote the projection onto base $\widecheck{\reg} (N\#\ssc)$ in the ball bundle $B_\de^\perp\widecheck{\reg}(N\#\ssc).$ Then on $B_\de(\widecheck{\reg}(N\#\ssc))$, define
	\begin{align*}
		{\pi}=\widecheck{\pi}\circ(\exp^\perp)\m.
	\end{align*}
	In other words ${\pi}$ is the nearest distance projection of $B_\de(\widecheck{\reg}(N\#\ssc))$ onto $\widecheck{\reg}(N\#\ssc).$ The projection ${\pi}$ provides a smooth orthogonal splitting of the tangent spaces to $M$ in $B_\de(\widecheck{\reg}(N\#\ssc))$ into
	\begin{align}\label{split}
		\ker \ed\pi\oplus \ker^\perp \ed\pi.
	\end{align}
	Here $\oplus$ is the Riemannian direct sum of vector spaces, and $\ker^\perp \ed\pi$ is the orthogonal complement of $\ker \ed\pi.$ By definition of $\pi,$ $\ker\ed\pi$ are tangent spaces to the level sets of $\pi,$ i.e., the submanifolds formed by geodesics orthogonal to $\reg \left(N\#\ssc\right).$ Thus, $\kn$ restricted to $\reg \left(N\#\ssc\right)$
	is the tangent space to $\reg \left(N\#\ssc\right)$. Heuristically, one should think of $\kn$ as the natural extension of the tangent space of $\reg \left(N\#\ssc\right)$ into the ball bundle.

	This also provides a natural smooth splitting of the Riemannian metric $h_\sing$ into
	\begin{align*}
		h_\sing=(h_\sing)|_{\ker \ed\pi}\oplus(h_\sing)|_{\ker^\perp \ed\pi}.
	\end{align*}Here $\oplus$ means the Riemannian sum of the metrics on orthogonal subspaces. 
	
	In $B_\de(\widecheck{\reg}(N\#\ssc))$
	the pullback of the volume form $\dvol_{N\#\ssc}$ of $\reg \left(N\#\ssc\right)$ under $\pi$, i.e., $${\pi}\du(\dvol_{N\#\ssc}),$$ is a closed simple form,  because the form $\dvol_{N\#\ssc}$ is a closed simple $d$-form and consequently its pullback under the projection ${\pi}$ is also closed and simple. With respect to the splitting (\ref{split}), the form ${\pi}\du(\dvol_{N\#\ssc})$ can be written as
	\begin{align*}
		{\pi}\du(\dvol_{N\#\ssc})=\big({\pi}\du(\dvol_{N\#\ssc})(\ker^\perp \ed\pi)\big)\big(\ker^\perp \ed\pi\big)\du.
	\end{align*}Here $\ker^\perp \ed\pi$ also denotes the unit simple $d$-vector representing $\ker^\perp \ed\pi$ and
	$(\ker^\perp \ed\pi)\du$ denotes the unit simple $d$-form Riemannian dual to the unit simple $d$-vector $\ker^\perp \ed\pi.$ Then we have the following fact:
	\begin{fact}\label{calreg}
		The smooth closed form $$\phi_\reg={\pi}\du(\dvol_{N\#\ssc})$$ is a calibration form defined on $$B_\de(\widecheck{\reg}(N\#\ssc))$$ that calibrates $ N\#\ssc$ in	$B_\de(\widecheck{\reg}(N\#\ssc))$, with respect to the metric
		$$h_\reg=(\phi_\reg(\ker^\perp \ed\pi))^{\frac{2}{d}}h_{\sing}.$$
	\end{fact}
	The above fact is taken from \cite[Remark 3.5]{YZa}. The proof is a simple calculation using definition of comass, i.e., bullet (\ref{cms2}) of Fact \ref{cmsvec}.
	\subsection{Calibrating $N\#\ssc$ in $U_\e(N\#\ssc)$ by gluing calibrations}
	Now we have come to the crux of the matter. We will glue the metrics and calibrations we have obtained so far to make $N\#\ssc$ area-minimizing in $U_\e(N\#\ssc).$ This subsection is taken from the proof of \cite[Theorem 4.6]{YZa}.   
	
	By equation (\ref{inc}), Fact \ref{calreg} and Fact \ref{calsing}, in the region
	\begin{align*}
		\big(U\left(\sing\left( N\#\ssc\right)\right)\setminus\hat{U}\left(\sing\left( N\#\ssc\right)\right)\big)\cap U_\e(N\#\ssc),
	\end{align*}
	we have two Riemannian metrics $h_\reg,h_\sing$ and two smooth calibrations forms $\phi_\reg,\phi_\sing$, so that $N\#\ssc$ is calibrated in this region by the corresponding forms in the corresponding metrics, respectively. Now our goal is to glue these pairs of Riemannian metrics and calibration forms together. The gluing is based on the following central fact:
	\begin{fact}\label{calk}
		Let $g$ be a positive definite quadratic form on $\R^{d+c}$, $\psi$ be a constant $d$-form on $\R^{d+c}$ and $\R^{d+c}=P\oplus Q$ be a direct sum decomposition of $\R^{d+c}$ into orthogonal complementary subspaces $P$ of dimension $d$ and $Q$ of dimension $c$ with respect to $g$, such that $\psi(P)>0.$ Then in the metric
		\begin{align*}
			h=(\psi(P))^{\frac{2}{d}}\big(g|_P\oplus (K g|_Q)\big),
		\end{align*}$\psi$ has comass
		\begin{align*}
			\cms_{h}\psi=1,
		\end{align*} 
		and calibrates $P,$ provided the real number $K$ satisfies
		\begin{align}\label{kge}
			K\ge \binom{d+c}{d}(\psi(P))^{-1}\cms_g\psi.
		\end{align}
	\end{fact}
	Here $\psi(P)$ mean the value of $\psi$ at the simple unit $d$-vector representing $P$ with respect to $g.$
	The above fact is just a restatement of \cite[Lemma 3.4, Remark 3.5]{YZa}, which is in turn a restatement of a classical result on comass by Harvey-Lawson \cite[Lemma 2.14]{HLf}. 
	
	Roughly speaking, in the glued calibrations, the only component  that matters to us is $\psi(\kn)(\kn)\du.$ Then Fact \ref{calk} will enable us to basically annihilate the effects of other terms spanned by covectors in $\ker\ed\pi\du$ in the glued calibration form. 
	\subsubsection{Notation convention}
	From now on, we will restrict our attention to the region $\big(U\left(\sing\left( N\#\ssc\right)\right)\setminus\hat{U}\left(\sing\left( N\#\ssc\right)\right)\big)\cap U_\e(N\#\ssc)$. We adopt the notation:
	\begin{align*}
		U_{\operatorname{gluing}}=\big(U\left(\sing\left( N\#\ssc\right)\right)\setminus\hat{U}\left(\sing\left( N\#\ssc\right)\right)\big)\cap U_\e(N\#\ssc).
	\end{align*}
	Recall that $h_\sing$ is defined throughout $M$ and $h_\reg$ is obtained by rescaling $h_\sing$ (Fact \ref{calreg}). 
	
	Thus, from now on, we will use $h_\sing$ as our base metric and use the symbols
	\begin{align*}
		h_{\ker \ed\pi}=&(h_\sing)|_{\ker \ed\pi}, \\h_{\ker^\perp \ed\pi}=&(h_\sing)|_{\ker^\perp \ed\pi}.
	\end{align*}
	Then by definition we have
	\begin{align*}
		h_\sing=\hkp \oplus \hkn.
	\end{align*}
	
	\subsubsection{Anti-derivatives of forms}
	By Lemma \ref{singssc}, $U(\sing\ssc)$ deformation retracts onto $\sing\ssc.$ On the other hand, by Assumption \ref{assumpnn}, $U(N\cap N)$ deformation retracts onto $N\cap N.$ Thus, $U\left(\sing\left( N\#\ssc\right)\right)$ deformation retracts onto $\sing\left( N\#\ssc\right).$ Since $\sing\left( N\#\ssc\right)$ is a subcomplex of dimension less than $d,$ we deduce that
	\begin{align*}
		H_d(U\left(\sing\left( N\#\ssc\right)\right),\Z)=0.
	\end{align*}
	By the universal coefficient theorem \cite[Theorem 3.2]{AH}, this implies that
	\begin{align*}
		H^d(U\left(\sing\left( N\#\ssc\right)\right),\R)=0.
	\end{align*}
	By De Rham's theorem for smooth manifolds with boundary \cite[Theorem 2.48]{RM}, we deduce that there is a smooth $(d-1)$-form $\Phi_\sing$ on $U\left(\sing\left( N\#\ssc\right)\right)$ such that
	\begin{align*}
		d\Phi_\sing=\phi_\sing.
	\end{align*}
	
	On the other hand, by transversality, we can assume that $\ug\cap \reg \left(N\#\ssc\right)$ is a smooth $d$-dimensional manifold with boundary, thus
	\begin{align*}
		H_d(\ug\cap \reg \left(N\#\ssc\right),\Z)=0.
	\end{align*} Again, by the universal coefficient theorem and De Rham's theorem for manifolds with boundary, we deduce that $\dvol_{N\#\ssc}$ is exact on $\ug\cap \reg \left(N\#\ssc\right)$. Since $\phi_\reg$ is the pullback of $\dvol_{N\#\ssc}$ under $\pi,$ we deduce that there is a smooth-$(d-1)$ form $\Phi_\reg$ on $\ug$, such that
	\begin{align*}
		d\Phi_{\reg}=\phi_\reg.
	\end{align*}
	\subsubsection{Dividing $\ug$ into four regions}
	By Fact \ref{uinner}, we have
	\begin{align*}
		U\left(\sing\left( N\#\ssc\right)\right)\setminus\hat{U}\left(\sing\left( N\#\ssc\right)\right)&\overset{\Lambda}{\cong}\pd\big(U\left(\sing\left( N\#\ssc\right)\right)\big)\times (0,1],\\
		\pd\big(U\left(\sing\left( N\#\ssc\right)\right)\big)&\overset{\Lambda}{\cong}\pd\big(U\left(\sing\left( N\#\ssc\right)\right)\big)\times\{0\},\\
		\pd\big(\hat{U}\left(\sing\left( N\#\ssc\right)\right)\big)&\overset{\Lambda}{\cong}\pd\big(U\left(\sing\left( N\#\ssc\right)\right)\big)\times\{1\}.
	\end{align*}
	Let $\rh$ be the smooth function defined by
	\begin{align*}
		\rh=\pi_{(0,1]}\circ\Lambda.
	\end{align*}
	Here $\pi_{(0,1]}$ is the projection onto the $(0,1]$ factor of $\pd\big(U\left(\sing\left( N\#\ssc\right)\right)\big)\times (0,1].$
		\begin{figure}[h]
		\centering
		\def\svgwidth{0.8\paperwidth}
		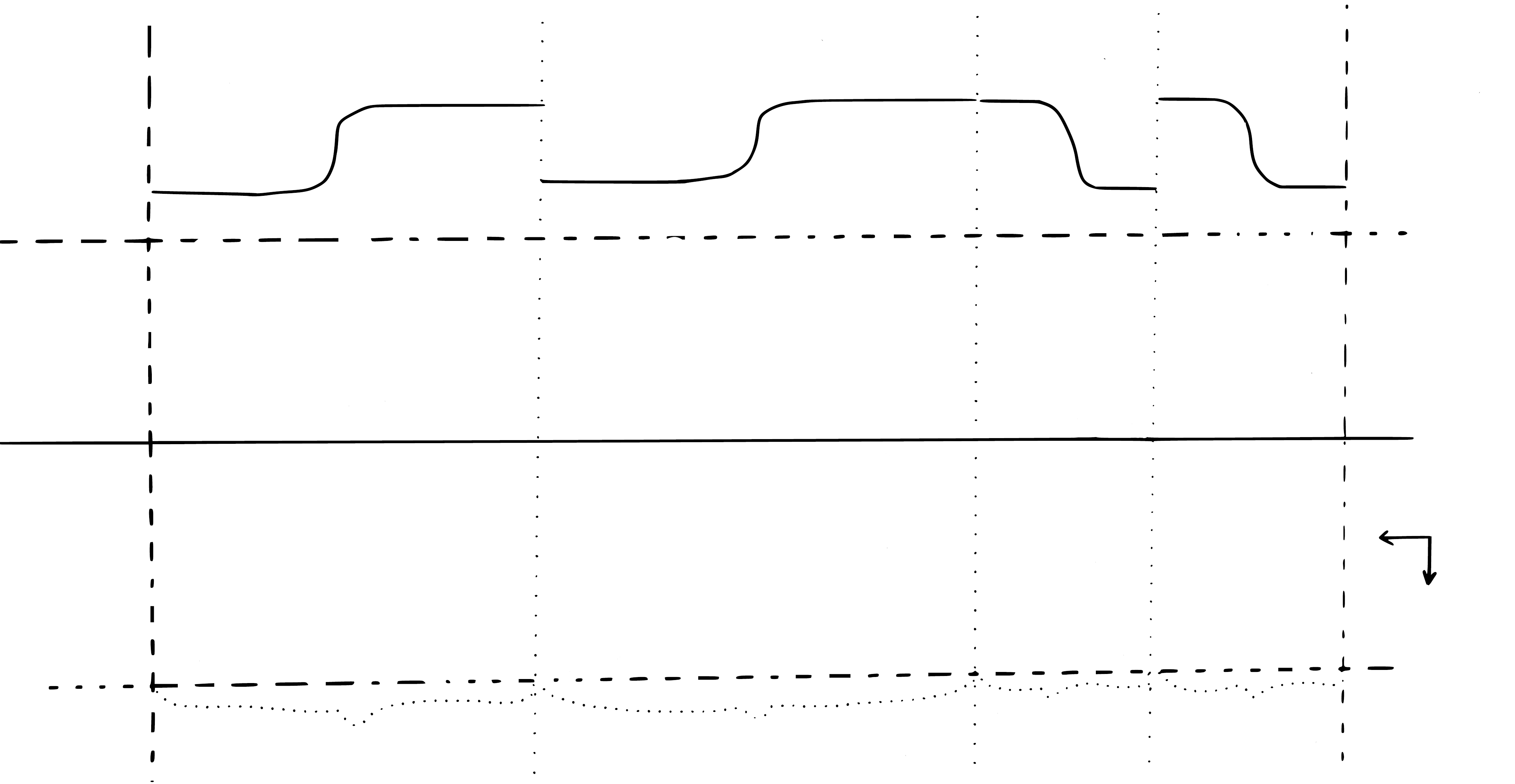
		\caption{Gluing of forms and the change of metrics}
		\label{figgluing}
	\end{figure}
	\begin{defn}
		For any interval $(a,b)\s (0,1]$, define
		\begin{align*}
			U_{(a,b)}=\rh\m(a,b)\cap \ug.
		\end{align*}
	\end{defn}
	We will single out four regions of $\ug$,\newcommand{\ua}{U_{(0,\frac{1}{3})}}\newcommand{\ub}{U_{(\frac{1}{3},\frac{2}{3})}}\newcommand{\uc}{U_{(\frac{2}{3},\frac{5}{6})}}\newcommand{\ud}{U_{(\frac{5}{6},1)}}
	\begin{align*}
		U_{(0,\frac{1}{3})}&=\rh\m \bigg(0,\frac{1}{3}\bigg)\cap\ug,\\
		U_{(\frac{1}{3},\frac{2}{3})}&=\rh\m\bigg(\frac{1}{3},\frac{2}{3}\bigg)\cap \ug,\\
		U_{(\frac{2}{3},\frac{5}{6})}&=\rh\m\bigg(\frac{2}{3},\frac{5}{6}\bigg)\cap \ug,\\
		U_{(\frac{5}{6},1)}&=\rh\m\bigg(\frac{5}{6},1\bigg)\cap \ug.
							\end{align*}The reader should now consult Figures  \ref{fignbhd} and \ref{figgluing} for intuition.
	\subsubsection{Auxiliary functions}\label{auxf}
	Let $\ai:\R\to\R$ be a smooth monotonic function that is $0$ on $(-\infty,0]$ and $1$ on $[1,\infty)$. Define
	\begin{align*}
		\de_0=&\ai\circ(3{\rh}{}),\\
		\lam=&\ai\circ(3\rh-1),\\
		\de_\perp=&\ai\circ(-6\rh+4),\\
		\de_K=&\ai\circ(-6\rh+5).
	\end{align*}
	Then $\de_0,\lam,\de_\perp,\de_K$ are smooth functions defined on $\ug$ with values between $[0,1]$, and not equal to $0,1$ only on $\ua,\ub,\uc,\ud$ respectively. For an intuitive illustration, please refer to Figure \ref{figgluing}.
	\subsubsection{Gluing of forms}
	In $\ub$, define\begin{align*}
		\phi_\lam=d(\lam\Phi_\sing+(1-\lam)\Phi_\reg)=\lam\phi_\sing+(1-\lam)\phi_\reg+d\lam\w (\Phi_\sing-\Phi_\reg).
	\end{align*}
	Then by definition of $\lam,$ it is straightforward to verify that 
	\begin{align}\label{defp}	
		\phi=
		\begin{cases}
			\phi_\sing,&\text{on } \ov{\uc}\cup\ov{\ud}\cup \hat{U}\left(\sing\left( N\#\ssc\right)\right),\\
			\phi_\lam,&\text{on }\ub,\\
			\phi_\reg,&\text{on } \left(B_\de(\widecheck{\reg}N\#\ssc)\setminus \ug\right)\cup \ov{\ua},
		\end{cases}
	\end{align}
	gives a closed smooth $d$-form when restricted to $U_\e(N\#\ssc)$. The reader should now consult Figure \ref{figgluing} for intuition.
	\subsubsection{Gluing of metrics}
	In the following, when we write $\ker^\perp \ed\pi$ as the input of a differential form, we always mean the unit simple $d$-vector corresponding to $\kn$ in the base metric $h_\sing.$ Similarly, when we write $(\kn)\du$, we always mean the unit simple $d$-vector Riemannian dual to the unit simple $d$-vector $\kn$ in the base metric $h_\sing.$
	
	In $\ug,$ with respect to the orthogonal decomposition $\kp\oplus\kn$, we can always write uniquely
	\begin{align*}
		\phi=(\phi(\kn))(\kn)\du+\psi,
	\end{align*}with $\psi$ spanned by covectors in $\kp.$ Note that $\psi\equiv 0$ on $\ua.$ 
	
	By construction, 
	\begin{align*}
		\phi(\kn)\equiv 1,
	\end{align*}on $\ug\cap \reg \left(N\#\ssc\right).$
	By taking $\e$ in $U_\e(N\#\ssc)$ small, we can assume that
	\begin{align}\label{philbd}
		\phi(\kn)\ge \frac{1}{2}
	\end{align}  in $\ug.$ Let
	\begin{align}
		\label{defk}	K=\binom{d+c}{d}\frac{\cms_g\phi|_{\ug}}{\inf_{p\in \ug}\phi_p(\kn)}.
	\end{align}
	We define a Riemannian metric $\hgd$ defined by:
	\begin{align}h_{\operatorname{glued}}=
		\begin{cases}
			h_\sing=\hkn\oplus \hkp,&\text{on }\hat{U}\left(\sing\left( N\#\ssc\right)\right),\\
			\hkn\oplus \big((1+K\de_K)\hkp\big),&\text{on }\ov{\ud},\\
			\left((1-\de_\perp)+\de_\perp(\phi_\sing(\kn))\right)^{\frac{2}{d}}\Big(\hkn\oplus \big((1+K)\hkp\big)\Big),&\text{on }{\uc},\\
			(\phi_\lam(\kn)^{\frac{2}{d}})\Big(\hkn\oplus \big((1+K)\hkp\big)\Big),&\text{on }\ov{\ub},\\
			(\phi_\reg(\kn)^{\frac{2}{d}})\Big(\hkn\oplus \big((1+\de_0 K)\hkp\big)\Big),&\text{on }{\ua},\\
			h_\reg=(\phi_\reg(\kn)^{\frac{2}{d}})\big(\hkn\oplus \hkp\big),&\text{on } B_\de(\widecheck{\reg}(N\#\ssc))\setminus \ug.				
		\end{cases}
	\end{align}
	\begin{claim}\label{hgdc}
		The metric $h$	is a smooth Riemannian metric when restricted to $U_\e(N\#\ssc)$ and $\phi$ defined in (\ref{defp}) is a calibration form in metric $h$ that calibrates $N\#\ssc$. 
	\end{claim}
	The reader should consult Figure \ref{figgluing} for intuition.
	
	To verify the above claim, the smoothness of $\hgd$ follows directly from the definitions of $\de_K,\de_\perp,\lam,\de_0$ in Section \ref{auxf} and the definition of $\phi$ in (\ref{defp}).
	
	Note that $\phi\equiv\phi_\sing,\hgd\equiv h_\sing$ on $\hat{U}\left(\sing\left( N\#\ssc\right)\right)$ and $\phi\equiv \phi_\reg,\hgd\equiv h_\reg$ on $B_\de(\widecheck{\reg}(N\#\ssc))\setminus \ug$. Thus, by Fact \ref{calsing} and Fact \ref{calreg}, $\phi$ is a calibration and calibrates $N\#\ssc$ with respect to the metric $\hgd$ in $\hat{U}\left(\sing\left( N\#\ssc\right)\right)$ and $B_\de(\widecheck{\reg}(N\#\ssc))\setminus \ug$.
	
	On $\ov{\ud},$ we have $$\hgd=\hkn\oplus \big((1+K\de_K)\hkp\big)\ge \hkn\oplus \hkp=h_\sing,$$ as quadratic forms and $\phi=\phi_\sing.$ By bullet (\ref{cms1}) in Fact \ref{cmsvec}, the form $\phi$ remains a calibration on $\ov{\ud}.$ When restricted to $\ov{\ud}\cap\reg \left(N\#\ssc\right),$ the metrics $\hgd$  and $h_\sing$ are equal to each other. Thus, we deduce that in $\ov{\ud}$, the form $\phi=\phi_\sing$ still calibrates $N\#\ssc$  with respect to $\hgd.$
	
	On ${\uc}$, first note that $\frac{1}{2}\le\phi_\sing(\ker^\perp\ed\pi)\le 1$ by (\ref{philbd}) and Fact \ref{calsing}. Thus,  on ${\uc}$ we have
	\begin{align*}
		\hgd=&	\bigg((1-\de_\perp)+\de_\perp(\phi_\sing(\kn))\bigg)^{\frac{2}{d}}\Big(\hkn\oplus \big((1+K)\hkp\big)\Big)\\\ge& \big(\phi_\sing(\ker^\perp\ed\pi)^{\frac{2}{d}}\big)\big(\hkn\oplus((1+K)\hkp)\big)\\
		\ge& \big(\phi_\sing(\ker^\perp\ed\pi)^{\frac{2}{d}}\big)\big(\hkn\oplus(K\hkp)\big),
	\end{align*}as quadratic forms.
	By bullet (\ref{cms1}) in Fact \ref{cmsvec}, inequality (\ref{kge}) of Fact \ref{calk}, and (\ref{defk}), we deduce that on ${\uc}$, the form $\phi$ remains a calibration in $\hgd.$ Again, the restrictions of $h_\sing$ and $\hgd$ to $\uc\cap\reg \left(N\#\ssc\right)$ are equal, so $N\#\ssc$ is calibrated by $\phi=\phi_\sing$ in ${\uc}$. 
	
	On $\ov{\ub},$ we have
	\begin{align*}
		\hgd=&(\phi_\lam(\kn)^{\frac{2}{d}})\Big(\hkn\oplus \big((1+K)\hkp\big)\Big)\\\ge& \big(\phi_\lam(\ker^\perp\ed\pi)^{\frac{2}{d}}\big)\big(\hkn\oplus(K\hkp)\big),
	\end{align*}as quadratic forms.
	By bullet (\ref{cms1}) in Fact \ref{cmsvec}, inequality (\ref{kge}) of Fact \ref{calk}, and (\ref{defk}), we deduce that on $\ov{\ub}$, the form $\phi$ is a calibration form in $\hgd.$ Also when restricted to $\ov{\ub}\cap\reg \left(N\#\ssc\right)$, we have $\hgd=h_\reg$, and $\phi=\phi_\reg$. This implies that on $\ov{\ub},$ the form $\phi$ calibrates $N\#\ssc$ with respect to $\hgd.$
	
	Finally, on $\ua$, we have \begin{align*}
		\hgd=&	(\phi_\reg(\kn)^{\frac{2}{d}})\Big(\hkn\oplus \big((1+\de_0 K)\hkp\big)\Big)\\\ge&(\phi_\reg(\kn)^{\frac{2}{d}})\big(\hkn\oplus\hkp\big)\\=& h_\reg.
	\end{align*} as quadratic forms. By bullet (\ref{cms1}) of Fact \ref{cmsvec}, this means that $\phi=\phi_\reg$ is a calibration in ${\ua}$ with respect to $\hgd.$ Again, the restrictions of $\hgd$ and $h_\reg$ to $\ua\cap \reg \left(N\#\ssc\right)$ are equal, so we deduce that $N\#\ssc$ is calibrated by $\phi=\phi_\reg$ in $\ua$.  
	
	Now extend $\hgd|_{U_\e(N\#\ssc)}$  to a smooth Riemannian metric defined throughout $M.$ We will still use $\hgd$ to denote this extension of $\hgd|_{U_\e(N\#\ssc)}$.
	\subsubsection{Adding a bump to make $N\#\ssc$ uniquely calibrated}
	Use Lemma \ref{lemzs} to construct a smooth function  $\be$ on $M$ with value between $[0,\frac{1}{2}]$ and \begin{align}\label{bez}
		\be\m(0)=\supp N\#\ssc\cup \ov{U\left(\sing\left( N\#\ssc\right)\right)}.
	\end{align}
	
	Then in the metric $(1+\be)\hgd$, by bullet (\ref{cms1}) of Fact \ref{cmsvec}, $\phi$ is still a calibration form. By (\ref{bez}), $\phi$ still calibrates $N\#\ssc$ in $(1+\be)\hgd$. 
	
	Now suppose $T$ is a homologically non-trivial integral cycle on $M$ that is calibrated by $\phi$ with respect to $(1+\be)\hgd$ in $U_\e(N\#\ssc)$. Since on the complement of $\supp\left(N\#\ssc\right)\cup \ov{U\left(\sing\left( N\#\ssc\right)\right)}$ in $M,$ $\phi$ has comass strictly smaller than $1$ in $(1+\be)\hgd$ by (\ref{bez}) and bullet (\ref{cms2}) in Fact \ref{cmsvec}, we deduce that $T$ must be supported in $\supp\left(N\#\ssc\right)\cup \ov{U\left(\sing\left( N\#\ssc\right)\right)}$.
	
	Note that $T$ cannot be supported in $\ov{U\left(\sing\left( N\#\ssc\right)\right)}$, since $$H_d(U\left(\sing\left( N\#\ssc\right)\right),\Z)=0$$ yet $[T]\not=0\in H_d(M,\Z).$
	
	Now by Fact \ref{const}, we deduce that $T$ restricted to the non-empty closed set $\supp N\#\ssc\cap \ov{U\left(\sing\left( N\#\ssc\right)\right)}\cp$ equals an integer multiple $k$ of $N\#\ssc$ restricted to $\supp N\#\ssc\cap \ov{U\left(\sing\left( N\#\ssc\right)\right)}\cp,$  with $k\not=0.$ By Lemma \ref{suct} and the fact that having infinite order tangent to each other is a relatively closed condition, we deduce that $T=k N\#\ssc.$ Since $T$ is calibrated by $\phi,$ from $$0<\ms_{(1+\be)\hgd}(T)=\phi(T)=k\phi( N\#\ssc),$$ we deduce that $k>0.$ 
	
	To sum it up, we have proved that if $T$ is calibrated by $\phi$, then $T$ is a positive integer multiple of $N\#\ssc.$
	
	From $H_d(U_\e(N\#\ssc),\Z)=\Z[N\#\ssc]$, we can conclude that:
	\begin{fact}\label{uni}
		For integer $k\not=0,$ the integral current $kN\#\ssc$ is the unique integral current calibrated by $\frac{k}{|k|}\phi $ in the homology class $k[N\#\ssc]$ in the region $U_\e(N\#\ssc)$ with respect to the metric $(1+\be)\hgd.$
	\end{fact}
	\subsubsection{Making the metric large away from $\supp N\#\ssc$}This subsubsection is taken from the proof of \cite[Theorem 4.1]{YZj}.
	By the monotonicity formula for stationary varifolds \cite[Section 5.1]{WA}, there is a constant $A$ such that for any stationary integral current $T$ on $M$ with respect to the metric $(1+\be)\hgd$, we have
	\begin{align}\label{mon}
		\ms_{(1+\be)\hgd}(T|_{B_r(p)})\ge Ar^d,
	\end{align}
	provided $p\in \supp T$ and $r<\operatorname{InjRad} M$. Here by $B_r(p)$ we mean the ball of radius $r$ centered at $p$  in the metric $(1+\be)\hgd$. The symbol  $\operatorname{InjRad} M$ denotes the injectivity radius of $M$. By stationary current we mean that the mass, i.e., area, of $T$ is a critical point with respect to ambient diffeomorphism perturbations that fix $\pd T$. \newcommand{\ue}{U_{\e'}(N\#\ssc)}\newcommand{\uee}{U_{\e''}(N\#\ssc)}
	
	By Lemma \ref{tube}, there exists $0<\e''
	<\e'<\e<\operatorname{InjRad}M$ such that \begin{align*}
		\uee\s \ue\s U_\e(N\#\ssc),
	\end{align*}
	with the distance between $\pd U_{\e'}(N\#\ssc)$ and $\pd U_\e(N\#\ssc)$ at least $\e''$, and the distance between $\pd\uee$ and $\pd  \ue$ at least $\e''.$
	
	Then let $\ga$ be a smooth function that equals to  $1$ on $\uee$ and equals to 
	\begin{align}\label{mxm}
		\max\Bigg\{1,\frac{\sqrt[d]{2\ms_{(1+\be)\hgd}(N\#\ssc)}}{\sqrt[d]{A} \e''}\Bigg\},
	\end{align}
	on the complement of $\ue$, and takes value between $1$ and (\ref{mxm}) on $\ue\setminus \uee.$ 
	\begin{defn}\label{defnh}
		Define
		\begin{align*}
			h=\ga^2(1+\be)\hgd.
		\end{align*}
	\end{defn}
	Then by bullet (\ref{cms1}) in Fact \ref{cmsvec}, $\phi$ remains a calibration in $h$ and still calibrates $N\#\ssc$. 
	
	Now let $T$ be an area-minimizing integral current in $[\Si]$ with respect to metric $h$ on $M.$ There are two possibilities, either $\supp T$ is contained in $U_\e(N\#\ssc)$, or $\supp T$ contains a point in the complement of $U_\e(N\#\ssc)$. The latter cannot happen. To see this, note that $h$ is a constant multiple of $\hgd$ in the complement of $\ue,$ so $T$ restricted to $U_{\e'}(N\#\ssc)\cp$ is also a stationary integral current with respect to $\hgd$. Thus, one can apply (\ref{mon}) with $p\in U_\e(N\#\ssc)\cp\cap\supp T$ and $r=\e''$ so that $B_r(p)\s U_{\e'}(N\#\ssc)\cp$ to conclude that
	\begin{align*}
		&\ms_{h}(T)\\\ge& \ms_{h}(T|_{B_{\e''}(p)})=\max\Bigg\{1,\frac{\sqrt[d]{2\ms_{(1+\be)\hgd}(N\#\ssc)}}{\sqrt[d]{A} \e''}\Bigg\}^d\ms_{(1+\be)\hgd}(T|_{B_{\e''}(p)})\\\ge&\max\Bigg\{1,\frac{{2\ms_{(1+\be)\hgd}(N\#\ssc)}}{{A} (\e'')^d}\Bigg\} A(\e'')^d\\
		=&\max\{ A(\e'')^d,2\ms_{(1+\be)\hgd}(N\#\ssc)\}\\
		=&\max\{ A(\e'')^d,2\ms_h(N\#\ssc)\}.
	\end{align*}
	In other words, $T$ has area at least two times that of $N\#\ssc$ in metric $h$, so $T$ cannot be area-minimizing.
	
	 To sum it up, we have the proven the following:
	\begin{fact}\label{mint}
		In the smooth Riemannian metric $h$ (Definition \ref{defnh}) on $M,$ the smooth form $\phi$ (\ref{defp}) on $U_\e(B\#\ssc)$ calibrates $N\#\ssc.$ Moreover, any other area-minimizing integral current $T$ in $[\Si]$ must have $\supp T\s U_\e(B\#\ssc).$
	\end{fact}
	\subsection{Wrapping up proof of Lemma \ref{nsmin}}
	Let $T$ be an area-minimizing integral current in $[\Si].$ By Fact \ref{mint}, $T$ is supported in $U_\e(N\#\ssc).$ However, since $U_\e(N\#\ssc)$ is a subset of $M,$ we need to redetermine the relation between $[T]$ and $[N\#\ssc]$ when restricted to $U_\e(N\#\ssc).$
	
	Since we have $H_d(U_\e(N\#\ssc))=\Z[N\#\ssc]$ by Lemma \ref{tube}, we deduce that $[T]= k[N\#\ssc]$ as homology classes in $U_\e(N\#\ssc)$ with $k\not=0.$ Thus, we have
	$$|k|\ms_h(N\#\ssc)=|kN\#\ssc(\phi)|=|T(\phi)|\le \ms_h(T).$$ Since $T$ is area-minimizing in $[\Si]$, we have
	\begin{align*}
		\ms_h(N\#\ssc)\ge \ms_h(T).
	\end{align*}
	Thus we must have $k=\pm 1$ and $\ms_h(N\#\ssc)=\ms_h(T).$ This implies that $T$ is calibrated by $\frac{k}{|k|}\phi.$ Now by Fact \ref{uni} and $h\ge \hgd$, we deduce that $T=\pm N\#\ssc$ as claimed. 
	
	The bullet in our prompt follows because we have never changed the metric $h_\sing$ inside $$\hat{U}\left(\sing\left( N\#\ssc\right)\right)\cap {\uee}.$$ Thus, for some $\rr>0$ small, in a neighborhood ${U}_\rr(\sing\ssc)$ diffeomorphic to $B_{\rr}^{d-s+c}\times S^s$ via $\Ga$ the metric $h$ is the pullback of the standard product metric on $B_{\rr}^{d-s+c}\times S^s$. 
	\begin{defn}\label{defnur}
		For $0<\eta\le\rr,$ we define  $U_\eta(\sing\ssc)$, the radius $\eta$ tubular neighborhood of $\sing\ssc,$ by
		\begin{align*}
			U_\eta(\sing\ssc)=\Ga\m(B_{\eta}^{d-s+c}\times S^s).
		\end{align*}
	\end{defn}
	\section{Orthogonal transverse immersions have persistent singular sets}\label{lawlor}
	Now that we have constructed the area-minimizing representative $N\#\ssc$ in Lemma \ref{nsmin}, we need to choose suitable $C$ and $s$ to create persistent singular sets. As explained in our plan of proof Section \ref{planpf}, we have two kinds of singular sets that we want to add, transverse intersections from Fact \ref{fctgl} and products of $C(\cpt)$ with spheres from Fact \ref{fctcb}. Transverse intersections from Fact \ref{fctgl} contributes a singular set of dimension $(d-c)$, and products of $C(\cpt)$ with spheres from Fact \ref{fctcb} contributes a singular set of dimension $(d-5)$.
	
	In this section we will focus on adding the singular sets in Fact \ref{fctgl}. 
	
	Now take two copies of $\R^c,$ with standard coordinate systems  $x=(x_1,\dots,x_c),$ and $y=(y_1,\dots,y_c),$ respectively. Then  the product structure of $$\R^{2c}=\R^c\times \R^c$$ gives us a coordinate system
	\begin{align*}
		(x,y)=(x_1,\dots,x_c,y_1,\dots,y_c)
	\end{align*}on $\R^{2c}.$
	
	Set $\R^c\times\{0\}^c$ with parameterization $(x,0)$ and $\{0\}^c\times\R^c$ with parameterization $(0,y).$ In this section we always use
	\begin{defn}\label{imc}(\textbf{In this section only})
		\begin{align*}
			C=&\,\R^c\times\{0\}^c+\{0\}^c\times \R^c,\\
			\phi_C=&\,dx_1\w\dots\w dx_c+dy_1\w\dots\w dy_c,\\
			s=&\,d-c.
		\end{align*}	
	\end{defn}
	By \cite[Chapter II, Proposition 7.12]{HL}, the form $\phi_C$ is a calibration form that calibrates $C$, provided $c\ge 2.$
	
	Our goal in this section is to prove the following lemma:
	\begin{lem}\label{nst}
		The subset  $\sing\sdc$ of the the singular set of the area-minimizing integral current $N\#\si^{d-c}(C)$ obtained in Lemma \ref{nsmin} is a persistent singular subset (Definition \ref{psing}), provided  $3\le c\le d.$
	\end{lem}
	Before delving into the proof, it is beneficial to review Lawlor's example in \cite[Section 6.4]{GL}:
	\begin{exam}\label{l3e}
		The area-minimizing cone $C=\R^3\times\{0\}^3+\{0\}^3\times \R^3$ restricted to the unit ball $B_1^6$ of $\R^6$ has a persistent singular set among Euclidean metric perturbations.
		
		More precisely, let $\Psi$ be a smooth self-diffeomorphism of $\R^6.$ Then area-minimizing integral currents $T$ with boundary $\pd T=\Psi\pf(\pd (C|_{B_1^6}))$  decompose into the sum of two minimal $3$-balls intersecting transversely near the origin, provided $\Psi$ is close to identity map.
	\end{exam}
	Note that above result is totally false if we replace $C$ with $\R^2\times\{0\}^2+\{0\}^2\times \R^2,$ since White has proved that $2$ dimensional area-minimizing integral currents are generically smooth in \cite{BWgt}.
	
	Lawlor's idea is to construct explicitly an integral current with boundary $\Psi\pf(\pd (C|_{B_1^6}))$ that is the sum of two minimal $3$-balls $X$ and $Y$ and prove that the sum $X+Y$ is the unique area-minimizing integral current with boundary $\Psi\pf(\pd(C|_{B_1^6}))$ using an ingenious calibration argument.
	\begin{figure}[h]
		\centering
		\def\svgwidth{0.8\paperwidth}
		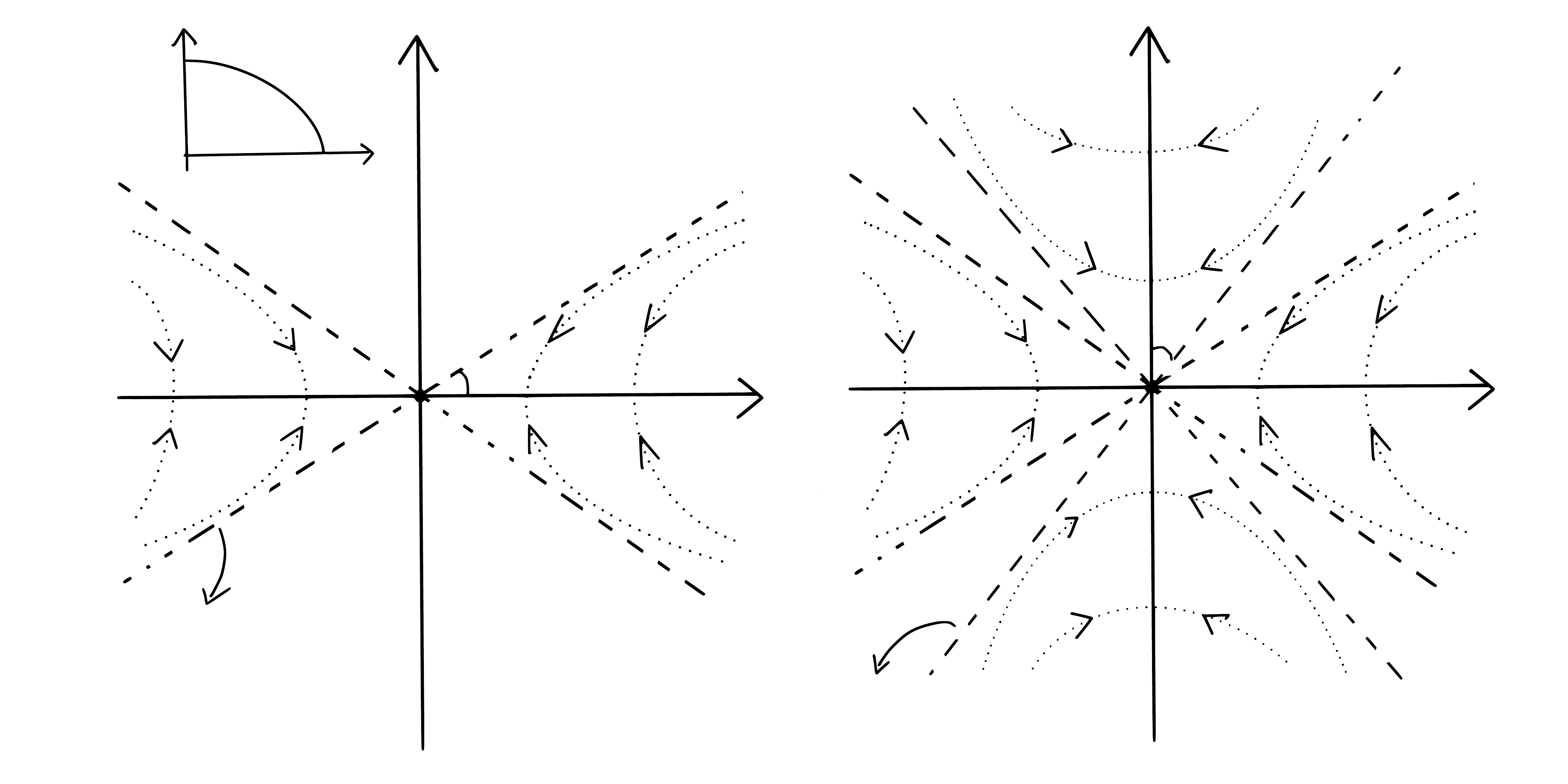
		\caption{The retractions and the two wedges}
\label{figlaw}
	\end{figure}
	To explain his ideas, let us start with Lawlor's calibration of $C=\R^3\times\{0\}^3+\{0\}^3\times \R^3.$ 
	First, let us define a Lipschitz retraction $$\Pi_{\R^3\times\{0\}^3}:\R^6\to \R^3\times\{0\}^3,$$ 
	defined by
	\begin{align}\label{poe}
		{\po}(x,y)=\bigg(\bigg(\cg\bigg(\frac{|y|}{|x|}\bigg)\bigg)^{\frac{1}{3}}x,0\bigg).
	\end{align}
	Here $\cg:[0,\infty)\to \R$ is a non-negative Lipschitz function satisfying
	\begin{align}\label{lodee}
		\bigg(\cg(t)-\frac{t}{3}\cg'(t)\bigg)^2+\bigg(\frac{\cg'(t)}{3}\bigg)^2\le1-\frac{1}{16} t^2,
	\end{align} with  $\cg(0)=1$, $\cg(t)=0$ for $t\ge t_3$ and $\cg$ strictly decreasing on $[0,t_3],$ with
	\begin{align}\label{t0e}
		0<t_3<1.
	\end{align}
	In other words $\po$ is a Lipschitz retraction of $\R^6$ onto $\R^3\times\{0\}^3$, such that
	\begin{align}\label{vpe}
		{\po}(x,y)\equiv(0,0),\text{ if }\frac{|y|}{|x|}\ge t_3.
	\end{align} 
	\begin{defn}\label{defnwe}
		Define $W_{\ta}(\R^3\times\{0\}^3,\{0\}^6)$, the angle $\ta$ wedge neighborhood of $\R^3\times\{0\}^3$ with vertex $\{0\}^6$ as 
		\begin{align}\label{woe}
			W_{\ta}(\R^3\times\{0\}^3,\{0\}^6)=\bigg\{(x,y)\bigg|\frac{|y|}{|x|}\le \tan\ta\bigg\}.
		\end{align}
	\end{defn}
	Then (\ref{vpe}) says that $\po$ is non-vanishing only in $$W_{\arctan t_3}(\R^3\times\{0\}^3,\{0\}^6).$$
	
	The ODE (\ref{lodee}) guarantees that $\po$ is a $3$\textbf{-area-non-increasing} retraction onto $\R^3\times\{0\}$, which is defined as follows: 
	\begin{defn}\label{dani}
		We say a Lipschitz map $\Pi$ is $d$-area-non-increasing at $p$, if the differential $\ed\Pi$ of $\Pi$ at $p$ maps unit simple $d$-vectors to simple $d$-vectors of length at most $1.$ If $\Pi$ is $d$-area-non-increasing at every point where it is differentiable, then we say $\Pi$ is $d$-area-non-increasing.
	\end{defn}
	A direct consequence of area-non-increasing property is as follows:
	\begin{fact}\label{fctani}
		A Lipschitz retraction $\Pi$ onto a $d$-dimensional submanifold $N$ is $d$-area-non-increasing at $p$ if and only if the $d$-form $\Pi\du(\dvol_N)$ has comass at most $1$ at $p,$ i.e.,
		\begin{align}\label{anicms}
			\cms\big(\Pi\du(\dvol_N)\big)_p\le 1.
		\end{align}Furthermore, equality in (\ref{anicms}) holds if and only if $\ed\Pi$ maps the unit simple $d$-vector orthogonal to $\ker \ed\Pi$ at $p$ to a unit simple $d$-vector at $\Pi(p).$
	\end{fact}
	Indeed we have a more general fact.
	\begin{fact}\label{anical}
		If we have a smooth calibration form $\psi$ defined on a $d$-dimensional submanifold $N$, and a smooth $d$-area-non-increasing map $\Pi$ onto $N,$ then $\Pi\du\Psi$ is a calibration form.
	\end{fact}
	The above facts follow directly from the definition of comass, the simpleness of $\dvol_N$ and the orthogonal splitting of tangent spaces of the ambient manifold into $\ker \ed\Pi\oplus \ker^\perp\ed\Pi.$
	
	From Fact \ref{fctani}, we deduce that
	\begin{align*}
		\po\du(\dvol_{\R^3\times\{0\}}),
	\end{align*} is a degree $3$ calibration form with Lipschitz anti-derivative (Definition \ref{defnlip}), since $\dvol_{\R^3\times\{0\}^3}$ is exact. Direct calculation shows that $\po\du(\dvol_{\R^3\times\{0\}})$ is almost continuous (Definition \ref{defnlip}) on $\R^3\times\{0\}^3$ and calibrates only $\R^3\times\{0\}^3$.
	
	Now in the same vein, we can define
	\begin{align}\label{pte}
		\pt(x,y)=\bigg(0,\cg\bigg(\frac{|y|}{|x|}^{\frac{1}{3}}\bigg)y\bigg),
	\end{align}so that
	\begin{align}\label{vpte}
		{\pt}(x,y)\equiv(0,0),\text{ if }\frac{|x|}{|y|}\ge t_3.
	\end{align}
	We can similarly define $W_{\ta}(\{0\}^3\times\R^3,\{0\}^6)$,the angle $\ta$ wedge neighborhood of $\{0\}^3\times\R^3$ with vertex $\{0\}^6$ as 
	\begin{align*}
		W_{\ta}(\{0\}^3\times\R^3,\{0\}^6)=\bigg\{(x,y)\bigg|\frac{|x|}{|y|}\le \tan\ta\bigg\}.
	\end{align*}
	Then $\pt$ is non-vanishing only in $W_{\arctan t_3}(\{0\}^3\times\R^3,\{0\}^6).$
	
	Again $\pt\du(\dvol_{\{0\}^3\times\R^3})$ is a calibration form with Lipschitz anti-derivative and calibrates only $\{0\}^3\times \R^3.$
	
	Since the two wedges $W_{\arctan t_3}(\R^3\times\{0\}^3,\{0\}^6)$ and $W_{\arctan t_3}(\{0\}^3\times\R^3,\{0\}^6)$ only intersect at the origin due to (\ref{t0e}), by (\ref{vpe}) and (\ref{vpte}) we deduce that $\po$ and $\pt$ patch together Lipschitzly through the region
	\begin{align*}
		t_3<\frac{|y|}{|x|}<\frac{1}{t_3}
	\end{align*}into a Lipshictz $3$-area-non-increasing retraction:
	\begin{align*}
		\Pi_C:\R^6\to C.
	\end{align*}
	The reader should now consult Figure \ref{figlaw} for intuition. Consequently,
	\begin{align*}
		\Pi_C\du\dvol_C=\po\du(\dvol_{\R^3\times\{0\}^3})+\pt\du(\dvol_{\{0\}^3\times\R^3}),
	\end{align*}
	is a calibration form with Lipschitz anti-derivative that calibrates only $C.$ 
	
	At this point, all of these seem fairly uninteresting and we have just found a very artificial and complicated way to calibrate $C.$ However, a little twist to the above construction immediately proves the statements in Example \ref{l3e}. 
	
	Lawlor's ingenious idea is as follows. Note that $\pd(C|_{B_1^6})$ is the sum of two disjoint pieces 
	\begin{align*}
		S^2\times\{0\}^3+\{0\}^3\times S^2
	\end{align*}
	By the implicit function theorem of partial differential equations, $\Psi\pf(S^2\times \{0\}^3)$ bounds a minimally embedded $3$-ball $X$ and $\Psi\pf(\{0\}^3\times S^2)$ bounds a minimally embedded $3$-ball $Y,$ intersecting at a point $$X\cap Y=p$$ near the origin, provided $\Psi$ is close to the identity map.
	
	Now define a retraction $\Pi_X$ onto $X$ as follows in the spirit of (\ref{poe}) . For any point $q\in \R^6,$ Use $\pi_X(q)$ to denote the nearest distance projection of $q$ to $X$. On length minimizing geodesic $\ga\s X$ from $p$ to $\pi_X(q)$ on $X$ there is a unique point, $$\cg^{\frac{1}{3}}\bigg(\frac{\fd_{\R^6}(q,\pi_X(q))}{\fd_X(\pi_X(q),p)}\bigg)\pi_X(q),$$that is of distance precisely 
	\begin{align*}
		\cg^{\frac{1}{3}}\bigg(\frac{\fd_{\R^6}(q,\pi_X(q))}{\fd_X(\pi_X(q),p)}\bigg)\fd_X(\pi_X(q),p)
	\end{align*}
	to $p.$ Here $\fd_{\R^6}$ is the distance in $\R^6$ and $\fd_X$ is the Riemannian distance function on $X.$ Now define
	\begin{align*}
		\Pi_X(q)=\cg^{\frac{1}{3}}\bigg(\frac{\fd_{\R^6}(q,\pi_X(q))}{\fd_X(\pi_X(q),p)}\bigg)\pi_X(q).
	\end{align*}
	This gist of this definition is that in a Fermi coordinate chart centered at $p$ and adopted to $X,$ 
	\begin{align*}
		\Pi_X=\po.
	\end{align*}
	Then again using (\ref{lodee}), one can prove that $\Pi_X$ is a $3$-area-non-increasing retraction onto $X$, provided $\Psi$ is close to identity. Thus, we can prove that
$		\Pi_X^*\dvol_X
$ is a calibration form with Lipschitz anti-derivative and  calibrates only $X.$
	
	Similarly, one define $\Pi_Y$ as
	\begin{align*}
		\Pi_Y(q)=\cg^{\frac{1}{3}}\bigg(\frac{\fd_{\R^6}(q,\pi_Y(q))}{\fd_Y(\pi_Y(q),p)}\bigg)\pi_Y(q),
	\end{align*} and prove that $\Pi_Y$ is $3$-area-non-increasing. Then 
$		\Pi_Y^*\dvol_Y
$ is a calibration form with Lipschitz anti-derivative and calibrates only $Y.$
	
	Then again, one can define wedges $W_{\arctan t_3}(X,p)$ and $W_{\arctan t_3}(Y,p)$ similar to Definition \ref{defnwe}, such that $\Pi_X$ is not of constant value $p$ only on $W_{\arctan t_3}(X,p)$ and $\Pi_Y$ is not of constant value $p$ only on $W_{\arctan t_3}(Y,p)$. Now the two wedges $W_{\arctan t_3}(X,p)$ and $W_{\arctan t_3}(Y,p)$ do not intersect, provided $\Psi$ is close to identity. This implies that $\Pi_X$ patches Lipschitzly onto $\Pi_Y$ to give a Lipschitz $3$-area-non-increasing retraction onto $X+Y$, and the calibration form
	\begin{align*}
		\Pi_X\du \dvol_X+\Pi_Y\du\dvol_Y
	\end{align*} only calibrates $X+Y.$ This 
	implies that $X+Y$ is the unique area-minimizing integral current with boundary $\Psi\pf(\pd(C|_{B_1^6}))$ and we get the claims in our Example \ref{l3e}. 
	
	At this point the reader might wonder why the above argument fail in the case of $2$-dimensional surfaces, since $2$-dimensional area-minimizing integral currents are generally smooth by \cite{BWgt}. The main reason is that the ODE (\ref{lodee}) does not have a solution $\cg$ with $\cg(0)=1$ and $\cg(1)=0$ when the number $3$ in the denominators of (\ref{lodee}) is replaced by the number $2$. Thus, we cannot patch together $\Pi_X$ and $\Pi_Y$ as above.

	Our proof of Lemma \ref{nst} is precisely a Riemannian version of Lawlor's arguments. Let $g$ be a smooth metric close enough to $h$ obtained in Lemma \ref{nsmin}. (The closeness of $g$ to $h$ will be determined in the proof.) Our argument will involve five steps: 
	\begin{enumerate}
		\item With respect to the metric $g$, determine the structure of area-minimizing integral currents $T$ in $[\Si]$, when restricted to the complement of $U_\eta(\sing \sdc)$ (Definition \ref{defnur}).\label{im1}
		\item Using the above step, we construct embedded minimal $d$-dimensional submanifolds $X$ and $Y$ contained in $U_\eta(\sing \sdc),$ with transverse intersection $S^s$ such that \begin{align*}
			\pd(X+Y)=\pd\bigg( T|_{U_\eta(\sing \sdc)}\bigg).
		\end{align*}\label{im2}
		\item Construct the $d$-area-non-increasing retractions $\Pi_X$ and $\Pi_Y$ onto $X$ and $Y$.
		\item Prove that  $\Pi_X\du\dvol_X+\Pi_Y\du\dvol_Y$ is a calibration form only calibrating $X+Y.$\label{im3}
		\item Use $H_d(U_\eta(\sing\sdc),\Z)=0$ to conclude that
		\begin{align*}
			X+Y=T|_{U_\eta(\sing \sdc)}.
		\end{align*}\label{im4}
	\end{enumerate}
	The Step (\ref{im3}) is the hardest and will take up the majority of our efforts. The main difficulty is that we no longer have a Euclidean metric, so that we have to estimate carefully the effects of curvature on our argument. Our proof in this section will start with some preparational work and then go through the above five steps.
	
	\subsection{Lawlor's ODE}
	Lawlor's ODE in (\ref{lodee}) will be a crucial ingredient of our proof. In this subsection, we will formally define the Lipschitz functions $\cg.$
	\begin{defn}\label{defng}
		For every integer $c\ge 3,$ define a  function $\cg_c(t):[0,\infty)\to \R$ as follows:
		\begin{align*}\cg_c(t)=
			\begin{cases}
				1-\frac{2 c^2-9}{8}  t^2, &\text{for }t\in\Big[0,\frac{2 \sqrt{2}}{\sqrt{2 c^2-9}}\Big),\\
				0,&\text{for }t\in \Big[\frac{2 \sqrt{2}}{\sqrt{2 c^2-9}},\infty\Big).
			\end{cases}
		\end{align*} 
		We will also denote
		\begin{align}
		\label{defntc}	t_c&=\frac{2 \sqrt{2}}{\sqrt{2 c^2-9}},\\\label{defncc}\cc_c(t)&={\cg_c}(t)-\frac{t
			}{c}{\cg_c}'(t),\\\label{defncs}\cs_c(t)&=\frac{{\cg_c}'(t)}{c}.
		\end{align}
	\end{defn}
	\begin{fact}\label{fctcg}
		The  function $\cg_c$ is a monotone Lipschitz function on $[0,\infty),$ and is strictly monotonically decreasing when restricted to $[0,t_c]$, with $\cg_c(0)=1,\cg_c(t_c)=0$. Moreover, 
		\begin{align}\label{asymptc}
			1>t_c=\frac{2}{c}+O(c^{-2}).
		\end{align} for all $c\ge 3.$
	\end{fact}The above fact is obvious. The function $\cg_c$ satisfies a more general version of Lawlor's ODE as follows:
	\begin{lem}
		On the interval $[0,\ct_c),$ the Lipschitz function satisfies Lawlor's ODE:
		\begin{align}\label{lode}
			\frac{3}{4}\le \cc_c(t)^2+\cs_c(t)^2= \left({\cg_c}(t)-\frac{t
			}{c}{\cg_c}'(t)\right)^2+ \left(\frac{{\cg_c}'(t)}{c}\right)^2\le 1-\frac{1}{16}t^2,
		\end{align}and
		\begin{align}\label{csm}
			\cc_c(t)^2+\cs_c(t)^2=1\iff t=0.
		\end{align}
	\end{lem}
	\begin{proof}
		The proof is a direct calculation. A Mathematica verification will be provided at the end of the manuscript after bibliography. We will give a sketch here. We have
		\begin{align}\label{eqcs}
			\cc_c(t)^2+\cs_c(t)^2=a_c\cdot(t^2)^2+b_c\cdot(t^2)+1,
		\end{align}with
		\begin{align}\label{acsign}
			a_c=\frac{c^4}{16}-\frac{c^3}{4}-\frac{5
				c^2}{16}+\frac{81}{16 c^2}+\frac{9 c}{4}-\frac{81}{16
				c}-\frac{63}{64}=\frac{(c-2)^2 \left(2 c^2-9\right)^2}{64 c^2}>0,
		\end{align}
		and
		\begin{align}
			b_c=-\frac{c^2}{4}+\frac{81}{16 c^2}+c-\frac{9}{2 c}=-\frac{\left(2 c^2-9\right) \left(2 c^2-8 c+9\right)}{16 c^2}<0.
		\end{align}provided $c\ge 3.$
		Now treat the right hand side of equation (\ref{eqcs}) as a quadratic polynomial in variable $\tau=t^2,$ then the discriminant of the right hand side of (\ref{eqcs}) equals
		$$
			b_c^2-4a_c=-\frac{\left(2 c^2-9\right)^3 \left(6 c^2-16 c+9\right)}{256
				c^4}<0,
		$$
		for $c\ge 3.$ Negative discriminant implies $	a_c\cdot(t^2)^2+b_c\cdot(t^2)+1$ has no root, so $a_c>0$ implies that $	a_c\cdot(t^2)^2+b_c\cdot(t^2)+1$ is always positive. The minimum of 	$a_c(t^2)^2+b_c(t^2)+1$ is achieved at $t^2=\frac{-b_c}{2a_c}$. We have
		\begin{align*}	
			&	a_c\cdot(t^2)^2+b_c\cdot(t^2)+1\\
			\ge& a_c\left(\frac{-b_c}{2a_c}\right)^2+b_c\left(\frac{-b_c}{2a_c}\right)+1\\
			=&\frac{\left(2 c^2-9\right) \left(6 c^2-16 c+9\right)}{16 (c-2)^2
				c^2}\\
			=&\frac{3}{4}+\frac{(-3 + 2 c)^2 (-9 + 4 c)}{\left(2 c^2-9\right) \left(6 c^2-16 c+9\right)}\\
			>&\frac{3}{4}
		\end{align*}for  all $t\in\R.$ This finishes the proof of the left hand side of (\ref{lode}).
		
		On the other hand, when we restrict to $t\in[0,t_c),$ we always have $t^4<t_c^2t^2.$ Since $a_c>0$ due to (\ref{acsign}), we have
		\begin{align*}
			&	a_c\cdot(t^2)^2+b_c\cdot(t^2)+1\\\le& a_ct_c^2t^2+b_ct^2+1\\=& (a_ct_c^2+b_c)t^2+1\\
			=&\left(\frac{9}{16 c^2}-\frac{1}{8}\right) t^2+1\\
			\le&1-\frac{1}{16}t^2,
		\end{align*}provided $c\ge 3.$ This finishes the proof of the right hand side of (\ref{lode}).
		
		The claim (\ref{csm}) follows from the right hand side of (\ref{lode}) and $\cc_c(0)^2+\cs_c(0)^2=1.$

The referee has suggested a more human verifiable approach as follows. Setting $\mathbf{t}=\frac{2c^2-9}{4}t^2,$ then we have $\mathbf{t}\in[0,2)$. From $\cc_c^2(t)=\left(1-\left(\frac{1}{2}-\frac{1}{c}\right)\mathbf{t}\right)^2,\cs_c(t)^2=\frac{2c^2-9}{4c^2}\mathbf{t}$, the proof can be done as estimating the range of a quadratic function in $\mathbf{t}$ with domain $[0,2)$.
	\end{proof}
	\subsection{Geometric dependence}
	As Lemma \ref{nst} is essentially a statement about an open subset in the space of Riemannian metrics, we will frequently encounter quantitative estimates or qualitative statements that need a uniform bound depending on the metric. To this end, we give a formal definition as follows.
	\begin{defn}\label{defngdp}
		Let $X,Y,Z,$ etc., be submanifolds of $M$.	Suppose we have an object $A$, that is constructed from $X,Y,Z,$ etc., and the ambient metric $h.$ We say the quantity $A$  depends geometrically on $X,Y,Z,$ etc. and $h$, if the object $A$ varies continuously in its category if we vary $X,Y,Z,$ etc., and the ambient metric $h$ continuously in the \textbf{smooth topology}.
	\end{defn}
	Here by varying $X$, etc., we mean the following.
	\begin{defn}\label{defngm}
		Let $X$ be a compact smooth submanifold of $M.$ We regard $X$ as a point in the space $C^\infty(X,M)$ of smooth maps from $X$ to $M.$ Then we say a submanifold $X'$ is $\e$-close to $X$ as a submanifold, if $X'$ can be realized as a map in $C^\infty(X,M)$ of distance at most $\e$ away to $X$ in the standard metric on $C^\infty(X,M)$. \end{defn}
	Here the standard metric on $C^\infty(X,M)$ is the one obtained from the countably many $C^k$ semi-norms for $k\in\N$, e.g., \cite[Theorem A.4.7]{CW}.
	
	For instance, the second fundamental form tensor of a submanifold $N$ depends geometrically on $N$ and the ambient metric $h.$ The injectivity radius of $M$ depends geometrically on the ambient metric.
	
	In most cases, geometric dependence of $A$ on $X,Y,Z,$ etc., and $h$ is a straightforward calculation and we will not give a detailed proof.
	
	An easy consequence of geometric dependence is as follows. If a real number $A$ is positive on some $X,Y,Z$ etc., and $h,$ and depends geometrically on $X,Y,Z,$ etc., and $h,$ then there exists a neighborhood $\Om_h$ in the space of Riemannian metrics, and $\e>0,\de>0,$ such that whenever $X',Y',Z'$ etc., is $\de$ close to the $X,Y,Z$ in smooth topology and $h'\in \Om_h$ then $A\ge \e.$ 
	\subsection{Normal bundles and neighborhoods}
	As mentioned when introducing Lawlor's ideas, we need several geometric constructions, like retractions, wedges, etc. These concepts per se are not limited to minimal submanifolds. Thus, we collect these concepts together in the next several subsections and give the formal definitions. The constructions all depend geometrically on the ambient metric and submanifolds involved.
	
	Unfortunately, our geometric set up will constantly involve three levels of manifolds. We need a submanifold $S^{d-c}$ that are formed by the transverse intersection of two $d$-dimensional submanifolds $X$ and $Y$, all of which are again contained in the ambient manifold $M$. Thus frequently we have to deal with several different notions of normal bundles and neighborhoods. To this end, we give a formal definition.
	
	For a sequence of inclusion of submanifolds $$X_1\s X_2\s X_3\cdots\s X_n= M,$$ we use $\T X_j$ to denote tangent bundle of $X_j.$ Then we have a sequence of natural inclusion maps
	\begin{align*}
		\T X_1\s \T X_2\s\cdots \s \T M \end{align*}  When $X_i\s X_j$, we use \begin{align*}
		\T^\perp({X_i, X_j})
	\end{align*}to denote the normal bundle of $X_i$ as a Riemannian submanifold of $X_j.$ 
	
	We use $\exp^\perp_{X_i,X_j}:\T^\perp({X_i, X_j})\to X_j$ to denote the exponential map of the normal bundle of $X_i$ considered as a Riemannian submanifold of $X_j$. 
	\begin{defn}
		Define the radius $r$ ball-bundle of $X_i$ inside $X_j$ as
		\begin{align*}
			B_r^\perp(X_i, X_j)=\{(p,v)\,|\,p\in X_i,v\in \T_p X_j,v\perp \T_p X_i,\no{v}\le r\}\s \T^\perp({X_i, X_j}).
		\end{align*}
		If $\exp^\perp_{X_i, X_j}$ is \textbf{ injective wherever it is well defined} in $B_r^\perp(X_i, X_j)$, we define the radius $r$ neighborhood of $X_i$ inside $X_j$ to be
		\begin{align*}
			B_r(X_i, X_j)=\exp^\perp_{X_i, X_j}(B_r^\perp(X_i, X_j)).
		\end{align*}
		
		On $B_r(X_i,X_j)$ we denote the nearest distance projection of $B_r(X_i,X_j)$ onto $X_i$ by $
			\pi_{X_i,X_j}.$
	\end{defn}
	Note that by definition, $\pi_{X_i,X_j}$ is the composition of the inverse of $\exp^\perp_{X_i,X_j}$ with the projection from the normal bundle $\T^\perp({X_i,X_j})$ to its base $X_i.$
		\subsection{Allard's regularity theorem}
	Now that we have set up the basic concepts, it is time to carry out Steps (\ref{im1}) and (\ref{im2}). First recall Lemma \ref{nsmin}, Definition \ref{defnur} and Definition \ref{defngm}. 
	
	The following lemma does not need any assumption on $C$ beyond Assumption \ref{assumpc}. 
	\begin{lem}\label{localt}
		If $2[\Si]\not=0,$ for any $\ee>0$, there is an open subset $\Om_h$ containing $h$ in the space of Riemannian metrics, such that any area-minimizing integral current $T$ in $[\Si]$ with respect to metric $g\in\Om_h$ satisfies 
		\begin{align}\label{cst}
			\supp T\s U_\ee(N\#\ssc),
		\end{align} and $T$ restricted to $$\bigg(U(N\cap N)\cup U_\ee(\sing\ssc)\bigg)\cp$$ is a smooth submanifold with boundary and $\ee$-close (Definition \ref{defngm}) to $N\#\ssc$ restricted to  $\bigg(U(N\cap N)\cup U_\ee(\sing\ssc)\bigg)\cp$. 
		
		If $2[\Si]=0,$ replace every $N\#\ssc$ above with $\pm N\#\ssc,$ i.e., $N\#\ssc$ or $N\#\ssc$ with reverse orientation.
	\end{lem}
	\begin{proof}
	Let us first deal with the case of $2[\Si]\not=0.$	We claim that 
		\begin{claim}\label{ccst}
			For any $\e'>0,$ there is an open set $\Om_h$ containing $h$ in the space of Riemannian metrics such that any area-minimizing integral current $T$ in $[\Si]$ with respect to metric $g\in\Om_h$ is of distance at most $\e'$ away to $N\#\ssc$ in the following three senses with respect to  our base metric $h$:\begin{enumerate}
				\item $T$ and $N\#\ssc$ as integral currents are at most $\e'$ away in the Whitney flat distance in the space of currents,
				\item The induced mass measures of $T$ and $N\#\ssc$ are at most $\e'$ away in the total variation distance in the space of measures,
				\item The supports $\supp T$ and $\supp N\#\ssc$ are at most $\e'$ away in Hausdorff distance in the space of closed subsets.
			\end{enumerate}
		\end{claim}

		Note that the above claim already implies (\ref{cst}).
		
		Suppose the above claim is not true. Then we can find a sequence of Riemannian metrics $h_j\to h$ in smooth topology and a sequence of area-minimizing integral currents $T_j\in[\Si]$ with respect to metric $h_j$ that violates at least one of the three $\e'$ distance conditions.
		
		Let us first prove that $\{T_j\}$ as a sequence in integral currents must converge to $N\#\ssc$. Let $\{T_{j_k}\}$ be any subsequence of $\{T_j\}$. By the compactness theorem for integral currents \cite[Theorem 4.2.17]{HF}, we can find a subsequence $\{T_{j_{k_m}}\}$ of $\{T_{j_k}\}$ such that $\lim_{m\to\infty} T_{j_{k_m}}$ exists. Since each $T_j$ is area-minimizing and $N\#\ssc$ is the unique area-minimizing representative of $[\Si]$ in $h$ (Lemma \ref{nsmin}, $2[\Si]\not=0$) we deduce that $\lim_{m\to\infty}T_{j_{k_m}}= N\#\ssc$ as integral currents by \cite[Theorem 34.5]{LS}. Now recall the classical fact that if a every subsequence of a sequence has a subsequence that coverges to the the same limit, then the original sequence must converge to the same limit. We deduce that $T_j$ converges to $N\#\ssc$ as integral currents.

		Since $N\#\ssc$ and each $T_j$ are all area-minimizing, \cite[Theorem 34.5]{LS} also implies that $T_j\to N\#\ssc$ as their induced mass measures. Finally $\supp T_j$ converge to $\supp N\#\ssc$ in Hausdorff distance by upper semi-continuity of densities of stationary varifolds \cite[Corollary 17.8]{LS}. Thus we have arrived at a contradiction.
		
		To prove that $N\#\ssc$ and $T$ are $\ee$-close when restricted to 	$\bigg(U(N\cap N)\cup U_\ee(\sing\ssc)\bigg)\cp$ as submanifolds is a classical application of Allard's regularity theorem \cite[Section 8, Regularity Theorem]{WA}. We will only give a sketch. 
		
		Roughly speaking the core of Allard's regularity theorem \cite[Section 8, Theorem 8.19]{WA} is that volume close to balls implies close to balls as submanifolds. Roughly speaking, if the area of a $d$-dimensional stationary varifold $V$ restricted to a $(d+c)$-dimensional geodesic ball centered at $p$ on $M$ is at most $\kappa$ off from the area of a standard $d$-dimensional ball in $\R^{d+c}$, then the varifold $V$ is $\iota$ close as submanifolds to a $d$-dimensional standard ball in a neighborhood of $p$ in $M$.

		Precisely speaking, Allard's regularity theorem implies the following. For any $\iota>0,$ there is $R,\kappa>0,$ such that if we have a $d$-dimensional stationary varifold $V$ in $M$ with metric $h$ and the area of $V$ inside a radius $R$ geodesic ball centered at $p$ on $M$ is at most $\kappa$ off from the volume of a $d$-dimensional radius $R$ standard  ball  in Euclidean space, then in $B_{\frac{1}{2}R}(p)$ on $M$, $V$ is $\iota$ close as a submanifold (Definition \ref{defngm}) to the exponentiated image of the tangent plane $P$ to $V$ at $p.$ 
		(Allard's original theorem is stated with $C^{1,\ai}$ bounds for varifolds with bounded mean curvature in Euclidean space. However, routine elliptic PDE arguments like \cite[Theorem 5.2.15 (7)]{HF}, upgrades  $C^{1,\ai}$ estimates into $C^{k,\ai_k}$ estimate with all $k>0$.)
		
		Now cover $\bigg(U(N\cap N)\cup U_\ee(\sing\ssc)\bigg)\cp$ by finitely many radius $\e''$ balls such that none of the balls in the covering intersect $\sing\left( N\#\ssc\right)$. Since $N\#\ssc$ restricted to $\bigg(U(N\cap N)\cup U_\ee(\sing\ssc)\bigg)\cp$ is a smooth submanifold, for any $\e'>0,$ there is $\e''>0$ such that $N\#\ssc$ restricted to any of the radius $\e''$ balls in the covering has area at most $\e'$ off from that of a radius $\e''$ standard $d$-dimensional ball in $\R^{d+c}$.

		Now since the area-minimizing current $T$ is at most $\e'$ away from $N\#\ssc$ as induced mass measures by Claim \ref{ccst}, we obtain that $T$ restricted to each radius $\e''$ ball in the covering has area  at most $2\e'$  off from that of a $d$-dimensional radius $\e''$ standard ball in $\R^{d+c}$. Now we can apply Allard's regularity theorem to deduce that $T$ restricted to $\bigg(U(N\cap N)\cup U_\ee(\sing\ssc)\bigg)\cp$ is $\ee$-close to $N\#\ssc$ restricted to $\bigg(U(N\cap N)\cup U_\ee(\sing\ssc)\bigg)\cp$ as submanifolds, provided $\e'$ is small. By making $\e'$ arbitrarily small, we can take $\ee$ arbitrarily small. We are done.
		
		In the case of $2[\Si]=0,$ the above argument still goes through with every $N\#\ssc$ replaced by $\pm N\#\ssc,$ since by Lemma \ref{nsmin} the only area-minimizing integral currents in $[\Si]$ are $\pm T$ for $2[\Si]=0.$
	\end{proof}
	Now recall Definition \ref{defnur}.
	\begin{lem}\label{localti}(Candidates $X+Y$ for replacing $T|_{U_{{\eta}}(\sing \sdc)}$)
		If $2[\Si]\not=0,$ For any $\ees>0$, and any $0<\eta\le\rr$ (Definition \ref{defnur}) there is an open subset $\Om_h$ containing $h$ in the space of Riemannian metric, such that any area-minimizing integral current $T$ in $[\Si]$ with respect to metric $g\in\Om_h$ satisfies that
		\begin{align}\label{dti}
			\pd(T|_{U_{{\eta}}(\sing \sdc)})=\pd(X+Y),
		\end{align} where $X,Y$ are smooth compact \textbf{minimal} submanifolds with boundaries $\pd X,\pd Y$, respectively, and $X,Y$ are supported in $U_\eta(\sing\sdc)$ and are $\ees$-close as submanifolds to the union of 
		\begin{align}\label{pm1}
			\Ga\m({B_\eta^c\times\{0\}^c}\times S^{d-c})
		\end{align} and  \begin{align}\label{pm2}
		\Ga\m({\{0\}^c\times B_\eta^c}\times S^{d-c}),
		\end{align} respectively. And $X,Y$ intersect transversely in their interior along a $(d-c)$-dimensional submanifold 
		$	L,
		$ that is $\ees$-close to 
		\begin{align}\label{pm3}
			\Ga\m(\{0\}^{2c}\times S^{d-c})
		\end{align} as a submanifold.
		
		If $2[\Si]=0,$ then the above still holds by adding a $\pm$ sign to (\ref{pm1}) (\ref{pm2}) and (\ref{pm3}), i.e., either positively oriented or negatively oriented.
	\end{lem}
	Since  $\Ga$ provides an isometric diffeomorphism from $U_\rr(\sing\sdc)$ to $B_\rr^{2c}\times S^{d-c}$, we will omit the symbol $\Ga$ and $\Ga\m$ frequently from now on.
	\begin{proof}		The proof of this fact is similar to the proof of Lemma \ref{localt}. For $2[\Si]\not=0,$ by Lemma \ref{localt}, for any $0<\ee<\eta$, there is an open subset $\Om_h$ containing $h$ in the space of Riemannian metrics, such that any area-minimizing integral current $T$ in $[\Si]$ with respect to metric $g\in\Om_h$ satisfies 
		\begin{align*}
			\pd(T|_{U_\eta(\sing\sdc)})=A+B,
		\end{align*} where each summand  is $\ee$-close to the each summand of
		\begin{align*}
			\pd(N\#\ssc|_{U_\eta(\sing\sdc)})=\Ga\m\pf(S^{c-1}_\eta\times\{0\}^c\times S^{d-c})+\Ga\m\pf(\{0\}^c\times S^{c-1}_\eta\times S^{d-c}),
		\end{align*}as submanifolds, respectively. Here $A,B$ are disjoint submanifolds diffeomorphic to $S^{c-1}\times S^{d-c}.$ The symbol $\Ga\m\pf$ means the pushforward of integral currents under $\Ga\m,$ and the reader can just understand this as embeddings under the map $\Ga\m$.
		
		Now solve the Plateau problem separately for boundary $A$ and $B$ to find solutions $X,Y$ with $\pd X= A$ and $\pd Y= B$, respectively, i.e., finding integral currents $X,Y$ with boundary $A,B$ respectively of minimal possible area.
		
		Since $\Ga\m\pf({B_\eta^c\times\{0\}^c}\times S^{d-c})$ is the unique area-minimizing integral current with boundary $\Ga\m\pf(S^{c-1}_\eta\times\{0\}^c\times S^{d-c})$ (Lemma \ref{nsmin}), we can argue like in the proof of Lemma \ref{localt} to deduce the following:

			 $X$ is $\ees$-close to $\Ga\m\pf({B_\eta^c\times\{0\}^c}\times S^{d-c})$ as a submanifold, provided $\pd X$ and $\Ga\m\pf(S^{c-1}_\eta\times\{0\}^c\times S^{d-c})$ are $\ee$-close as submanifolds with $\ee$ small enough.

 This time we use both Allard's regularity theorem \cite{WA} and Allard's boundary regularity theorem \cite[Section 4 Regularity Theorem]{WAb}. Similarly we can conclude that $Y$ is $\ees$ close to $\Ga\m\pf(\{0\}^c\times S^{c-1}_\eta\times S^{d-c})$ as submanifolds. This implies that $X$ and $Y$ must be supported in $U_\eta(\sing\sdc),$ provided $\ee$ is small enough.
		
		Since $X$ and $Y$ are $\ees$-close to  $\Ga\m({B_\eta^c\times\{0\}^c}\times S^{d-c})$ and  $\Ga\m({\{0\}^c\times B_\eta^c}\times S^{d-c})$ as submanifolds, $X$ and $Y$ must intersect transversely in  their interior along a closed $(d-c)$-dimensional submanifold $L$ that is $\ees$ close to the transverse intersection of $\Ga\m({B_\eta^c\times\{0\}^c}\times S^{d-c})$ and  $\Ga\m({\{0\}^c\times B_\eta^c}\times S^{d-c})$ as submanifolds. We are done with the case of $2[\Si]\not=0.$

		In the case of $2[\Si]=0,$ the above argument still goes through with every $N\#\ssc$ replaced by $\pm N\#\ssc,$ since by Lemma \ref{nsmin} the unique area-minimizing integral currents are $\pm N\#\ssc$ for $2[\Si]=0.$ This induces a $\pm$ sign into (\ref{pm1}) (\ref{pm2}) and (\ref{pm3}).
	\end{proof}
	From now on in this section, based on Lemma \ref{localt} and Lemma \ref{localti}, we assume that \begin{assump}\label{assumplxm}
		$X$ and $Y$ are smooth $d$-dimensional compact submanifolds with boundary in $M$ with respect to an ambient metric $h\in\Om_h$ (Lemma \ref{localti}), and $L$ is a $(d-c)$-dimensional closed submanifold of $M$, such that $L$ is formed by the transverse intersection of $X$ and $Y$ in their interior and we have
		\begin{align*}
\injrad_X>\rr,&\textnormal{ }\injrad_Y>\rr,\\\focrad_{X,M}>\rr,&\textnormal{ }\focrad_{Y,M}>\rr,\\
\injrad L>\frac{\pi}{2}>\rr,\\
X\s B_\rr(L,X),&\textnormal{ } Y\s B_\rr(L,Y).
		\end{align*} 
	\end{assump}
	Here $\injrad_X$ is the injectivity radius of $X$. It  is defined to be the maximum of $\rh$ over all points $p\in X$ such that $\exp^\perp_{\{p\},X}$ is injective in $B_\rh^\perp(p,X)$ wherever it is defined. Also $\focrad_{X,M}$ is the focal radius of $X$ inside $M$, i.e., the maximum of $\rh$ over all points $p\in X$ such that the exponential map $\exp^\perp_{X,M}$ restricted to $B_\rh^\perp(X,M)$ at $p$ is well-defined and injective.
	
	In other words the above assumption roughly says that every point $q\in X$ is of distance at most $\rr$ to $L$ and is in the injective image of the exponential map of the normal bundle of $L$ inside $X.$ 
	\subsection{Double Fermi coordinates}
	In this section with Assumption \ref{assumplxm}, we will define the double Fermi coordinate chart at a point $p\in L$ adapted to $(L,X,M)$ with radius $(\rr,\rr,\rr)$ as follows. 
	\subsubsection{Normal coordinates of radius $\rr$ on $L$}\label{secl}
	For any point $p$ in the Riemannian manifold $L$, take a normal coordinate chart \begin{align*}
		l=(l_1,\dots,l_{d-c}).
	\end{align*}
	centered at $p$ with radius $r$. In other words we take an orthonormal frame 
	\begin{align}\label{ofl}
\{		\ell_1,\cd,\ell_{d-c}\}
	\end{align} in $\T_pL$, and set the coordinate point $l=(l_1,\dots,l_{d-c})$ to be the injective image of $\sum_j l_j\ell_j$ under the exponential map of $L$ at $p,$ with $\sum_j l_j^2\le\rr^2$ to utilize Assumption \ref{assumplxm}. Our coordinate chart is defined on the ball
	\begin{align*}
		B_\rr(p,L).
	\end{align*}
	
	If $d-c=0,$ then skip the above step. 
	
	Now define a radial vector field  $\pd_l$ in the $l$ directions by
	\begin{align*}
		\pd_l=\frac{l_j}{|l|}\pd_{l_j}.
	\end{align*}
	\begin{fact}\label{fctl}
		In the normal coordinate chart $(l_1,\dots,l_{d-c})$ we have the following:
		\begin{align*}
			\pd_l\perp_{|\cdot|}v&\iff\pd_l\perp v,\\
			\na_{\pd_l}\pd_l&=0,\no{\pd_l}=|\pd_l|=1,\\	
			\sqrt{\sum_j l_j^2}&=\fd_L\Big((l_1,\dots,l_{d-c}),\{0\}^{d-c}\Big),\\
			\na_{\pd_{l_i}}\pd_{l_j}(p)&=0,\\
			\{\pd_{l_1},\dots,\pd_{l_{d-c}}\}&\textnormal{ is an orthonormal frame at }(0,\dots,0).
		\end{align*}
	\end{fact}
	The symbol $\na$ means the Levi-Civita connection with respect to the ambient metric.
	Here $\pd_l\pp v$ at $l$ means that $\pd_l$ and $v$ are perpendicular as tangent vectors with respect to the quadratic norm $|\cdot|.$ In other words, if we write $v=v_j\pd_{l_j}$ with $v\in \T_{(l_1,\cd,l_{d-c})}L,$ then we have $\sum_j v_jl_j=0.$ On the other hand $\pd_l\perp v$ means $\pd_l$ is perpendicular to $v$ with respect to the Riemannian metric.
	
	All of the above facts are classical results about normal coordinates. For instance, the first line is just the Gauss lemma of exponential maps.
	\subsubsection{First Fermi coordinates of radius $(\rr,\rr)$}\label{secx} 
	Now take an orthonormal frame 
	\begin{align}\label{ofx}
		\{e_1,\dots,e_c\}
	\end{align} of the normal bundle $\T^\perp(L,X)$ at $p$. Parallel translating the frame $\{e_1,\dots,e_c\}$ along geodesics of $L$ passing through $p,$ we obtain a smooth orthonormal frame $\{e_1,\dots,e_c\}$ of $\T^\perp(L,X)$ restricted to $B_\rr(p,L),$ i.e., where the chart $(l_1,\dots,l_{d-c})$ is defined. By construction we have
	\begin{align*}
		\na_ve_j(p)=0,
	\end{align*}for any $v$ tangent to $L.$
	
	Now define a Fermi coordinate chart of radius $\rr$ \cite[Definition of Fermi coordinates on page 17]{AG}
	\begin{align*}
		(l,x)=(l_1,\dots,l_{d-c},x_1,\dots,x_c)
	\end{align*} adapted to $L$ in chart $(l_1,\dots,l_{d-c})$ with orthonormal frames $\{e_1,\dots,e_c\}$ and $\sum_j x_j^2<\boldsymbol{\rh}^2_{\operatorname{product}}.$ In other words the coordinate point $(l_1,\dots,l_{d-c},x_1,\dots,x_{c})$ is defined by the injective image of the vector $\sum_jx_je_j$  under the normal bundle exponential map $\exp^\perp_{L,X}$ in $\T^\perp(L,X)$ at point $(l_1,\dots,l_{d-c})$ in normal coordinate chart on $L$ centered at $p.$ 
	
	If $d-c=0,$ just take a normal coordinate chart $(x_1,\dots,x_{c})$ in this step. 
	
	By Assumption \ref{assumplxm}, every point $q$ in $X$ with $\pi_{L,X}(q)$ lying in chart $(l_1,\dots,l_{d-c})$ is contained in our coordinate chart $(l_1,\dots,l_{d-c},x_{1},\dots,x_c).$ 
	
	Define a radial vector field $\pd_x$ along the $x$-directions by
	\begin{align*}
		\pd_x=\frac{x_j}{|x|}\pd_{x_j}.
	\end{align*}
	\begin{fact}\label{fctx}In the Fermi coordinate chart $(l,x)$ we have the following properties:
		\begin{align*}
			\pd_x\pp v&\iff \pd_x\perp v,\\
			\na_{\pd_x}\pd_x&=0,\no{\pd_x}=|\pd_x|=1,\\
			\na_{\pd_{x_j}}\pd_{x_i}(p)&=\na_{\pd_{x_j}}\pd_{l_i}(p)=\na_{\pd_{l_i}}\pd_{x_j}(p)=0,\\	\pi_{L,X}(l,x)&=(l,0),\\
			\sqrt{\sum_jx_j^2}&=\fd_X\big((l,x),(l,0)\big),\\
		\{\pd_{x_1},\dots,\pd_{x_{d-c}}\}&\textnormal{ is an orthonormal frame on }\T^\perp(L,X).
		\end{align*}
	\end{fact}
	The proof of the above fact is a straightforward generalization of Fact \ref{fctl} about normal coordinates, and can be found in \cite[Section 2]{AG}.
	
	Fact \ref{fctl} holds in chart $(l,x)$ when restricted to $\{|x|=0\}.$
	\subsubsection{Double Fermi coordinate charts of radius $(\rr,\rr,\rr)$}\label{secy}
	Now let $\{f_1,\dots,f_c\}$ be an orthonormal frame of $\T^\perp(X,M)$ at point $p.$ Again, we can parallel translate $\{f_1,\dots,f_c\}$ along geodesics of $X$ passing through $p$ to obtain a smooth frame of $\T^\perp(X,M)$ on chart $(l_1,\dots,l_{d-c},x_1,\dots,x_c)$, such that
	\begin{align*}
		\na_ef_j(p)=0,
	\end{align*}whenever $e$ is tangent to $X$. Then define a Fermi coordinate chart with radius $\rr$: 
	\begin{align*}
		(l,x,z)=(l_1,\dots,l_{d-c},x_1,\dots,x_c,z_1,\dots,z_c),
	\end{align*}
	adapted to $X$ with chart $(l_1,\dots,l_{d-c},x_1,\dots,x_c)$ and frames $\{f_1,\dots,f_c\}$ of $\T^\perp(X,M)$ with $\sum_j z_j^2\le \rr^2.$ Any point $q$ in $B_\rr(X,M)$ with $\pi_{X,M}(q)$ lying in chart $(l_1,\cd,l_{d-c},x_1,\cd,x_c)$ is contained in our chart $(l,x,z)$
	\begin{defn}\label{defndfermi}
		We call the coordinate system $(l_1,\dots,l_{d-c},x_1,\dots,x_c,z_1,\dots,z_c)$ the double Fermi coordinate chart of radius $(\rr,\rr,\rr)$ adapted to $(L,X,M).$
	\end{defn}
	
	Now define a radial vector field $\pd_z$ along the $z$-directions by
	\begin{align*}
		\pd_z=\frac{z_i}{|z|}\pd_{z_i}.
	\end{align*}
	\begin{fact}\label{fcty}
		In the double Fermi coordinate $(l,x,z)$ we have the following properties
		\begin{align*}
			\pd_z\pp v&\iff\pd_z\perp v,\\
			\na_{\pd_z}\pd_z&=0,\no{\pd_z}=|\pd_z|=1,\\
			\na_{\pd_{z_i}}\pd_{z_j}(p)=\na_{\pd_{l_i}}\pd_{z_j}(p)&=\na_{\pd_{z_j}}\pd_{l_i}(p)=\na_{\pd_{x_i}}\pd_{z_j}(p)=\na_{\pd_{z_j}}\pd_{x_i}(p)=0,\\	\pi_{X,M}(l,x,z)&=(l,x,0),\\
			\sqrt{\sum_j z_j^2}&=\fd_M\big((l,x,z),(l,x,0)\big),\\
		\{\pd_{z_1},\dots,\pd_{z_{d-c}}\}&\textnormal{ is an orthonormal frame on }\T^\perp(X,M).
		\end{align*}
	\end{fact}
	Again the proof of the above fact is straightforward and the reader can refer to \cite[Chapter 2]{AG}.
	
	Fact \ref{fctl} holds in chart $(l,x,z)$ when we restrict to $\{|x|=|z|=0\}.$ Fact \ref{fctx} holds in chart $(l,x,z)$ when we restrict to $\{|z|=0\}.$
\subsection{Geometric dependence of the double Fermi coordinate chart}
At first sight, it may seem that the double Fermi coordinate chart centered at $p$ defined in Definition \ref{defndfermi} depends on a lot of arbitrary choices. However, it turns out the double Fermi coordinate chart is really intrinsically defined up to orthogonal transformations:
\begin{fact}
		The double Fermi coordinate chart defined in Definition \ref{defndfermi} is unique up to an $$\operatorname{O}(d-c)\times\operatorname{O}(c)\times\operatorname{O}(c)$$ matrix multiplication in coordinate.
\end{fact}
\begin{proof}
	By construction, our construction of double Fermi coordinate chart only depends on the choice of orthonormal frames $\{\ell_1,\cd,\ell_{d-c}\}$ of $\T_pL$, $\{e_1,\dots,e_c\}$ of $\T_p^\perp(L,X)$ and $\{f_1,\dots,f_c\}$ of $\T_p^\perp (X,M)$. All choices of orthonormal frames is equivalent up to an $\operatorname{O}(d-c)\times\operatorname{O}(c)\times\operatorname{O}(c)$ action and by definition of the double Fermi coordinate chart. This $\operatorname{O}(d-c)\times\operatorname{O}(c)\times\operatorname{O}(c)$ action extends to all points of our coordinates.
\end{proof}
An immediate corollary of the above fact is that
\begin{fact}\label{dfcgd}
The coordinate lengths $|l|,|x|,|z|,$ the vector fields $\pd_l,\pd_x,\pd_z$, and all ambient metric inner products among coordinate vector fields and the length of the following multivectors and the length of their wedge products
\begin{align*}
	&\pd_{l_1}\w\cd\w \pd_{l_{d-c}},\\
	&\pd_{x_1}\w\cd\w \pd_{x_c},\\
	&\pd_{z_1}\w\cd\w \pd_{z_c},
\end{align*} all depend geometrically (Definition \ref{defngdp}) on $L,X,M$ and the ambient metric.
\end{fact}
Without loss of generality, we can assume that
\begin{align*}
	U_{\frac{\rr}{2}}(\sing\sdc)
\end{align*}is contained in the open cover of double Fermi coordinate charts of radius $(\rr,\rr,\rr)$ centered at $p$ with $p$ ranging over all points in $L.$ Since $	U_{\frac{\rr}{2}}(\sing\sdc)$ has compact closure, we can take a finite subcover.
\begin{assump}\label{assumpdfc}
	There exists a finite subset $\mathbf{F}\s L$ such that the union of double Fermi coordinate charts of radius $(\rr,\rr,\rr)$ centered at points $p\in\mathbf{F}$ covers $	U_{\frac{\rr}{2}}(\sing\sdc)$.
\end{assump}
Form now on when we mention a double Fermi coordinate chart, we always mean a double Fermi coordinate chart of radius $(\rr,\rr,\rr)$ centered at a point $p\in\bff.$
	\subsection{Angled wedges and retractions}
\begin{figure}[h]
		\centering
		\def\svgwidth{0.85\paperwidth}
		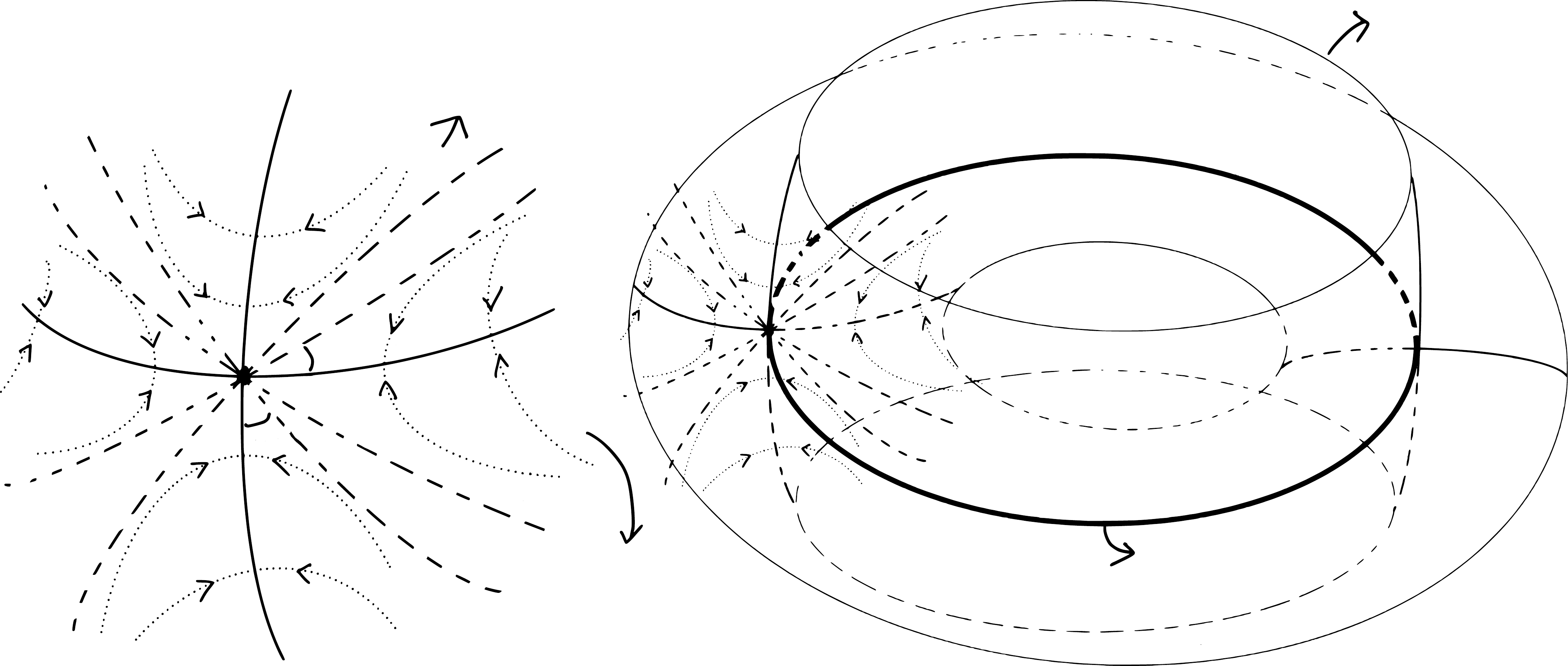
		\caption{A schematic illustration of $\Pi_{X,L},\Pi_{Y,L}$}
	\label{figpxy}
	\end{figure}

	\begin{defn}\label{defnw}
		With Assumption \ref{assumplxm}, we define the radius $\rh$, angle $\ta$, wedge neighborhood of $X$ with vertex $L$ to be
		\begin{align*}
			W_{(\ta,\rh)}(X,L)=\left\{q\in B_\rh(X,M)\bigg|\frac{\fd_M(q,\pi_{X,M}(q))}{\fd_X(\pi_{X,M}(q),\pi_{L,X}\circ\pi_{X,M}(q))}\le\tan\ta\right\},
		\end{align*} provided $0<\rh<\rr.$
	\end{defn}
	By Facts \ref{fctl}, \ref{fctx}, and \ref{fcty}, the geometric meaning of $W_{(\ta,r)}(X,L)$ is that when restricted to the double Fermi coordinate charts $(l,x,z)$ in Definition \ref{defndfermi}, we have
	\begin{align*}
		W_{(\ta,\rh)}(X,L)=\left\{(l,x,z)\bigg|\frac{|z|}{|x|}\le \tan\ta\right\}.
	\end{align*}
	The reader should consult Figure \ref{figpxy} and (\ref{woe}) for intuition.

	With Assumption \ref{assumplxm},	define Lawlor's retraction $\Pi_{X,L}$ as follows. 
	\begin{defn}\label{defnpl}
		For any point $q\in B_\rr(X,M)$, there is a unique point on the length minimizing geodesic of $X$ from $\pi_{X,M}(q)$ to $\pi_{L,X}\circ\pi_{X,M}(q)$, such that the point is precisely of $\fd_X$ distance
		\begin{align*}
			\cg_c\bigg(\frac{\fd_M(q,\pi_{X,M}(q))}{\fd_X\big(\pi_{X,M}(q),\pi_{L,X}\circ\pi_{X,M}(q)\big)}\bigg)^{\frac{1}{c}}\fd_X\big(\pi_{X,M}(q),\pi_{L,X}\circ\pi_{X,M}(q)\big)
		\end{align*}away from $\pi_{L,X}\circ\pi_{X,M}(q)$. Define
		\begin{align*}
			\Pi_{X,L}(q)
		\end{align*}to be this unique point.
	\end{defn}
	Here $\cg_c$ is defined in Definition \ref{defng}.
	By Facts \ref{fctl}, \ref{fctx} and \ref{fcty}, the geometric significance of $\pi_{X,L}$ is that in the double Fermi coordinates, we have
	\begin{align}\label{pxlcdnt}
		\Pi_{X,L}(l,x,z)=\bigg(l,\cg_c^{\frac{1}{c}}\bigg(\frac{|z|}{|x|}\bigg)x,0\bigg).
	\end{align}
	The reader should consult Figure \ref{figpxy} and (\ref{poe}) for intuition.
		\subsubsection{Calculating $\ed\Pi_{X,L}$}
	With the coordinate expression (\ref{pxlcdnt}) of $\Pi_{X,L}$ in double Fermi coordinate charts, we can calculate the differential $\ed\Pi_{X,L}.$ 
	
	At any point $$q=(l,x,z)$$ in our double Fermi coordinate chart centered at $p,$  use $\nu^x$ to denote any vector in $\T_qM$ with only nonzero components among $x_1,\dots,x_c$ axes directions and $$\nu^x\pp \pd_x.$$ We can then extend $\nu^x$ to a smooth vector field on our double Fermi coordinate chart by keeping all coefficients constant. In other words, for $\nu^x=\sum_j a^j\pd_{x_j}$ at $q=(l,x,z)$, we define $\nu^x=\sum_j a^j\pd_{x_j}$ everywhere.
	
	 Similarly, use $\nn^z$ to denote a vector in $\T_qM$ with only nonzero components among $z_1,\dots,z_c$ axes directions and $$\nn^z\pp \pd_z.$$ Then extend $\nn^z$ constantly in our double Fermi coordinate chart. By abusing the notations, we can calculate with Definition \ref{defng} in mind,
	\begin{lem}\label{dker}At point $q=(l,x,z)$, if $\frac{|z|}{|x|}>t_c$ (defined \ref{defntc}), we have
		\begin{align*}
			\ed \Pi_{X,L}(v,w,u)=(v,0,0),
		\end{align*}and if $\frac{|z|}{|x|}<t_c$, we have
		\begin{align*}
			\ed\Pi_{X,L}(v,0,0)=&(v,0,0),\\
			\ed\Pi_{X,L}(0,\pd_x,0)=&\bigg(0,\cg_c\bigg(\frac{|z|}{|x|}\bigg)^{\frac{1}{c}-1}\cc_c\bigg(\frac{|z|}{|x|}\bigg)\pd_x,0\bigg),\\
			\ed\Pi_{X,L}(0,\nu^x,0)=&\bigg(0,\cg_c\bigg(\frac{|z|}{|x|}\bigg)^{\frac{1}{c}}\nu^x,0\bigg),\\
			\ed\Pi_{X,L}(0,0,\pd_z)=&\bigg(0,\cg_c\bigg(\frac{|z|}{|x|}\bigg)^{\frac{1}{c}-1}\cs_c\bigg(\frac{|z|}{|x|}\bigg)\pd_x,0\bigg),\\
			\ed\Pi_{X,L}(0,0,\nn^z)=&(0,0,0),
		\end{align*}where the right hand sides always lie in $\T_{\Pi_{X,L}(q)}M.$
		Moreover the kernel subspace $\ker \ed\Pi_{X,L}$ is spanned by
		\begin{align}\label{kpxy}
			-\cs_c\bigg(\frac{|z|}{|x|}\bigg)\pd_x+\cc_c\bigg(\frac{|z|}{|x|}\bigg)\pd_z,
		\end{align}
		and all $\nn^z$ with only nonzero components among $z$ axes directions such that $$\nn^z\pp \pd_z.$$
	\end{lem}
	\begin{proof}
		By definition, $\nu^x(|x|)=\nu^x(|z|)=\nn^z(|z|)=\nn^z(|x|)=0,$ so only $\ed\Pi_{X,L}(0,\pd_x,0)$ and $\ed\Pi_{X,L}(0,0,\pd_z)$ need our attention.
		
		On the other hand, we have $\pd_x|x|=1,\frac{x}{|x|}=\pd_x$ and $\pd_z|z|=1,\frac{z}{|z|}=\pd_z.$ Thus, we have
		\begin{align*}
			\ed\Pi_{X,L}(0,\pd_x,0)=&\cgxyc \pd_x+\frac{1}{c}\cgxyco\cg_c'\xy\bigg(-\frac{|z|}{|x|^2}\bigg)x\\
			=&	\cgxyco\bigg(\cgxy-\frac{1}{c}\frac{|z|}{|x|}\cg_c'\xy\bigg)\pd_x\\
			=&\cgxyco\ccxy \pd_x.
		\end{align*}
		We also have
		\begin{align*}
			\ed\Pi_{X,L}(0,0,\pd_z)=&\frac{1}{c}\cgxyco\cg_c'\xy\bigg(\frac{1}{|x|}\bigg)x\\
			=&	\cgxyco\bigg(\frac{1}{c}\cg_c'\xy\bigg)\pd_x\\
			=&\cgxyco\csxy \pd_x.
		\end{align*}
		Then (\ref{kpxy}) is apparently mapped to $0$ by $\ed\Pi_{X,L}.$ We are done.
	\end{proof}
	Note that every construction in this section works by replacing $X$ with $Y$ in Assumption \ref{assumplxm}.
	
	In other words we can define
	\begin{align*}
		W_{(\ta,\rh)}(Y,L)=\left\{q\in B_\rr(Y,M)\bigg|\frac{\fd_M(q,\pi_{Y,M}(q))}{\fd_X(\pi_{Y,M}(q),\pi_{L,Y}\circ\pi_{Y,M}(q))}\le\tan\ta\right\},
	\end{align*}by replacing $X$ with $Y$
	in Definition \ref{defnw} and define \begin{align}\label{defnply}
		\Pi_{Y,L},
	\end{align}
	by replacing $X$ with $Y$ in Definition \ref{defnpl}.\subsection{Intersection angles and wedge bundles}
	As mentioned when introducing Lawlor's ideas, in the end we have to verify that the wedges  $W_{(\ta,\rr)}(X,L)$ and $W_{(\ta,\rr)}(Y,L)$  in the last subsection intersect only at $L$. Testing how wedges intersect is a bit tricky. To this end, we need an angle condition and an auxiliary wedge bundle. 
	\begin{fact}\label{fctdec}
		For every point $p\in L,$ a tangent vector $v$ in $\T_p M$ can be uniquely written as a sum of orthogonal vectors:
		\begin{align*}
			v=v_L+v_{L,X}+v_{X,M},
		\end{align*}
		where $v_L\in \T_pL,v_{L,X}\in \T_p^\perp(L,X),v_{X,M}\in \T_p^\perp(X,M)$.
	\end{fact}
	\begin{defn}\label{defnta}
		Recall Assumption \ref{assumplxm} and Fact \ref{fctdec}, we will define the intersection angles $\Ta(X,Y)$ between $X$ and $Y$. Geometrically speaking, $\Ta(X,Y)$ is the minimum possible angle between vectors in $\T^\perp(L,X)$ and $\T^\perp(L,Y)$:
		\begin{align*}
			\Ta(X,Y)=\inf_p\inf_{v\in\T^\perp_p(L,Y)}\arctan\frac{\no{v_{X,M}}}{\no{v_{L,X}}}. 
		\end{align*}
	\end{defn}
	The tangent space to a submanifold captures the infinitesimal behavior of the submanifold. Along the same line, we need to capture the the infinitesimal structures of wedge neighborhoods in Definition \ref{defnw}, $W_{(\ta,\rh)}(X,L)$. To this end, we define a wedge bundle  $W_{(\ta,\rh)}^\T(X,L)$ associated with the wedge neighborhood $W_{(\ta,\rh)}(X,L)$. Geometrically speaking, it just consists of all tangent vectors in $\T^\perp(L,M)$ that lie in the wedge neighborhood $W_{(\ta,\rh)}(X,L)$:
	\begin{defn}\label{defnwb}
		Define the angle $\ta$ radius $\rh$ wedge bundle of $X$ with vertex $L$ to be
		\begin{align*}
			W_{(\ta,\rh)}^\T(X,L)=\bigg\{(p,v)\in \T M\,|\,p\in L,v\in \T_p^\perp(L,M),\no{v_{X,M}}\le\rh,\frac{\no{v_{X,M}}}{\no{v_{L,X}}}\le\tan\ta\bigg\}.
		\end{align*}
	\end{defn}
	Note that for the definition of wedge bundles $W^\T_{(\ta,\rh)}(X,L),$ we no longer need $\rh\le\rr,$ and even $\rh=\infty$ can be allowed.
	Similarly one can define $$W^\T_{(\ta,\rh)}(Y,L)$$ with every $X$ replaced by $Y$ in Fact \ref{fctdec} and Definition \ref{defnwb} . A very nice thing about the wedge bundles is that it is always easy to determine the intersection set of two wedge bundles $W^\T_{(\ta,\rh)}(X,L)$ and $W^\T_{(\ta,\rh)}(Y,L)$, since one just has to do the algebraic problem in a vector bundle. We have the following fact.
	\begin{lem}\label{lemwi}
		With Assumption \ref{assumplxm} in mind, we have
		\begin{align*}
			\Ta(X,Y)=\Ta(Y,X).
		\end{align*}Furthermore,
		\begin{align*}
			W^\T_{(\ta,\rh)}(X,L)\cap		W^\T_{(\ta,\rh)}(Y,L)=L,
		\end{align*}provided
		\begin{align*}
			\ta<\min\left\{\frac{\Ta(X,Y)}{2},\frac{\pi}{4}\right\}.
		\end{align*}
	\end{lem}
	\begin{proof}
		Let us first prove that $\Ta(X,Y)=\Ta(Y,X).$ 
		
		For any point $p\in L,$ let $S_1^{2c-1}(p)$ denote the unit sphere in $\T_p^\perp(L,M).$ Then $\T_p^\perp(L,X)$ and $\T_p^\perp(L,Y)$ intersect $S_1^{2c-1}$ at two great circle (i.e. totally geodesic) spheres $S_X^{c-1}$ and $S_Y^{c-1}$, respectively, of dimension $c-1$. 
		
		Now let $v$ denote a unit vector in $S_Y^{c-1}.$ Since $S_X^{c-1}$ is totally geodesic, elementary spherical geometry shows that
		\begin{align}\label{atd}
			\arctan\frac{\no{v_{X,M}}}{\no{v_{L,X}}}=\fd_{S_1^{2c-1}(p)}(v,S_{X}^{c-1}).
		\end{align}
		Thus, we have
		\begin{align*}
			&\inf_{v\in\T^\perp_p(L,Y)}\arctan\frac{\no{v_{X,M}}}{\no{v_{L,X}}}\\=&\inf_{v\in S_Y^{c-1}}\inf_{w\in S_X^{c-1}}\fd_{S_1^{2c-1}(p)}(v,w)\\=&\inf_{w\in S_X^{c-1}}\inf_{v\in S_Y^{c-1}}\fd_{S_1^{2c-1}(p)}(v,w)\\
			=&\inf_{w\in\T^\perp_p(L,X)}\arctan\frac{\no{v_{Y,M}}}{\no{v_{L,Y}}}.
		\end{align*}
		Since $p$ is arbitrary, this implies that $\Ta(X,Y)=\Ta(Y,X).$
		
		Using interpretation (\ref{atd}), it is easy to show that when restricted to a point $p$, the two infinite radius wedge bundles
		\begin{align*}
			W^\T_{(\ta,\infty)}(X,L),
			W^\T_{(\ta,\infty)}(Y,L)
		\end{align*}are the cones over the closed domains
		\begin{align*}
			B_{\ta}(S_X^{c-1},S^{2c-1}_1(p))=	\{v\in S^{2c-1}_1(p)|\fd_{S^{2c-1}_1(p)}(v,S_X^{c-1})\le \ta\},\\
			B_{\ta}(S_Y^{c-1},S^{2c-1}_1(p))=	\{v\in S^{2c-1}_1(p)|\fd_{S^{2c-1}_1(p)}(v,S_Y^{c-1})\le \ta\},
		\end{align*}respectively. Thus, if $z\in 
		W^\T_{(\ta,\infty)}(X,L)\cap
		W^\T_{(\ta,\infty)}(Y,L)$ at $p,$ then
		\begin{align*}
			\fd_{S^{2c-1}_1(p)}(\no{z}\m z,S_X^{c-1})\le\ta,\\\fd_{S^{2c-1}_1(p)}(\no{z}\m z,S_Y^{c-1})\le\ta.
		\end{align*}
		By the triangle inequality for distance, this implies that there are points $a\in S_X^{c-1}$, and $b\in S_Y^{c-1}$ such that
		\begin{align*}
			\fd_{S^{2c-1}_1(p)}(a,b)\le \fd_{S^{2c-1}_1(p)}(a,\no{z}\m z)+\fd_{S^{2c-1}_1(p)}(\no{z}\m z,b)\le 2\ta.
		\end{align*}
		This is a contradiction to the fact that  $$\inf_{w\in S_X^{c-1}}\inf_{v\in S_Y^{c-1}}\fd_{S_1^{2c-1}(p)}(v,w)\ge \Ta(X,Y)>2\ta.$$ Thus, restricted to the point $p,$ we must have $W^\T_{(\ta,\infty)}(X,L)\cap
		W^\T_{(\ta,\infty)}(Y,L)=p.$ Since $p$ is arbitrary in $L$, we deduce that
		\begin{align*}
			L\s 	W^\T_{(\ta,\rh)}(X,L)\cap
			W^\T_{(\ta,\rh)}(Y,L)\s 	W^\T_{(\ta,\infty)}(X,L)\cap
			W^\T_{(\ta,\infty)}(Y,L)=L.
		\end{align*}We are done.
	\end{proof}
	With Lemma \ref{lemwi} at hand, a natural idea is to use wedge bundles to help us determine how wedge neighborhoods intersect.
	
	However, if a submanifold is not totally geodesic, then exponentiating its tangent bundle in the ambient manifold will not give us the submanifold. 
	Similarly, in general
	\begin{align*}
		\exp^\perp_{L,M}(W^\T_{(\ta,\rh)}({X,L}))\not= W_{(\ta,\rh)}(X,L).
	\end{align*} 
	Thus we need the following fact, which says that up to a negligible error, the wedge neighborhoods are contained in the exponentiated image of wedge bundles:
	\begin{fact}\label{fctwi}
		With Assumption \ref{assumplxm} in mind, we have
		\begin{align*}
			W_{(\ta,\rh)}(X,L)\s \exp^\perp_{L,M}\Big(W^\T_{\big(\big(1+O(\rr+\rh)\big)\ta,\big(1+O(\rr+\rh)\big)\rh\big)}(X,L)\Big), 
		\end{align*}where the term $O(\cdot)$ depends geometrically on $X,L,M,$ and the ambient metric. 
	\end{fact}
	\begin{proof}
		The proof is a direct calculation in the double Fermi coordinate adapted to $(L,X,M).$ At a point $p,$ let $v\in\T_p^\perp (L,M).$ Taylor's theorem with remainders gives
		\begin{align*}
			\fd_M(\exp^\perp_{L,M}v,\pi_{X,M}(\exp^\perp_{L,M}v))=&(1+O(\rr+\rh))\no{v_{X,M}},\\
			\fd_X(\pi_{X,M}(\exp^\perp_{L,M}v),\pi_{L,X}\circ \pi_{X,M}(\exp^\perp_{L,M}v))=&(1+O(\rr+\rh))\no{v_{L,X}},
		\end{align*}and thus
		\begin{align*}
			\frac{\fd_M(\exp^\perp_{L,M}v,\pi_{X,M}(\exp^\perp_{L,M}v))}{\fd_X(\pi_{X,M}(\exp^\perp_{L,M}v),\pi_{L,X}\circ \pi_{X,M}(\exp^\perp_{L,M}v))}=&(1+O(\rr+\rh))\frac{\no{v_{X,M}}}{\no{v_{L,X}}},
		\end{align*}with the term $O(\cdot)$ depending geometrically on $X,L,M$ and the ambient metric. 
	\end{proof}

	\subsection{Proving $\Pi_{X,L}$ is $d$-area-non-increasing}	
	Now with the basic definitions and facts in hand, we are ready to achieve Step (\ref{im3}). Recall Assumption \ref{assumpdfc} and adopt a double Fermi coordinate chart of radius $(\rh,\rh,\rh)$ centered at $p\in\bff$.
	
	Our goal in this section is to prove that the retraction $\Pi_{X,L}$ defined in Definition \ref{defnpl} on the minimal submanifold $X$ obtained in Lemma \ref{localti} is $d$-area-non-increasing as defined in Definition \ref{dani} with respect to  any ambient metric $g\in \Om_h$, with $\Om_h$ defined in Lemma \ref{localti}.
	
	For a point $q=(l,x,z)$, by (\ref{pxlcdnt}), for $\frac{|z|}{|x|}> t_c$ (defined in (\ref{defntc})), we have
	\begin{align*}
		\Pi_{X,L}(l,x,z)=(l,0,0).
	\end{align*}
	Thus, any $d$-dimensional simple vector is mapped to $0$ and $\Pi_{X,L}$ is automatically $d$-area-non-increasing for $\frac{|z|}{|x|}> t_c.$
	
	The set $\{\frac{|z|}{|x|}=t_c\}$ is where $\Pi_{X,L}$ is not continuously differentiable and thus irrelevant to discussing the area-non-increasing property.
	
	Thus, we only have to focus on the points of $q=(l,x,z)$ with \begin{align*}
		\frac{|z|}{|x|}<t_c.
	\end{align*} Then the kernel $\ker\ed\Pi_{X,L}$ has dimension precisely $c$.  Consequently $\ed\Pi_{X,L}$ smoothly splits the tangent bundle $\T M$ to $M$ restricted to $B_\rr(X,M)$ into a sum of a $d$-dimensional and a $c$-dimensional orthogonal subspaces
	\begin{align*}
		\kxn\oplus	\kxp.
	\end{align*}
	At any point $q,$ the orthogonal projection onto $\kxn$ in $\T_pM$ is $d$-area-non-increasing. Thus, to verify $\Pi_{X,L}$ is $d$-area-non-increasing, it suffices to show that $\ed\Pi_{X,L}$ maps $\kxn$ to a $d$-vector of length less than $1.$ In other words we want to show that
	\begin{claim}
		At any point $q\in B_\rr(X,M)$, $\Pi_{X,L}$ is $d$-area-non-increasing.
	\end{claim}
And the above claim follows from
	\begin{claim}\label{pxl} Recall Assumption \ref{assumpdfc}. In any double Fermi coordinate chart centered at $p\in\bff$, at any point $q=(l,x,z)$ with $\frac{|z|}{|x|}<t_c$, let $v_1,\dots,v_d$ be a basis of $\kxn$, then we have
		\begin{align*}
			\frac{\no{\ed\Pi_{X,L}(v_1)\w\dots\w \ed\Pi_{X,L}(v_d)}^2}{\no{v_1\w\dots\w v_d}^2}	\le 1, 
		\end{align*}and equality holds only on $X.$
	\end{claim}
	The subsection will be devoted to the proof of Claim \ref{pxl}.
	\subsubsection{The ideas of the calculations}\label{sss}
	The proof reduces to careful calculations in double Fermi coordinate charts defined in Definition \ref{defndfermi}. However, before delving into the detailed calculations, it is worthwhile to give a quick illustration of why the calculation works.
	
	The basis we plan to use is obtained as follows. By Lemma \ref{dker}, we already have a complementary subspace $W$ to $\kxp$. The complementary subspace $W$ is spanned by the following vectors $w_1,\dots,w_d.$
	
	For $1\le j\le d-c$, we set $$w_j=\pd_{l_j}.$$ From ${d-c+1}\le j\le{d-1},$ we let $w_j\in \T_qM$ be $c-1$ different orthonormal vectors in $|\cdot|$ in the $x$-directions that spans the orthogonal complement of $\pd_x$ in $|\cdot|.$ In the notation of Lemma \ref{dker}, $w_{d-c+1},\dots,w_{d-1}$ spans the space of $\nu^x$ vectors. Finally we let 
	\begin{align*}
		w_d=\ccxy \pd_x+\csxy\pd_z.
	\end{align*}
	We emphasize that the coordinate quadratic norm $|\cdot|$ in general does not agree with the ambient metric $g$. Thus, the complementary space $W$ to $\kxp$ is not orthogonal to $\kxp$. Instead, our basis $v_1,\dots,v_d$ will be defined as the images of the orthogonal projections of $w_1,\dots,w_d$ on to $\kxn$.
	
	Our calculation proceeds in five steps, where each term $O(\cdot)$ depends geometrically on $L,X,M$ and the ambient metric:
	\begin{enumerate}
		\item First we show that
		\begin{align*}
			\no{v_1\w\dots\w v_d}^2=\no{w_1\w\dots\w w_d}^2(1+O(|z|^2)).
		\end{align*}
		\item Next we show that
		\begin{align*}
			\no{w_1\w\dots\w w_d}^2=(1+O(|z|^2))\bigg(\ccxy^2+\csxy^2\bigg)\sqrt{\det(l,x)}(l,x,0).
		\end{align*}Here $\sqrt{\det(l,x)}$ is the coefficient of the volume form of $X$ in the double Fermi coordinate chart, i.e, $$\dvol_X=\sqrt{\det(l,x)}dl_1\w\dots\w dl_{d-c}\w dx_1\w\dots\w dx_c.$$ The symbol $\sqrt{\det(l,x)}(l,x,0)$ means evaluating this coefficient at $(l,x,0).$ Precisely at this step we use the fact that the mean curvature of $X$ is $0,$ i.e., $X$ being a minimal submanifold.
		\item Next, we show that
		\begin{align*}
			\no{\ed\Pi_{X,L}(v_1\w\dots\w v_d)}^2=\bigg(\ccxy^2+\csxy^2\bigg)^2\sqrt{\det(l,x)}\bigg(l,\cgxyc x,0\bigg)
		\end{align*}
		\item Then we compare the values of $\sqrt{\det(l,x)}$
		\begin{align*}
			\frac{\sqrt{\det(l,x)}\bigg(l,\cgxyc x,0\bigg)}{\sqrt{\det(l,x)}(l,x,0)}=1+O\bigg(1-\cgxyc \bigg)|x|=1+O\bigg(\frac{|z|^2}{|x|}\bigg),
		\end{align*}
		\item Finally, using all of the above steps, and $O(|z|^2)=|x|^2O\big(\frac{|z|^2}{|x|^2}\big)$, we deduce from (\ref{lode}) that
		\begin{align*}
			&\frac{\no{\ed\Pi_{X,L}(v_1\w\dots\w v_d)}^2}{\no{v_1\w\dots\w v_d}^2}\\=&\bigg(1+|x|O\xy^2\bigg)\bigg(\ccxy^2+\csxy^2\bigg)\\
			\le&\bigg(1+O(|x|)\xy^2\bigg)\bigg(1-16 \xy^2\bigg)\\
			=&1-\bigg(\frac{1}{16}+O(|x|)\bigg)\xy^2.
		\end{align*}
		After some basic preparational work, our calculation will be carried out in the order of the above five steps.
	\end{enumerate}
	\subsubsection{Constructing $\{v_1,\dots,v_d\}$ from $\{w_1,\dots,w_d\}$}
	Now that we have introduced our ideas of calculations in Section \ref{sss}, we will formally define $\{w_1,\dots,w_d\}$ and $\{v_1,\dots,v_d\}$ to prepare for our calculations.
	
	Adopt the double Fermi coordinate chart in Definition \ref{defndfermi}. Fix an ambient metric $g\in \Om_h$ as obtained in Fact \ref{localti} and a point $q=(l,x,z)\in B_\de(X,M),$ with $\frac{|z|}{|x|}<t_c.$ 
	
	Choose $c-1$ vectors $\nu_1^x,\dots,\nu_{c-1}^x$ in $\T_{(l,x,z)}M$ with only non-vanishing components among the $x$-directions, such that for all $j,$ \begin{align*}
		\nu_j^x\pp\pd_x,\,|\nu_j^x|=1
	\end{align*} and for all $i\not=j$
	\begin{align*}
		\nu_i^x\pp\nu_j^x.
	\end{align*}
	In other words if we write $$\nu_j^x=a_j^i\pd_{x_i}.$$ Then at $\T_{(l,x,z)}M$, we have\begin{align}\label{defnvx}
		\sum_{i=1}^ca_j^ia_k^i=&\dd_{jk},\\\sum_{i=1}^c a_j^i\frac{x_i}{|x|}=&0.\label{defnvx1}
	\end{align}
	
	Now extend $\{\nu_1^x,\dots,\nu_{c-1}^x\}$ constantly in the double Fermi coordinate chart, i.e., setting $\nu_j^x=a_j^i\pd_{x_i}$ everywhere in our chart.
	
	Similarly, in $\T_{(l,x,z)}$, choose $c-1$ vectors $\nn_1^z,\dots,\nn_{c-1}^z$ with only non-vanishing components among $z$-directions such that for all $j,$ \begin{align*}
		\nn_j^z\pp\pd_z,\,|\nn_j^z|=1
	\end{align*} and for all $i\not=j$
	\begin{align*}
		\nn_i^z\pp\nn_j^z.
	\end{align*}In other words if we write $$\nn_j^z=b_j^i\pd_{z_i}.$$ Then at a point $q=(l,x,z)$, we have \begin{align*}
		\sum_{i=1}^cb_j^ib_k^i=&\dd_{jk},\\\sum_{i=1}^c b_j^i\frac{z_i}{|z|}=&0.	
	\end{align*}Now extend $\{\nn_1^z,\dots,\nn_{c-1}^z\}$ constantly in the double Fermi coordinate chart, i.e., setting $\nn_j^x=b_j^i\pd_{z_i}$ everywhere in our chart.
	
	Set \begin{align*}	\nu_0^x=&\ccxy\pd_x+\csxy\pd_z.	\\
		\nn_0^z=&-\csxy\pd_x+\ccxy\pd_z.
	\end{align*}
\begin{defn}\label{defnww}
		Define \begin{align*}
			&w_1=\pd_{l_1},\cdots,w_{d-c}=\pd_{l_{d-c}},\\&w_{d-c+1}=\nu_1^x,\cdots,w_{d-1}=\nu_{c-1}^x,\\&w_d=\nu_0^x.
		\end{align*}
	\end{defn} 
Let us recall  Lemma \ref{dker}.		In the notation we just adopted, we have
	\begin{align*}
		\kxp=\operatorname{span}\{\nn_0^z,\dots,\nn_{c-1}^z\}.
	\end{align*}
	Moreover, a  subspace \begin{align}\label{comp}
		W,
	\end{align} complementary to $\kxp$ is spanned by the vectors $w_1,\dots,w_d$.

	Unfortunately, $W=\operatorname{span}\{w_1,\cd,w_d\}$ is not orthogonal to $\kxp$, so we have to modify $W$ a bit.
	\begin{defn}
		Use $$\pi_\kxp$$ to denote the orthogonal projection of tangent spaces to $M$ onto $\kxp.$ 
	\end{defn}
	Set
	\begin{align}\label{defnvw}
		v_1=w_1-\pi_\kxp(w_1),\dots,v_d=w_d-\pi_\kxp(w_d).
	\end{align}
	Then applying Lemma \ref{dker}, we obtain the following fact: 
	\begin{fact}\label{imkxn}The vectors $\{v_1,\dots,v_d\}$ form a basis of $\kxn$ in $\T_qM$ with $q=(l,x,z),$ and
		\begin{align*}
			\ed\Pi_{X,L}(v_1\w\dots\w v_d)=\bigg(\ccxy^2+\csxy^2\bigg)(\pd_{l_1}\w\dots\w \pd_{l_{d-c}})\w(\nu_1^x\w\dots\w \nu_{c-1}^x)\w\pd_x.
		\end{align*}
	\end{fact}
	\subsubsection{Calculating $\pi_\kxn$}
	We want to write out $v_1,\dots,v_d$ explicitly and give a nice asymptotic of their expression in terms of $|x|,|z|.$ To achieve this, we first need to calculate the orthogonal projection $\pi_\kxp.$ Let us start with some basic calculations. Use $\dd_{ij}$ to denote the Kronecker delta symbol, i.e.,
	$\dd_{ij}=
	\begin{cases}
		1,&\textnormal{if }i=j,\\
		0,&\textnormal{if }i\not=j.
	\end{cases}
	$\begin{lem}\label{lemlxm}For $1\le i,j\le c,1\le k,m\le d-c$, at our point $q=(l,x,z)$, we have in $\T_{q}M$,
		\begin{align*}
			\min\{\frac{3}{4},\frac{3}{4}\no{\pd_x}^2\}\le&g(\nn_0^z,\nn_0^z)=\csxy^2\no{\pd_x}^2+\ccxy^2\le\max\{1,\no{\pd_x}^2\},\\	&g(\nn_i^z,\nn_j^z)=\dd_{ij}+O(|z|),\\	
			&g(\nn_i^z,\nn_0^z)=-\csxy g(\nn_i^z,\pd_x)=O\bigg(\frac{|z|^2}{|x|}\bigg),\\
			&g(\nu_i^x,\nn_0^z)=-\csxy g(\nu_i^x,\pd_x)=O\bigg(\frac{|z|^2}{|x|}\bigg),\\
			&g(\nu_0^x,\nn_0^z)=\ccxy\csxy(1-\no{\pd_x}^2)=O\bigg(\frac{|z|^2}{|x|}\bigg),\\
			&g(\nu_i^x,\nn_j^z)=O(|z|),\\
			&g(\nu_0^x,\nn_j^z)=O(|z|),\\
			&g(\pd_z,\nn_j^z)=0,\\
			&g(\pd_{l_k},\nn_0^z)=-\csxy g(\pd_{l_k},\pd_x)=O\bigg(\frac{|z|^2}{|x|}\bigg),\\
			&g(\pd_{l_k},\nn_j^z)=O(|z|),\\
			\min\{\frac{3}{4},\frac{3}{4}\no{\pd_x}^2\}\le&g(\nu_0^x,\nu_0^x)=\ccxy^2\no{\pd_x}^2+\csxy^2\le\max\{1,\no{\pd_x}^2\},\\
			&g(\nu^x_i,\nu_j^x)=\dd_{ij}+O(|x|+|z|),\\
			&g(\pd_z,\nu_j^x)=0,\\
			&g(\pd_x,\nu_j^x)=O(|z|),\\
			&g(\pd_{l_k},\nu_0^x)=\ccxy g(\pd_{l_k},\pd_x)= O(|z|)\\&g(\pd_x,\pd_z)=0,\\ &g(\pd_{l_k},\pd_x)=O(|z|),\\&g(\pd_{l_k},\nu_i^x)=O(|x|+|z|),\\
			&g(\pd_{l_k},\pd_{l_m})=O(|l|+|x|+|z|).
		\end{align*}Here the term $O(\cdot)$ only depends on $L,X,M$ and the ambient metric geometrically. To avoid repeating the above sentence, from now on in this section, whenever we use a Big $O,$ it is understood that the term $O(\cdot)$ depends geometrically on $L,X,Y,M$  and the ambient metric.
	\end{lem}
	\begin{proof}
		The proof reduces to calculations using Taylor's Theorem with remainders in double Fermi coordinate chart, plus Definition \ref{defng} and equation (\ref{lode}).
		
		First of all, recall Facts \ref{fctl}. \ref{fctx}, \ref{fcty}, especially $\pd_z\pp v\iff \pd_z\perp v$ holding everywhere and $\pd_x\pp v\iff\pd_x\perp v$ holding on $X.$
		
		At the center point $(0,0,0)$, we have
		\begin{align*}
			g(\pd_{l_k},\pd_{l_m})=\dd_{km}.
		\end{align*}
		Now Taylor's theorem with remainder gives that at $(l,x,z)$
		\begin{align*}
			g(\pd_{l_k},\pd_{l_m})(l,x,z)=&g(\pd_{l_k},\pd_{l_m})(0,0,0)+O(|l|+|x|+|z|)\\
			=&\dd_{km}+O(|l|+|x|+|z|)
		\end{align*}with the term $O(\cdot)$ bounded by the first derivatives of $	g(\pd_{l_k},\pd_{l_m})$ in our coordinate chart. However, the first derivatives of $	g(\pd_{l_k},\pd_{l_m})$ depend geometrically on $L,X,M$ and the ambient metric, since the double Fermi coordinate chart is essentially unique and depends geometrically on $L,X,M$ and the ambient metric. Thus, we deduce that the term $O(\cdot)$ depends geometrically on $L,X,M$ and the ambient metric.
		
		Now run Taylor's expansion theorem with remainder on all other combinations of coordinate vectors, we deduce that
		\begin{align*}
			&g(\pd_{l_k},\pd_{x_i})=O(|x|+|z|),\\
			&g(\pd_{l_k},\pd_x)=O(|z|),\\
			&g(\pd_x,\pd_x)=1+O(|z|),\\
			&g(\pd_{x_i},\pd_{x_j})=
			\dd_{ij}+O(|x|+|z|),\\
			&g(\pd_{l_k},\pd_{z_i})=O(|z|),\\
			&g(\pd_{l_k},\pd_z)=0,\\	&g(\pd_{x_i},\pd_{z_j})=O(|z|),\\
			&g(\pd_{x_i},\pd_z)=0,\\
			&g(\pd_z,\pd_z)=1,\\
			&g(\pd_{z_i},\pd_{z_j})=\dd_{ij}+O(|z|),
		\end{align*}where the term $O(\cdot)$ depends geometrically on $L,X,M$ and the ambient metric due to Fact \ref{dfcgd}.
Our lemma follows by direct calculations using the above list of metric expansions.
	\end{proof}
	\begin{defn}
		Define a $c\times c$ matrix $Z$, with the $ij$-th entry, $0\le i,j\le c-1,$ being 
		\begin{align*}
			Z_{ij}=g(\nn_i^z,\nn_j^z).
		\end{align*}We will also use $Z^{ij}$ to denote the $ij$-th entry of the inverse matrix $Z\m.$
	\end{defn}
	Let us derive an asymptotic for $Z_{ij}$ and $Z^{ij}.$ 
	
	Applying Lemma \ref{lemlxm}, we have deduced that the matrix $Z$ can be decomposed into blocks\begin{align*}
		Z=\begin{bmatrix}
			\csxy^2\no{\pd_x}^2+\ccxy^2 & O\bigg(\frac{|z|^2}{|x|}\bigg)\\
			O\bigg(\frac{|z|^2}{|x|}\bigg)&I_{c-1}+O(|z|)\\
		\end{bmatrix}.
	\end{align*}
	Here $I_{c-1}$ is the $(c-1)\times (c-1)$ identity matrix.
	
By estimate in the first line of Lemma \ref{lemlxm}, we see that $Z$ is invertible. By the inversion formula for a block matrix
	\begin{align*}
		\begin{bmatrix}
			A & B\\
			C&D\\
		\end{bmatrix}^{-1}=	\begin{bmatrix}
			(A-BD\m C)\m & -(A-BD\m C)\m BD\m\\
			-D\m C(A-BD\m C)\m &D\m+D\m C(A-BD\m C)\m BD\m\\
		\end{bmatrix}
	\end{align*}we can deduce the decomposition of $Z\m$ as a block matrix
	\begin{align}\label{zmmt}
		Z\m=\begin{bmatrix}
			\bigg(\csxy^2\no{\pd_x}^2+\ccxy^2\bigg)\m +O\bigg(\frac{|z|^4}{|x|^2}\bigg)& O\bigg(\frac{|z|^2}{|x|}\bigg)\\
			O\bigg(\frac{|z|^2}{|x|}\bigg)&I_{c-1}+O(|z|)+O\bigg(\frac{|z|^4}{|x|^2}\bigg)\\
		\end{bmatrix}.
	\end{align}
	Since we eventually will restrict our attention to a wedge neighborhood, we can assume that $\frac{|z|}{|x|}\le t_c$ (Definition \ref{defng}). We deduce that $O\bigg(\frac{|z|^4}{|x|^2}\bigg)=O(|z|^2)\le O(|z|),$ so $I_{c-1}+O(|z|)+O\big(\frac{|z|^4}{|x|^2}\big)=I_{c-1}+O(|z|).$

	Now direct calculation shows that the orthogonal projection $\pi_\kxp$ onto $\kxp$ has the following expression:\begin{fact}\label{defnkxp}
		We have
		\begin{align*}
			\pi_\kxp(w)=\sum_{0\le i,j\le c-1}Z^{ij}g(w,\nn_i^z)\nn_j^z
		\end{align*}
	\end{fact}
	\subsubsection{Calculating $v_1,\dots,v_d$}
	Let us calculate $\pi_\kxp(w_j)$ for all $j.$ 
	
	For $1\le m\le d-c,$ by Lemma \ref{lemlxm}, we have
	\begin{align*}
		&\pi_\kxp(\pd_{l_m})\\
		=&\sum_{0\le i,j\le c-1}Z^{ij}g(\pd_{l_m},\nn_i^z)\nn_j^z\\
		=&Z^{00}g(\pd_{l_m},\nn_0^z)\nn_0^z+\sum_{1\le j\le c-1}Z^{0j}g(\pd_{l_m},\nn_0^z)\nn_j^z+\sum_{1\le j\le c-1}Z^{j0}g(\pd_{l_m},\nn_j^z)\nn_0^z+\sum_{1\le i,j\le c-1}Z^{ij}g(\pd_{l_m},\nn_i^z)\nn_j^z\\
		=&O\bigg(\frac{|z|^2}{|x|}\bigg)\nn_0^z+\sum_{1\le j\le c-1}O\bigg(\frac{|z|^4}{|x|^2}\bigg)\nn_j^z+\sum_{1\le j\le c-1}O\bigg(\frac{|z|^3}{|x|}\bigg)\nn_0^z+\sum_{1\le i,j\le c-1}(\dd_{ij}+O(|z|))O(|z|)\nn_j^z\\
		=&O\bigg(\frac{|z|^2}{|x|}\bigg)\nn_0^z+\sum_{1\le j\le c-1}O(|z|)\nn_j^z.
	\end{align*}
	Similar calculations using Lemma \ref{lemlxm} show that for $1\le i\le c-1$, we have
	\begin{align*}
		\pi_\kxp(\nu_i^x)	=&O\bigg(\frac{|z|^2}{|x|}\bigg)\nn_0^z+\sum_{1\le j\le c-1}O(|z|)\nn_j^z,\\
		\pi_\kxp(\nu_0^x)	=&O\bigg(\frac{|z|^2}{|x|}\bigg)\nn_0^z+\sum_{1\le j\le c-1}O(|z|)\nn_j^z.
	\end{align*}
	Thus, we can conclude that for $1\le n\le d,$ we always have
	\begin{align}\label{vd}
		v_n=w_n+O\bigg(\frac{|z|^2}{|x|}\bigg)\nn_0^z+\sum_{1\le j\le c-1}O(|z|)\nn_j^z.
	\end{align}
	\subsubsection{Comparing $\no{v_1\w\dots\w v_d}^2$ to $\no{w_1\w\dots\w w_d}^2$}
	By definition, we have
	\begin{align*}
		\no{v_1\w\dots\w v_d}^2=\det([g(v_m,v_n)]),
	\end{align*}where we regard $[g(v_m,v_n)]$ as the matrix with the $m,n$-th entry being $g(v_m,v_n).$
	
	Thus, we want to compare $g(v_m,v_n)$ with $g(w_m,w_n).$
	
	Let us calculate $g(v_d,v_d)$ first. By (\ref{vd}), We have
	\begin{align*}
		&g(v_d,v_d)\\
		=&g\bigg(w_d+O\bigg(\frac{|z|^2}{|x|}\bigg)\nn_0^z+\sum_{1\le j\le c-1}O(|z|)\nn_j^z,w_d+O\bigg(\frac{|z|^2}{|x|}\bigg)\nn_0^z+\sum_{1\le j\le c-1}O(|z|)\nn_j^z\bigg)\\
		=&g(w_d,w_d)+O\bigg(\frac{|z|^2}{|x|}\bigg)g(w_d,\nn_0^z)+O(|z|)\sum_{1\le j\le c-1}g(w_d,\nn_j^z)\\
		&+O\bigg(\frac{|z|^4}{|x|^2}\bigg)g(\nn_0^z,\nn_0^z)+O(|z|^2)\sum_{1\le i,j\le c-1}g(\nn_i^z,\nn_j^z)+O\bigg(\frac{|z|^3}{|x|}\bigg)\sum_{1\le j\le c-1}g(\nn_j^z,\nn_0^z)\\
		=&g(w_d,w_d)+O\bigg(\frac{|z|^4}{|x|^2}\bigg)+O(|z|^2)+O\bigg(\frac{|z|^4}{|x|^2}\bigg)+O(|z|^2)+O\bigg(\frac{|z|^4}{|x|}\bigg)\\
		=&g(w_d,w_d)+O(|z|^2).
	\end{align*}
	Similar calculations shows that for any $1\le m,n\le d,$ we always have
	\begin{align*}
		g(v_m,v_n)=g(w_m,w_n)+O(|z|^2).
	\end{align*}
	Thus, we have
	\begin{align}
		&\no{v_1\w\dots\w v_d}^2\\
		=&\det([g(v_m,v_n)])\\
		=&\det\big([g(w_m,w_n)+O(|z|^2)]\big)\\
		=&\det([g(w_m,w_n)])\det(I_d+O(|z|^2))\\
		=&(1+O(|z|^2))\det([g(w_m,w_n)]).\label{detvd}
	\end{align}
	\subsubsection{Calculating $\no{w_1\w\dots\w w_d}^2$}
	Let us first calculate $\no{w_1\w\dots\w w_d}^2=\det([g(w_m,w_n)])$.
	
	Since $\pd_z$ is always perpendicular to $w_1,\dots,w_{d-1},\pd_x$, we have
	\begin{align}
		&\no{w_1\w\dots\w w_d}^2\\
		=&g(w_1\w\dots\w w_d,w_1\w\dots\w w_d)\\
		=&g\bigg(w_1\w\dots\w w_{d-1}\w\bigg(\ccxy\pd_x+\csxy\pd_z\bigg),w_1\w\dots\w w_{d-1}\w\bigg(\ccxy\pd_x+\csxy\pd_z\bigg)\bigg)\\
		=&\ccxy^2\no{w_1\w\dots\w w_{d-1}\w \pd_x}^2+\csxy^2\no{w_1\w\dots\w w_{d-1}\w\pd_z}^2.\label{wde}
	\end{align}
	At our point $q=(l,x,z)$, recall the definition of $\nu_j^x$, i.e, (\ref{defnvx}) and (\ref{defnvx1}). Since $\nu_j^x$ are of constant value throughout our coordinate chart and $\pd_x$ is the radial vector field in the $x$ coordinate, at any point $(l',ax,y')$ with $a\in \R$ we have
	\begin{align}\label{voleq}
		\nu_1^x\w\dots\w \nu_{c-1}^x\w \pd_x=\pm\pd_{x_1}\w\dots\w \pd_{x_c}.
	\end{align}
	This implies that
	\begin{align}\label{wd}
		\no{w_1\w\dots\w w_{d-1}\w \pd_x}^2=\no{\pd_{l_1}\w\dots\w \pd_{l_{d-c}}\w\pd_{x_1}\w\dots\w \pd_{x_c}}^2.
	\end{align}
	On the other hand, since $\pd_z\perp w_j$ for $1 \le j \le d-1$ and $\no{\pd_z}=1,$ using determinant of block matrix, we have
	\begin{align}\label{wdy}	\no{w_1\w\dots\w w_{d-1}\w\pd_z}^2
		=\no{w_1\w\dots\w w_{d-1}}^2.
	\end{align}
	Now we want to compare $\no{w_1\w\dots\w w_{d-1}}^2$ with $\no{w_1\w\dots\w w_{d-1}\w \pd_x}^2.$
	
	By the Laplace expansion of determinants, we have
	\begin{align*}
		&\no{w_1\w\dots\w w_{d-1}\w \pd_x}^2\\
		=&\no{\pd_x}^2\no{w_1\w\dots\w w_{d-1}}^2+(-1)^{1+d}g( w_2\w \dots\w w_{d-1}\w \pd_x,w_1\w\dots\w w_{d-1})g(w_1,\pd_x)\\&+\dots+(-1)^{d-1+d}g(w_1\w\dots\w w_{d-2}\w \pd_x,w_1\w\dots\w w_{d-1})g(w_{d-1},\pd_x).
	\end{align*}
	Using the Laplace expansion again, we have
	\begin{align*}
		&g( w_2\w \dots\w w_{d-1}\w \pd_x,w_1\w\dots\w w_{d-1})g(w_1,\pd_x)\\
		=&(-1)^{d-1+1}g(w_2\w\dots\w w_{d-1},w_2\w\dots\w w_{d-1})g(\pd_x,w_1)g(w_1,\pd_x)\\
		&+\dots+ (-1)^{d-1+d-1}g(w_2\w\dots\w w_{d-1},w_1\w\dots\w w_{d-2})g(\pd_x,w_{d-1})g(w_1,\pd_x)\\
		=&O(|z|^2).
	\end{align*}
	Similar calculation shows that
	\begin{align*}
		g(w_1\w\dots\w \ov{w_j}\w\dots\w w_{d-1}\w \pd_x,w_1\w\dots\w w_{d-1})g(w_j,\pd_x)=O(|z|^2),
	\end{align*}where $\ov{w_j}$ means omitting the term $w_j.$
	
	Thus, we can conclude 
	\begin{align*}
		\no{w_1\w\dots\w w_{d-1}\w\pd_x}^2=\no{\pd_x}^2\no{w_1\w\dots\w w_{d-1}}^2+O(|z|^2).
	\end{align*}
	In other words
	\begin{align}\label{wdonly}
		\no{w_1\w\dots\w w_{d-1}}^2=	\no{\pd_x}^{-2}	\no{w_1\w\dots\w w_{d-1}\w\pd_x}^2+O(|z|^2).
	\end{align}
	
	Now, combine (\ref{wdonly}), (\ref{wdy}), (\ref{wd}) into (\ref{wde}), and using $\no{\pd_x}^2=1+O(|z|)$ and (\ref{lode}), we deduce that
	\begin{align}
		&\no{w_1\w\dots\w w_d}^2\\
		=&\ccxy^2\no{\pd_{l_1}\w\dots\w \pd_{l_{d-c}}\w\pd_{x_1}\w\dots\w \pd_{x_c}}^2\\&+\csxy^2\big(\no{\pd_x}^{-2}\no{\pd_{l_1}\w\dots\w \pd_{l_{d-c}}\w\pd_{x_1}\w\dots\w \pd_{x_c}}^2+O(|z|^2)\big)\\
		=&\bigg(\ccxy^2+\csxy^2\bigg)\no{\pd_{l_1}\w\dots\w \pd_{l_{d-c}}\w\pd_{x_1}\w\dots\w \pd_{x_c}}^2\\&+O\bigg(\frac{|z|^2}{|x|^2}\bigg)O(|z|)\no{\pd_{l_1}\w\dots\w \pd_{l_{d-c}}\w\pd_{x_1}\w\dots\w \pd_{x_c}}^2\\	=&\bigg(1+O\bigg(\frac{|z|^3}{|x|^2}\bigg)\bigg)\bigg(\ccxy^2+\csxy^2\bigg)\no{\pd_{l_1}\w\dots\w \pd_{l_{d-c}}\w\pd_{x_1}\w\dots\w \pd_{x_c}}^2\label{svd}.
	\end{align}
	Now bring (\ref{svd}) into (\ref{detvd}), we deduce that
	\begin{align}
		&\no{v_1\w\dots\w v_d}^2\\
		=&\bigg(1+O\bigg(|z|^2+\frac{|z|^3}{|x|^2}\bigg)\bigg)\bigg(\ccxy^2+\csxy^2\bigg)\\&\cdot\no{\pd_{l_1}\w\dots\w \pd_{l_{d-c}}\w\pd_{x_1}\w\dots\w \pd_{x_c}}^2(l,x,z).\label{vdult}
	\end{align}Here the last $(l,x,z)$ is to signify that $\no{\pd_{l_1}\w\dots\w \pd_{l_{d-c}}\w\pd_{x_1}\w\dots\w \pd_{x_c}}^2$ is evaluated at $(l,x,z)$.
	\subsubsection{Calculating $\no{\ed\Pi_{X,L}(v_1\w\dots\w v_d)}^2$}
	By Fact \ref{imkxn}, we have
	\begin{align*}
		&\no{\ed\Pi_{X,L}(v_1\w\dots\w v_d)}^2\\
		=&\bigg(\ccxy^2+\csxy^2\bigg)^2\no{\pd_{l_1}\w\dots\w \pd_{l_{d-c}}\w \nu_1^x\w\dots\w \nu_{c-1}^x\w\pd_x}^2.
	\end{align*}
	By (\ref{voleq}), we deduce that
	\begin{align}
		&\no{\ed\Pi_{X,L}(v_1\w\dots\w v_d)}^2\\
		=&\bigg(\ccxy^2+\csxy^2\bigg)^2\no{\pd_{l_1}\w\dots\w \pd_{l_{d-c}}\w\pd_{x_1}\w\dots\w \pd_{x_c}}^2\bigg(l,\cgxyc x,0\bigg)\label{pult}.
	\end{align}Here the last $\bigg(l,\cgxyc x,0\bigg)$ is to signify that we evaluate $\no{\pd_{l_1}\w\dots\w \pd_{l_{d-c}}\w\pd_{x_1}\w\dots\w \pd_{x_c}}^2$ at $\bigg(l,\cgxyc x,0\bigg)$.
	\subsubsection{$(l,x)$ volume element}
	At this stage, from (\ref{vdult}) and (\ref{pult}) it is clear that in order to prove $\Pi_{X,L}$ is $d$-area-non-increasing, the key is to estimate 
	\begin{align}
		\no{\pd_{l_1}\w\dots\w \pd_{l_{d-c}}\w\pd_{x_1}\w\dots\w \pd_{x_c}}^2.\label{lxvol}
	\end{align}
	\begin{defn}
		We call the above expression (\ref{lxvol}) the $(l,x)$ volume elements.
	\end{defn}
	We use the name $(l,x)$ volume element because the volume form of the $d$-dimensional submanifolds $\{z=a\}$ with $a\in \R^{c}$ is precisely
	\begin{align*}
		\dvol_{\{z=a\}}=
		\no{\pd_{l_1}\w\dots\w \pd_{l_{d-c}}\w\pd_{x_1}\w\dots\w \pd_{x_c}}dl_1\w\dots\w dl_{d-c}\w dx_1\w\dots\w dx_c.
	\end{align*}
	Estimating the $(l,x)$ volume element involves two steps. 
	
	First, let us prove that
	\begin{claim}
		We have
		\begin{align}\label{keylx}
			\no{\pd_{l_1}\w\dots\w \pd_{l_{d-c}}\w\pd_{x_1}\w\dots\w \pd_{x_c}}^2(l,x,z)=(1+O(|z|^2))\no{\pd_{l_1}\w\dots\w \pd_{l_{d-c}}\w\pd_{x_1}\w\dots\w \pd_{x_c}}^2(l,x,0).
		\end{align}
	\end{claim}
	The above claim is the only place in this section that uses the fact of $X$ being a minimal submanifold.
	\begin{proof}
		Again, we want to use Taylor's theorem with remainders. To simplify our notation, set
		\begin{align*}
			&u_1=\pd_{l_1},\dots,u_{d-c}=\pd_{l_{d-c}},\\
			&u_{d-c+1}=\pd_{x_1},\dots,u_d=\pd_{x_c}.
		\end{align*}
		At $(l,x,0)$, using formulas for derivatives of determinants, we can calculate
		\begin{align}
			&\pd_{z_k}		\no{\pd_{l_1}\w\dots\w \pd_{l_{d-c}}\w\pd_{x_1}\w\dots\w \pd_{x_c}}^2\\
			=&\no{\pd_{l_1}\w\dots\w \pd_{l_{d-c}}\w\pd_{x_1}\w\dots\w \pd_{x_c}}^2\cdot\tr\big([g(u_m,u_n)]\m [\pd_{z_k}g(u_{m'},u_{n'})]		\big)			.\label{zmn}
		\end{align}Here $[g(u_m,u_n)]$ means the matrix with the $mn$-th entry being $g(u_m,u_n)$ and $[g(u_m,u_n)]\m$ is the inverse matrix of $[g(u_m,u_n)]$.
		
		We have
		\begin{align}
			&\pd_{z_k}g(u_m,u_n)\\
			=&g(\na_{\pd_{z_k}}u_m,u_n)+g(u_m,\na_{\pd_{z_k}}u_n)\\
			=&g(\na_{u_m}{\pd_{z_k}},u_n)+g(u_m,\na_{u_n}{\pd_{z_k}})\label{excd}\\
			=&u_m\big(g({\pd_{z_k}},u_n)\big)-g({\pd_{z_k}},\na_{u_m}u_n)+u_n\big(g(u_m,{\pd_{z_k}})\big)-u_ng(\na_{u_n}u_m,{\pd_{z_k}})\label{vanin}\\
			=&-2g(\A_X(u_m,u_n),\pd_{z_k}),\label{amn}
		\end{align}where $\A_X$ is the second fundamental form tensor of $X$ inside $M$. At (\ref{excd}), we have used the fact that coordinate vector fields have zero lie bracket, so $0=[\pd_{z_k},u_m]=\na_{\pd_{z_k}}u_m-\na_{u_m}\pd_{z_k}$, etc. At (\ref{amn}), we have used the fact that $\pd_{z_k}$ is always orthogonal to $u_j$ on $X.$
		
		Thus, from (\ref{amn}), we deduce that
		\begin{align}
			&	\tr\big([g(u_m,u_n)]\m [\pd_{z_k}g(u_{m'},u_{n'})]\\	=&-2\tr\bigg([g(u_m,u_n)]\m [g(\A_X(u_{m'},u_{n'}),\pd_{z_k})]\bigg)\\
			=&-2g(\tr_X\A_X(\cdot,\cdot),\pd_{z_k})\\
			=&-2g(\hh_X,\pd_{z_k}),\label{hzk}
		\end{align}
		where $\tr_X$ means the Riemannian trace of the tensor $\A_X$ with respect to the metric on $X$ and $\hh_X$ is the mean curvature of $X$ inside $M.$
		
		Since $X$ is a minimal submanifold, the mean curvature of $X$ vanishes, i.e., $\hh_X=0.$	Thus, bringing (\ref{hzk}) to (\ref{zmn}), we deduce that
		\begin{align*}
			\pd_{z_k}		\no{\pd_{l_1}\w\dots\w \pd_{l_{d-c}}\w\pd_{x_1}\w\dots\w \pd_{x_c}}^2=0.			
		\end{align*}
		Since $l,x$ and $k$ are arbitrary, we deduce that on $X$ all the first derivatives of $$		\no{\pd_{l_1}\w\dots\w \pd_{l_{d-c}}\w\pd_{x_1}\w\dots\w \pd_{x_c}}^2
		$$ along $z$-directions vanish. Thus, applying Taylor's theorem with remainders, we deduce that
		\begin{align*}
			\no{\pd_{l_1}\w\dots\w \pd_{l_{d-c}}\w\pd_{x_1}\w\dots\w \pd_{x_c}}^2(l,x,z)=\no{\pd_{l_1}\w\dots\w \pd_{l_{d-c}}\w\pd_{x_1}\w\dots\w \pd_{x_c}}^2(l,x,0)+O(|z|^2).	
		\end{align*}Since $\no{\pd_{l_1}\w\dots\w \pd_{l_{d-c}}\w\pd_{x_1}\w\dots\w \pd_{x_c}}^2
		$ is bounded below geometrically from $0$, we are done.
	\end{proof}
	Next we shall prove that
	\begin{claim}
		We have
		\begin{align}&\no{\pd_{l_1}\w\dots\w \pd_{l_{d-c}}\w\pd_{x_1}\w\dots\w \pd_{x_c}}^2\bigg(l,\cgxyc x,0\bigg)\\=&\bigg(1+O\bigg(\frac{|z|^2}{|x|}\bigg)\bigg)\no{\pd_{l_1}\w\dots\w \pd_{l_{d-c}}\w\pd_{x_1}\w\dots\w \pd_{x_c}}^2(l,x,0).\label{keypg}
		\end{align}
	\end{claim}
	\begin{proof}
		Define
		\begin{align*}
			F(\tau)=\no{\pd_{l_1}\w\dots\w \pd_{l_{d-c}}\w\pd_{x_1}\w\dots\w \pd_{x_c}}^2\bigg(l,\tau \frac{x}{|x|},0\bigg),
		\end{align*}where here $\frac{x}{|x|}$ uses the $x$ in $q=(l,x,z).$ 
		Then we have $F'=O(1)$. Thus,
		\begin{align*}
			F(|x|)-F\bigg(\cgxyc|x|\bigg)=O\bigg(1-\cgxyc\bigg)|x|.
		\end{align*}
		By the factorization,\begin{align*}
			1-\cgxyc=\frac{1-\cgxy}{\sum_{j=0}^{c-1}\mathbf{g}^{\frac{j}{c}}\left(\frac{|z|}{|x|}\right)},
		\end{align*}
		 from Definition \ref{defng}, we deduce that
		\begin{align*}
			1-\cgxyc\le 1-\cgxy=\frac{2c^2-9}{8}\frac{|z|^2}{|x|^2}.
		\end{align*}
		This implies that
		\begin{align*}
			F(|x|)-F\bigg(\cgxyc|x|\bigg)=O\bigg(\frac{|z|^2}{|x|}\bigg). 
		\end{align*}We are done.
	\end{proof}
	\subsubsection{Wrapping up the proof}
	Combining (\ref{vdult}) (\ref{pult}) (\ref{keylx}) (\ref{keypg}), we deduce that
	\begin{align*}
		&\frac{\no{\Pi_{X,L}(v_1\w\dots\w v_d)}^2}{\no{v_1\w\dots\w v_d}^2}\\
		=&\bigg(\ccxy^2+\csxy^2\bigg)\frac{\bigg(1+O(\frac{|z|^2}{|x|})\bigg)\no{\pd_{l_1}\w\dots\w \pd_{l_{d-c}}\w\pd_{x_1}\w\dots\w \pd_{x_c}}^2(l,x,0)}{(1+O(|z|^2))\bigg(1+O\bigg(|z|^2+\frac{|z|^3}{|x|^2}\bigg)\bigg)\no{\pd_{l_1}\w\dots\w \pd_{l_{d-c}}\w\pd_{x_1}\w\dots\w \pd_{x_c}}^2(l,x,0)}\\
		&=\bigg(\ccxy^2+\csxy^2\bigg)\bigg(1+O\bigg(|z|^2+\frac{|z|^2}{|x|}+\frac{|z|^3}{|x|^2}\bigg)\bigg).
	\end{align*}
	By (\ref{lode}), we deduce that
	\begin{align}
		&	\frac{\no{\Pi_{X,L}(v_1\w\dots\w v_d)}^2}{\no{v_1\w\dots\w v_d}^2}\le \bigg(1-\frac{1}{16}\bigg(\frac{|z|^2}{|x|^2}\bigg)\bigg)\bigg(1+(|x|^2+|x|+|z|)O\bigg(\frac{|z|^2}{|x|^2}\bigg)\bigg)\\
		=&1-\bigg(\frac{1}{16}+O(|x|)\bigg)\bigg(\frac{|z|^2}{|x|^2}\bigg).\label{cmspi}
	\end{align}
Thus, if $O(|x|)$ is small enough, then we have
	\begin{align}\label{est1}
		\frac{\no{\Pi_{X,L}(v_1\w\dots\w v_d)}^2}{\no{v_1\w\dots\w v_d}^2}\le 1,
	\end{align}and
	\begin{align}\label{est2}
		\frac{\no{\Pi_{X,L}(v_1\w\dots\w v_d)}^2}{\no{v_1\w\dots\w v_d}^2}=1\iff |z|=0.
	\end{align}
		
	This smallness of $O(|x|)$ can be achieved as follows. 
	
	In our setting, when $2[\Si]\not=0,$ by making $\ees$  small in Lemma \ref{localti}, $X$ can be made arbitrarily close to \begin{align}\label{pm4}
		\Ga\m({B_\eta^c\times\{0\}^c}\times S^{d-c})
	\end{align} geometrically. By the bullet in Lemma \ref{nsmin}, our base metric $h$ restricts to the standard product metric on \begin{align}\label{pm5}
	\Ga\m(B_\eta^{2c}\times S^{d-c}).
	\end{align} This implies that 
	\begin{align}\label{splitt}
		\Pi_{B_\eta^c\times\{0\}^c\times S^{d-c},\{0\}^{2c}\times S^{d-c}}(u,v)=\left(\Pi_{\R^c\times\{0\}^c,\{0\}^{2c}}(u),\pi_{S^{d-c}}(v)\right),
	\end{align}for $(u,v)\in \left(B_\eta^c\times\{0\}^c\right)\times S^{d-c}$,
	where $\pi_{S^{d-c}}$ denotes the projection on to the factor $S^{d-c}$ in $\Ga\m(B_\eta^{2c}\times S^{d-c}).$
	
Now if we run our calculations again in $h,$ then all contributions from curvatures disappear due to the splitting (\ref{splitt}), and it is very quick to arrive that in our base metric $h$, we have term $O(\cdot)$ in (\ref{cmspi}) equal to $0.$ 
	 Since the term $O(\cdot)$ in (\ref{cmspi}) depends geometrically on $L,X,M$ and the ambient metric, we deduce that as long as we set the $\ees$ in Lemma \ref{localti} small enough, we have $|O(|x|)|<\frac{1}{32}$ in (\ref{cmspi}). When $2[\Si]=0,$ we add a $\pm$ sign to the (\ref{pm4}), (\ref{pm5}) and (\ref{splitt}), meaning positively or negatively oriented. However, the argument  goes through unchanged, so we can still deduce that $|O(|x|)|<\frac{1}{32}$ in (\ref{cmspi}).
	
Our choice of point $q=(l,x,z)$ with $\frac{|z|}{|x|}<t_c$ is arbitrary, so the estimates (\ref{est1}) and (\ref{est2}) holds for all $q=(l,x,z)$ with $\frac{|z|}{|x|}<t_c$. Now recall our Assumption \ref{assumpdfc}. Our choice of base point $p$ in the finite set $F$ is also arbitrary. Setting $$\eta=\frac{\rr}{4},$$ we deduce that $\Pi_{X,L}$ is a $d$-area-non-increasing as claimed. 	To sum it up, we have prove the following
	\begin{fact}\label{fctpx}
		Recall Lemma \ref{localti}.	For any metric $g\in \Om_h$, we have the following
		\begin{enumerate}
			\item In $U_{\frac{\rr}{4}}(\sing\sdc),$ the retraction $\Pi_{X,L}$ in Definition \ref{defnpl} is defined everywhere,
			\item $\Pi_{X,L}$ is $d$-area-non-increasing, wherever $\Pi_{X,L}$ is continuously differentiable and $\ees$ in Lemma \ref{localti} is small enough,
			\item $\Pi_{X,L}$ maps unit simple $d$-vectors to unit simple $d$-vectors only when restricted to $X.$
		\end{enumerate}
	\end{fact}
	\begin{rem}\label{remim}
		Another way to achieve the smallness of $O(|x|)$ in (\ref{cmspi}) works in more general settings. We can first make the radius $\eta$ in $U_{\eta}(\sing\sdc)$ small in Lemma \ref{localti}. This makes the maximum of $|x|$ small. Then make $\ees$ small in Lemma \ref{localti}, we can always achieve $|O(|x|)|<\frac{1}{32}$ by geometric dependence on $L,X,M$ and the ambient metric. However, this double parameter smallness argument is not needed for our case.
	\end{rem}
	\subsection{Proving that $\pi_{Y,L}$ is $d$-area-non-increasing}
	The same argument as in the previous subsection gives.
	\begin{fact}\label{fctpy}
		Recall Lemma \ref{localti}.	For any metric $g\in \Om_h$, we have
		\begin{enumerate}
			\item In $U_{\frac{\rr}{4}}(\sing\sdc),$ the retraction $\Pi_{Y,L}$ in Definition \ref{defnply} is defined everywhere,
			\item $\Pi_{Y,L}$ is $d$-area-non-increasing, wherever $\Pi_{Y,L}$ is continuously differentiable and $\ees$ in Lemma \ref{localti} is small enough,
			\item $\Pi_{Y,L}$ maps unit simple $d$-vectors to unit simple $d$-vectors only when restricted to $Y.$
		\end{enumerate}
	\end{fact}
	\subsection{Proving that $\Pi_{X,L}\du\dvol_X+\Pi_{Y,L}\du\dvol_Y$ is a calibration}
	First, let us show that \begin{fact}\label{fctxyad}
		$$\Pi_{X,L}\du\dvol_X+\Pi_{Y,L}\du\dvol_Y$$  is a $d$-form with Lipschitz antiderivative (Definition \ref{defnlip}).
	\end{fact}
	\begin{proof}
		Since $X$ is a smooth submanifold with boundary, we have $H_d(X,\Z)=0.$ By the universal coefficient theorem, we have $H^d(X,\R)=0.$ By de Rham's theorem, we deduce that $\dvol_X$ is an exact form. 
		
		Since $\Pi_{X,L}$ is Lipschitz, it is straightforward to verify that $\Pi_{X,L}\dvol_X$	is a form with Lipschitz anti-derivative. 
		
		Replace $X$ with $Y,$ we deduce that $\Pi_{Y,L}\du\dvol_Y$ is a form with Lipschitz antiderivative. Thus, the sum $\Pi_{X,L}\du\dvol_X+\Pi_{Y,L}\du\dvol_Y$ is also a form with Lipschitz antiderivative.
	\end{proof}
	Next, let us show that
	\begin{fact}\label{fctxynz}
		On $U_{\frac{\rr}{4}}(\sing\sdc)$, the form $
		\Pi_{X,L}\du \dvol_{X}+\Pi_{Y,L}\du\dvol_Y,
		$ is nonzero only when restricted to interior of 
		$$		W_{(\arctan t_c,\frac{\rr}{2})}(X,L)\cap U_{\frac{\rr}{4}}(\sing\sdc),$$ and $$		W_{(\arctan t_c,\frac{\rr}{2})}(Y,L)\cap U_{\frac{\rr}{4}}(\sing\sdc),$$ where it equals $\Pi_{X,L}\du \dvol_{X}$ and $\Pi_{Y,L}\du\dvol_Y$ respectively.
	\end{fact}
	\begin{proof}Recall (\ref{defntc}), we always have $\arctan t_c<\frac{\pi}{4}.$
		With $\ees$ in  Lemma \ref{localti} small enough, we can assume that the intersection angle between $X,Y$ satisfies
		\begin{align}\label{taineq}
		\frac{\pi}{2}\ge	\Ta(X,Y)>\frac{\frac{\pi}{2}+2\arctan t_c}{2}.
		\end{align}
		In our base metric $h$ restricted to $U_\rr(\sing\sdc),$ which is just the product of flat metric and standard metric on unit $(d-c)$-sphere, we have
		\begin{align*}
			&	W_{(\arctan t_c,\frac{\rr}{2})}(B_1^{c}\times\{0\}^c\times S^{d-c},\{0\}^{2c}\times S^{d-c})\\&=\exp^\perp W^\T_{(\arctan t_c,\frac{\rr}{2})}(B_1^{c}\times\{0\}^c\times S^{d-c},\{0\}^{2c}\times S^{d-c}),
		\end{align*}
		where $W^\T$ is the wedge bundle defined in Definition \ref{defnwb}.
		
		By geometric dependence in Fact \ref{fctwi}, we deduce that with $\ees$ in Lemma \ref{localti} small, we have
		\begin{align}\label{wx}
			W_{(\arctan t_c,\frac{\rr}{2})}(X,L)\s W^\T_{(\frac{\frac{\pi}{4}+\arctan t_c}{2},{\rr})}(X,L).
		\end{align}
		Similarly, we have
		\begin{align}\label{wy}
			W_{(\arctan t_c,\frac{\rr}{2})}(Y,L)\s W^\T_{(\frac{\frac{\pi}{4}+\arctan t_c}{2},{\rr})}(Y,L).
		\end{align}
		By (\ref{wx}), (\ref{wy}), (\ref{taineq}) and Lemma \ref{lemwi}, we deduce that
		\begin{align}\label{wxyi}
			W_{(\arctan t_c,\frac{\rr}{2})}(X,L)\cap	W_{(\arctan t_c,\frac{\rr}{2})}(Y,L)=L.
		\end{align}

		By Definition of $\Pi_{X,L}$ (\ref{defnpl}), $$\Pi_{X,L}\du\dvol_X\not=0$$ only when restricting to in the interior of wedge neighborhoods (Definition \ref{defnw})
		\begin{align*}
			W_{(\arctan t_c,\frac{\rr}{2})}(X,L)\cap U_{\frac{\rr}{4}}(\sing\sdc).
		\end{align*}
		Similarly,
		$$\Pi_{Y,L}\du\dvol_Y\not=0$$ only when restricting to in the interior of wedge neighborhoods (Definition \ref{defnw})
		\begin{align*}
			W_{(\arctan t_c,\frac{\rr}{2})}(Y,L)\cap U_{\frac{\rr}{4}}(\sing\sdc).
		\end{align*}
		Thus, applying (\ref{wxyi}), we are done.
	\end{proof}
	Now we will arrive at our final conclusion of this subsection.
	\begin{fact}\label{fctxyuni}
		$\Pi_{X,L}\du\dvol_X+\Pi_{Y,L}\du\dvol_Y$ has comass at most $1$ and calibrates only $aX+bY$ with $a,b\in\Z_+.$
	\end{fact}
	\begin{proof}
		In the interior of $		W_{(\arctan t_c,\frac{\rr}{2})}(X,L)\cap U_{\frac{\rr}{4}}(\sing\sdc),
		$ $\Pi_{X,L}$ is $d$-area-non-increasing by Fact \ref{fctpx}, and sends unit simple $d$-vectors to unit simple $d$-vectors only when restricted to $X$.  By Fact \ref{fctani}, this implies that  $\Pi_{X,L}\du\dvol_X$ is a calibration form and calibrates only natural number multiples of $X$ by Lemma \ref{const} applied to the interior of $X.$
		
		Repeating the argument with $X$ replaced by $Y,$ we deduce that  $\Pi_{Y,L}\du\dvol_Y$ is a calibration form and calibrates only integer multiples of $Y.$ 
	\end{proof}
	\subsection{Concluding the proof of Lemma \ref{nst}}
		It is straight forward to verify that $\Pi_{X,L}\du \dvol_X+\Pi_{Y,L}\du\dvol_Y$ is almost continuous (Definition \ref{defnlip}) on $X+Y.$ By $\pd(X+Y)=\pd T|_{U_{\frac{\rr}{4}}(\sing\sdc)}$, Fact \ref{fctxyuni} and Lemma \ref{fcal}, we deduce that 
		\begin{align}\label{txyineq}
			\ms(X+Y)\le\ms	\big(	T|_{U_{\frac{\rr}{4}}(\sing\sdc)}\big).
		\end{align} 
	
	Recall Lemma \ref{localti}. Since $H_d(U_{\frac{\rr}{4}}(\sing\ssc),\Z)=0,$ we deduce that
	\begin{align*}
		T|_{U_{\frac{\rr}{4}}(\sing\sdc)\cp}+X+Y
	\end{align*}
	is homologous to $T.$

Since $T$ is area-minimizing, we must  have $\ms(X+Y)\ge\ms	\big(	T|_{U_{\frac{\rr}{4}}(\sing\sdc)}\big).$ Combined with (\ref{txyineq}), this implies that $\ms(X+Y)=\ms	\big(	T|_{U_{\frac{\rr}{4}}(\sing\sdc)}\big),$ By Lemma \ref{fcal}, $T|_{U_{\frac{\rr}{4}}(\sing\sdc)}$ is also calibrated by $\Pi_{X,L}\du \dvol_X+\Pi_{Y,L}\du\dvol_Y$.

	By Fact \ref{fctxyuni} and $\pd(		T|_{U_{\frac{\rr}{4}}(\sing\sdc)})=\pd(X+Y)$, we deduce that
	\begin{align*}
		T|_{U_{\frac{\rr}{4}}(\sing\sdc)}=X+Y.
	\end{align*}
	We are done.\qed
	\section{$\cpt$ cone products give persistent singular sets}\label{seccpt}
	Recall from our plan of proof that we will use two type of singular sets in our area-minimizing representative $N\#\ssc$ in Lemma \ref{nsmin}: transverse intersections from Fact \ref{fctgl} and products of $C(\cpt)$ with spheres in Fact \ref{fctcb}. Our proof of Lemma \ref{nst} in the last section deals with the Fact \ref{fctgl}. In this section we will deal with Fact \ref{fctcb}.
	
	Before stating our goal in this section, it is beneficial to give a definition of $C(\cpt)$. It is well known that the Veronese embedding of $\cpt$ is an area-minimizing cone (\cite{CJg}). However, for reader's convenience, we produce an explicit smooth calibration of $C(\cpt)$ as follows.
	
	Take the standard embedding of the Lie algebra $\sut$ into $3\times 3$ complex matrices, i.e.,
	\begin{align*}
		\sut=\{U\textnormal{ is a }3\times 3\textnormal{ matrix over }\C\space\,|\,U+U^\dagger=0,\tr U=0\}.
	\end{align*} As a real vector space $\sut$ has dimension $8$,
	\begin{align*}
		\dim_\R\sut=8.
	\end{align*}
	\begin{defn}\label{defnka}
		Define a Riemannian metric $\ka$ on $\sut$ by
		\begin{align*}
			\ka(v,w)=-\frac{\tr(vw)}{2}.
		\end{align*}
	\end{defn}
	In other words $\ka$ is $-\frac{1}{12}$ times the Killing form. For the unexperienced reader, we will define the Killing form later in (\ref{defnkf}).
	\begin{defn}\label{defnsu2}
		We will use $\su(2)$ to denote the subspace in $\sut$ generated by $\vv_2,\vv_3,\vv_4$ below 
		\begin{align*}
			\vv_2=\begin{bmatrix}
				i&0&0\\
				0&-i&0\\
				0&0&0\\
			\end{bmatrix},	\vv_3=\begin{bmatrix}
				0&1&0\\
				-1&0&0\\
				0&0&0\\
			\end{bmatrix},	\vv_4=\begin{bmatrix}
				0&i&0\\
				i&0&0\\
				0&0&0\\
			\end{bmatrix}.
		\end{align*}
	\end{defn}
	The reader can easily verify that $\vv_2,\vv_3,\vv_4$ indeed form an $\su(2)$ subalgebra of $\sut.$
	\begin{fact}\label{fctpsi}The $3$-form
		\begin{align*}
			\psi(u,v,w)=\frac{1}{2}\ka([u,v],w)
		\end{align*} is a calibration form and calibrates the adjoint action images of $\su(2)$ under $\sus.$
	\end{fact}
	The above fact is \cite[Theorem 5]{HT}. For the convenience of the reader, we will later reproduce the entire proof.
	\begin{defn}
		Define
		\begin{align*}
			\vv_1=	\frac{1}{\sqrt{3}}\begin{bmatrix}
				i&0&0\\
				0&i&0\\
				0&0&-2i\\
			\end{bmatrix}.
		\end{align*}
	\end{defn}
	Note that $\vv_1$ lies on the unit sphere in $\sut$ with respect to the metric $\ka.$
	\begin{lem}\label{fctcpt}
		The orbit of $\vv_1$ under adjoint action of $\sus$ is an embedded $\cpt$ inside the unit sphere $S^7$ in $\sut.$ The cone $C(\cpt)$ over the orbit of $\vv_1$ is calibrated by $\ast\psi.$
	\end{lem}
	Here $\ast$ means the Hodge dual. The experienced reader might identify $\ast\psi$ as the standard $5$-form in $\operatorname{PSU}(3)$ special holonomy. However, we will never use any connection to the $\operatorname{PSU}(3)$ geometry.
	\begin{defn}(\textbf{In this section only})
		For $d\ge 5,c\ge 3$, take the standard isometric embedding of $\R^8\cong \sut$ into $\R^{c+5}\cong  \R^8\times\R^{c-3}$ by mapping $x\mapsto (x,\{0\}^{c-3}).$ Then define 
		\begin{align*}
			C=&\textnormal{embedded image of }C(\cpt),\\
			\phi_C=&\pi_{\R^8}\du(\ast\psi).
		\end{align*}		
	\end{defn}
	Here $\pi_{\R^8}$ is the orthogonal projection from $\R^{c+5}\cong \R^8\times \R^{c-3}$ on to the $\R^8$ factor. Note $\pi_{\R^8}$ is $5$-area-non-increasing, and by Fact \ref{anical}, $\pi\du_{\R^8}(\ast\psi)$ is a smooth calibration form that calibrates $C.$

	Our goal in this section is to achieve the following.
	\begin{lem}\label{cpst}
		The subset $\sing\sfc$ of the area-minimizing integral current $N\#\si^{d-5}(C)$ obtained in Lemma \ref{nsmin} is a persistent singular subset (Definition \ref{psing}), provided  $d\ge 5,c\ge 3.$
	\end{lem}
	The proof is actually quite easy. Roughly speaking, recall Lemma \ref{localt}. In $U_{\ee}(\sing\sfc),$ for $(d-5)$-dimensional Hausdorff measure almost every $p\in \{0\}^{c+5}\times S^{d-5}$, the intersection
	\begin{align*}
		N\#\sfc\cap \pi_{S^{d-5}}\m(p),
	\end{align*}
	gives an integral current $T_p$ with boundary $\cpt$ and $T_p$ is smooth near its boundary. This implies that $T_p$ must has an interior singular point, since $\cpt$ cannot bound a smooth $5$-manifold. 
	
	Our plan for this section is as follows. First, we will prove Fact \ref{fctpsi} and Lemma \ref{fctcpt}. Then we will prove Lemma \ref{cpst} using the argument in the previous paragraph.
	\subsection{Proving Fact \ref{fctpsi}}
	For $v,w\in \sut,\xi\in \sus$ define the adjoint action of $\xi$ and $v$ on $w$ by
	\begin{align*}
		\ad_v(w)&=[v,w],\\
		\add_\xi(w)&=\xi w\xi\m.
	\end{align*}where $[\cdot,\cdot]$ is the Lie bracket and all products are taken as matrices.
	
	Define the Killing form on $\sut$ to be
	\begin{align}\label{defnkf}
		\fb(v,w)=\tr(\ad_v\circ\ad_w).
	\end{align}
	It is well known that $\fb$ is a negative definite symmetric blinear form on $\sut$ and equals $6$ times the matrix trace of the matrix product of $vw$, i.e.,
	\begin{align*}
		\fb(v,w)=6\tr (vw).
	\end{align*}
	For instance, this is proved in \cite[Proposition 21.24]{GQ}
	
	The associated $3$-form to $\fb$ is defined by
	\begin{align*}
		\fb([u,v],w).
	\end{align*}
	Straightforward calculations using Jacobi identity shows that $\fb([u,v],w)$ is indeed a constant coefficient $3$-form on $\sut.$ It is also easy to show that $\fb([u,v],w)$ is $\add$-invariant.
	
	Thus, the metric $\ka$ defined in Definition \ref{defnka}
	\begin{align*}
		\ka(v,w)=-\frac{1}{12}\fb(v,w)
	\end{align*}
	is indeed a positive definite metric, and the form 
	\begin{align*}
		\psi(u,v,w)=-\frac{1}{24}\fb([u,v],w)
	\end{align*}
	is a well-defined $3$-form.
	
	To prove that $\psi$ calibrates $\su(2)$ defined in Definition \ref{defnsu2}, we can use  \cite[Theorem 5]{HT}. Indeed, just above the statement of  \cite[Theorem 5]{HT} the author explicitly mentions that the calibrated $3$-planes of $\psi$ corresponds to tangent spaces of  $\operatorname{SU}(2)$ subgroups.
	However, for the readers who are not familiar with the Lie theoretic arguments therein, we provide a completely self-contained proof here, which also sets the proper preliminaries for the proof of Lemma \ref{fctcpt}.

	Let us first provide an orthonormal basis of $\sut$. 
	\begin{defn}Define $\vv_1,\dots,\vv_8$ as follows
		\begin{align*}
			&\vv_1=\frac{1}{\sqrt{3}}\begin{bmatrix}
				i&0&0\\
				0&i&0\\
				0&0&-2i\\
			\end{bmatrix},	\vv_2=\begin{bmatrix}
				i&0&0\\
				0&-i&0\\
				0&0&0\\
			\end{bmatrix},	\vv_3=\begin{bmatrix}
				0&1&0\\
				-1&0&0\\
				0&0&0\\
			\end{bmatrix},	\vv_4=\begin{bmatrix}
				0&i&0\\
				i&0&0\\
				0&0&0\\
			\end{bmatrix},\\&
			\vv_5=\begin{bmatrix}
				0&0&1\\
				0&0&0\\
				-1&0&0\\
			\end{bmatrix},
			\vv_6=\begin{bmatrix}
				0&0&i\\
				0&0&0\\
				i&0&0\\
			\end{bmatrix},
			\vv_7=\begin{bmatrix}
				0&0&0\\
				0&0&1\\
				0&-1&0\\
			\end{bmatrix},
			\vv_8=\begin{bmatrix}
				0&0&0\\
				0&0&i\\
				0&i&0\\
			\end{bmatrix}.
		\end{align*}
	\end{defn}
	Direct calculations can show that $\vv_1,\dots,\vv_8$ form an orthonormal basis of $\sut$ with respect ot $\ka$. We will also provide Mathematica verifications at then end of the paper.
	
	Next, let us calculate the Lie bracket of $\vv_1,\vv_2$ with the entire basis.
	\begin{align}\label{liemt}
		\begin{bmatrix}
			0&0&0&0&\sqrt{3}\vv_6&-\sqrt{3}\vv_5&\sqrt{3} \vv_8&-\sqrt{3}\vv_7\\
			0&0&2\vv_4&	-2\vv_3&\vv_6&-\vv_5&-\vv_8&\vv_7
		\end{bmatrix}
	\end{align}
	Here the $ij$-th entry in the matrix means $[\vv_i,\vv_j].$
	
	Now we are ready to verify that $\psi$ is a calibration form. 
	
	\begin{fact}
		Any element $v\in \sut$ can be sent to a linear combination of $\vv_1$ and $\vv_2$ by the $\add$-action of $\sus.$ 
	\end{fact}
	\begin{proof}
		To see this, note that $iv$ is a Hermitian matrix, so $iv$ is unitarily diagonalizable, i.e.,
		\begin{align*}
			iv=\zeta \Lambda\zeta^\dagger,
		\end{align*}with $\zeta\in \operatorname{U}(3)$ and $\Lambda$ a zero trace real diagonal matrix. This implies that
		\begin{align*}
			v=((\det\zeta)^{-\frac{1}{3}}\zeta)(-i\Lambda)((\det\zeta)^{-\frac{1}{3}}\zeta)^\dagger,
		\end{align*}
		where $(\det\zeta)^{-\frac{1}{3}}$ is any complex root of $\det\zeta.$  Since $((\det\zeta)^{-\frac{1}{3}}\zeta)\in \sus$, and $-i\Lambda$ is a zero trace imaginary diagonal matrix. We are done.
	\end{proof}
	Since $\add$-actions preserve the metric $\ka$, to evaluate the comass of $\psi,$ we only have to calculate $\psi(P)$ where $P$ is a $3$-dimensional vector subspace containing a vector in $\operatorname{span}\{\vv_1,\vv_2\}.$ To be precise, set
	\begin{align*}
		u=\cos\ta \vv_1+\sin\ta \vv_2,\textnormal{ }v=\sum_{3\le j\le 8}a_j\vv_j,\textnormal{ }w=\sum_{3\le j\le 8}b_j\vv_j,
	\end{align*}with $u,v,w$ an orthonormal frame of $P$.
	
	Then by Cauchy-Schwarz, we have
	\begin{align*}
		|\psi(u,v,w)|^2=\frac{1}{4}\ka([u,v],w)^2\le\frac{1}{4} \ka([u,v])^2\ka(w)^2=\frac{1}{4}\ka([u,v])^2.
	\end{align*}
	On the other hand, by (\ref{liemt}) we have
	\begin{align*}
		&\ka([u,v])^2\\
		=&\ka\bigg(\sum_{3\le j\le 8}a_j\big(\cos\ta[\vv_1,\vv_j]+\sin\ta[\vv_2,\vv_j]\big) ,\sum_{3\le j\le 8}a_j\big(\cos\ta[\vv_1,\vv_j]+\sin\ta[\vv_2,\vv_j]\big) \bigg)^2\\
		=&4\sin^2\ta (a_3^2+a_4^2)+(3\cos^2\ta+\sin^2\ta)(a_5^2+a_6^2+a_7^2+a_8^2)\\
		\le&4\sum_{j=3}^8a_j^2\\
		=&4.
	\end{align*}
	Thus, we must have
	\begin{align*}
		|\psi(u,v,w)|\le 1.
	\end{align*}
	This shows that $\psi$ has comass at most one. 
	
	However, we have
	\begin{align*}
		\psi(\vv_2,\vv_3,\vv_4)=\frac{1}{2}\kappa([\vv_2,\vv_3],\vv_4)=\frac{1}{2}\kappa(2\vv_4,\vv_4)=1.
	\end{align*}
	Thus the $\su(2)$ subalgebra generated by $\vv_2,\vv_3,\vv_4$ is indeed calibrated by $\psi$.
	\subsection{Proof of Lemma \ref{fctcpt}}
	To determine the orbit of $\vv_1$ under $\add$-action, let us first calculate the stabilizer subgroup of $\vv_1$ in $\add$-action.
	
	Let $\zeta$ be an arbitrary $3\times 3$ matrix. We can write $\zeta$ as a block matrix with $2\times 2,2\times 1,1\times 2,1\times 1$ blocks $E,F,G,H,$  i.e.,
	\begin{align*}
		\zeta=\begin{bmatrix}
			E&F\\G&H
		\end{bmatrix}.
	\end{align*}
	If $\zeta \vv_1\zeta\m=\vv_1,$ then we have
	\begin{align*}
		0=&[\zeta,\vv_1]=\frac{1}{\sqrt{3}}\begin{bmatrix}
			E&F\\G&H
		\end{bmatrix}\begin{bmatrix}
			iI_2&0\\0&-2i
		\end{bmatrix}-\frac{1}{\sqrt{3}}\begin{bmatrix}
			iI_2&0\\0&-2i
		\end{bmatrix}\begin{bmatrix}
			E&F\\G&H
		\end{bmatrix}\\
		=&\begin{bmatrix}
			0&-\sqrt{3}iF\\\sqrt{3}iG&0
		\end{bmatrix}.
	\end{align*}
	This implies that $F=0,G=0,$ i.e., $\zeta$ is block diagonal with $2\times 2,1\times 1$ blocks. 
	
	Thus, any stabilizer $\zeta$ of $\vv_1$ under $\add$-action of $\sus$ consists of special unitary matrices that are $2\times 2,1\times 1$ block diagonal of the form
	\begin{align*}
		\zeta=	\begin{bmatrix}
			E&0\\0&H
		\end{bmatrix}.
	\end{align*}
	Since $\zeta\zeta^\dagger=\id,$ we deduce that $E\in \operatorname{U}(2)$ and $F\in \operatorname{U}(1).$ Since $\det \zeta=1,$ we deduce that the stabilizer of $\vv_1$ is a subgroup $$H=\operatorname{S}(\operatorname{U}(2)\times \operatorname{U}(1)).$$
	
	However, it is well known that $\cpt$ as a homogenous space is isometric to $
	\sus/H,$ e.g., \cite[p.452 last paragraph]{SH}.
	
	Thus the orbit of $\vv_1$ under $\add$-action of $\sus$ is indeed an embedded $\cpt$ in the unit sphere of $\sut.$ We will abuse the notations and also use $\cpt$ to denote this orbit of $\vv_1$.
	
	Since both $C(\cpt)$ and $\ast\psi$ are $\add$-invariant, to prove $C(\cpt)$ is calibrated by $\ast\psi$, by by bullet (\ref{cms3}) of Fact \ref{cmsvec}, it suffices to show that the tangent space to $C(\cpt)$ at $\vv_1$ is orthogonal to the $\su(2)$ generated by $\vv_2,\vv_3,\vv_4.$
	
	The tangent space to $C(\cpt)$ at $\vv_1$ are generated by $\vv_1$ and vectors of the form $[w,\vv_1]$ with $w\in\sut$, by differentiating $\add$-actions. However, checking (\ref{liemt}), it is clear to see that both $\vv_1$ and $[w,\vv_1]$ for any $w\in\sut$ are always orthogonal to $\vv_2,\vv_3,\vv_4.$ We are done.
	\subsection{Proving Lemma \ref{cpst}}
	Recall Lemma \ref{localt}. Take an area-minimizing integral current $T\in [\Si]$ with respect to a metric $g\in\Om_h.$  Restricted to $U_{2\ee}(\sing\sfc),$ consider the projection
$
		\pi_{S^{d-5}},
$ defined by the orthogonal projection onto the $S^{d-5}$ factor in $B_{2\ee}^{c+5}\times S^{d-5}$. Then by \cite[4.3.6,4.3.13]{HF}, for $(d-5)$ dimensional Hausdorff measure almost every $p\in S^{d-5}$, the slicing of $T|_{U_{2\ee}(\sing\sfc)}$ by $\pi_{S^{d-5}}$,
	\begin{align*}
		T_p=\ri{T,\pi_{S^{d-5}},p}
	\end{align*} is a $5$-dimensional integral current. Intuitively, the reader can just understand this as saying that the intersection of $T|_{U_{2\ee}(\sing\sfc)}$ with $\pi_{S^{d-5}}\m(p)$ is a $5$-dimensional integral current 
	\begin{align*}
		T_p,
	\end{align*} 
	for $(d-5)$-dimensional Hausdorff measure almost every $p\in S^{d-5}.$

	Suppose for some $p\in S^{d-5}$, every point $q$ with $$q\in \pi_{S^{d-5}}\m(p)\cap\supp T|_{U_{2\ee}(\sing\sfc)}$$ belongs to the regular set of $T$. Then since the regular set of $T$ is relatively open, we deduce that there is an open neighborhood $U_p$ of $p$ inside $S^{d-5}$ such that 
	\begin{align*}
	\pi_{S^{d-5}}\m(U_p)\cap\supp T|_{U_{2\ee}(\sing\sfc)}
	\end{align*}belongs to the regular set of $T.$ 	 By Sard's Theorem (\cite[Theorem 6.10]{JL}) this implies that for $(d-5)$-dimensional Hausdorff measure almost every $q\in U_p$, the current $T_q$ is a smooth submanifold with boundary. However, this is impossible.
	
To see this, by Lemma \ref{localt}, $T$ is smooth in $$\bigg(U(N\cap N)\cup U_{\ee}(\sing\sfc)\bigg)\cp.$$ By Sard's Theorem (\cite[Theorem 6.10]{JL}) this implies that $T_p$ is smooth in $$U_{\ee}(\sing\sfc)\cp$$ for $(d-5)$-dimensional Hausdorff measure almost every $p\in S^{d-5},$ and such $T_p$ necessarily has boundary $\pd T_p$ diffeomorphic to $\cpt.$ As $\cpt$ cannot be the boundary of a $5$-manifold \cite[page 203]{MScs}, this implies that for $(d-5)$-dimensional Hausdorff measure almost every  $p\in S^{d-5}$, the current $T_p$ must have a singular set in its interior. We have arrived at a contradiction.

In other words we have proved that $\pi_{S^{d-5}}\m(p)\cap \sing T$ has  zero dimensional Hausdorff measure at least $1$ for $(d-5)$-dimensional Hausdorff measure almost every $p\in S^{d-5}$. Now integrate this estimate using Eilenberg's inequality \cite[2.10.25]{HF}, we deduce that $\sing T$ has positive $(d-5)$-dimensional Hausdorff measure. We are done.
	\section{Conclusions and remarks}\label{conc}
	Our Theorem \ref{mtc} clearly follows from Lemma \ref{nst} and Lemma \ref{cpst}.
	\subsection{About the non-orientable case}
	Our main theorem Theorem \ref{mtc} holds if we drop the orientbility hypothesis of $M$ and replace integral homology with mod $2$ homology and integral currents with mod $2$ currents. The proof relies on replacing the gluing calibrations in Section \ref{zhang} by gluing area-non-increasing retractions. We will pursue the generalizations of Theorem \ref{mtc} to the non-orientable case in another manuscript.
	\subsection{About singular sets of area-minimizing currents coming from transverse intersections}
	Using Remark \ref{remim}, our proof of Lemma \ref{nst} can be easily used to imply the following fact:
	\begin{fact}\label{fctea}
		Let $T$ be the unique $d$-dimensional area-minimizing integral current in $[\Si]$ on $M^{d+c}$ with respect to a metric $h.$ If $T$ restricted to a smooth open set $U$ equals as an integral current the sum of two smooth minimal submanifolds $X,Y$  intersecting transversely along their interior, then there exists an open subset $\Om_{h}$ containing $h$ in the space of Riemannian metrics, such that any area-minimizing integral current in $[\Si]$ with respect to any metric $g\in\Om_h$  decomposes into the sum of two transversely intersecting minimal submanifolds when restricted to $U,$ provided
		\begin{align*}
			\frac{\pi}{2}\ge \Ta(X,Y)>2\arctan t_c=\frac{4}{c}+O(c^{-2}),\textnormal{ and }d\ge c\ge 3.
		\end{align*}
	\end{fact}
	As mentioned in the previous subsection the condition of the orientability can be dropped if we replace integral currents with mod $2$ currents. Moreover, the sum decomposition $T|_U=X+Y$ can be relaxed to the condition that $T|_U$ is an immersion with clean intersections, without a sum decomposition or transversality. We will purse it in another manuscript. Here a point in the intersection set is called clean if there is a coordinate chart centered at that point such that the immersion decomposes into two coordinate axis subspaces intersecting along coordinate axis subspaces in the chart.
	
	A hidden point in our argument is that the self-intersection set $N\cap N$ of area-minimizing representative already has some points being transverse intersection, although in general, triple, quadruple intersections, etc., do appear in $N\cap N.$ In case $N\cap N$ is non-empty, why do we still have to add $\sing\sdc$ manually to our area-minimizing representative?
	
	The reason is that Lawlor's ideas do not quite work for intersections with multiplicities more than $2.$ At the least the author does not know how to construct the area-non-increasing retractions onto each sheet in the intersection.
	\subsection{About singular sets of area-minimizing currents coming from $C(\cpt)$}
	Similarly our proof of Lemma \ref{cpst} shows that
	\begin{fact}\label{fcteb}
		Let $T$ be the unique $d$-dimensional area-minimizing integral current in $[\Si]$ on $M^{d+c}$ with respect to a metric $h.$ If $T$ restricted to a smooth open set equals an $(d-s)$-dimensional-cone $C(V)$ bundle over a base submanifold of dimension $s,$ then there exists an open subset $\Om_{h}$ containing $h$ in the space of Riemannian metrics, such that any area-minimizing integral current in $[\Si]$ with respect to any metric $g\in\Om_h$ has a singular set of positive Hausdorff $s$-dimensional measure, provided that $V$ belongs to a nontrivial class in the bordism ring.\end{fact}
	Here a $(d-s)$-dimensional-cone $C(V)$ bundle over a base submanifold means a fiber bundle over the base submanifold with each fiber a $C(V).$ Such area-minimizing cones $C(V)$ with non-bounding links $V$ are abundant.
	\subsection{About other spaces of metrics}
	Our Theorem \ref{mtc} is proved for the space of smooth Riemannian metrics. Inspecting the proof, it is easy to see that we have only used the continuity of derivatives of the Riemannian metric up to degree $2$. Thus, our Theorem \ref{mtc} can be proved as well when we replace the space of smooth metrics in Theorem \ref{mtc} with the space of $C^k$ Riemannian metrics for any $k\ge 2.$
	
	At the other end of differentiability of metrics, if we consider the space of real analytic metrics, then Theorem \ref{mtc} also holds. To see this, recall the classical fact that real-analytic metrics are dense in the space of smooth metrics. Thus, a non-empty open subset of the space of smooth metrics necessarily gives a non-empty open set in the space of real analytic metrics.
	
	There are also some conjectures about the analogues of Theorem \ref{mtc} in the space of specially calibrated integral currents in metrics of special holonomy, for example, special Lagrangians on Calabi-Yau manifolds. Our method per se does not yield a way to construct such metrics with special holonomy. However, Fact 	\ref{fctea} and Fact \ref{fcteb} still apply and it is worth noting both types of singularities in these facts do appear in special Lagrangian tangent cones.
	\printbibliography
	\includepdf[pages=-]{verifications.pdf}
\end{document}